\documentclass[11pt]{article}
\usepackage{latexsym}
\usepackage{amssymb}
\usepackage{amsfonts}

\input{epsf}

\def\proof{{\sc Proof. }}
\def\qed{\hfill$\Box$}

\newcommand{\RL}{{\mathbb R}}

\newcommand{\IN}{{\mathbb Z}}

\newcommand{\IND}{{\mathbb I}}

\newcommand{\VAR}{\mbox{\rm Var}}

\newcommand{\LA}{\Lambda}

\newcommand{\la}{\lambda}

\newcommand{\al}{\alpha}

\def\be{\begin{eqnarray}}
\def\ee{\end{eqnarray}}
\def\ben{\begin{eqnarray*}}
\def\een{\end{eqnarray*}}

\def\barlambda{\oo\lambda}

\def\LV{L_\infty^V}

\def\haxi{\hat\xi}
\def\barb{\oo b}

\def\barxi{\oo\xi}

\def\barF{\oo{F}}

\def\ddt{\frac{d}{dt}}

\def\Ebox#1#2{%
\begin{center}    
 \parbox{#1\hsize}{\epsfxsize=\hsize \epsfbox{#2}}
\end{center}}

\def\tlabel#1{\label{t:#1}}
\def\slabel#1{\label{s:#1}}
\def\flabel#1{\label{f:#1}}
\def\elabel#1{\label{e:#1}}

\def\sq{\hbox{\rlap{$\sqcap$}$\sqcup$}}
\def\qed{\ifmmode\sq\else{\unskip\nobreak\hfil
\penalty50\hskip1em\null\nobreak\hfil\sq
\parfillskip=0pt\finalhyphendemerits=0\endgraf}\fi}

\newsavebox{\junk}
\savebox{\junk}[1.6mm]{\hbox{$|\!|\!|$}}
\def\lll{{\usebox{\junk}}}

\def\limsup{\mathop{\rm lim\ sup}}
\def\liminf{\mathop{\rm lim\ inf}}

\def\state{{\sf X}}

\newcommand{\field}[1]{\mathbb{#1}}

\def\Re{\field{R}} 
\def\TT{\field{T}}

\def\nat{\field{Z}_+}

\def\Co{\field{C}}

\def\ind{\field{I}}

\def\One{\hbox{\large\bf 1}}

\def\cpi{\check{\pi}}

\def\cP{{\check{P}}}

\newcommand{\cf}{\check f}
\newcommand{\cg}{\check g}

\newcommand{\cF}{\check F}

\newcommand{\cclA}{{\check{\cal A}}}

\newcommand{\cmu}{{\check{\mu}}}

\newcommand{\cV}{{\check{V}}}

\newcommand{\haF}{{\widehat{F}}} 
\newcommand{\haM}{{\widehat{M}}}
\newcommand{\haQ}{{\widehat{Q}}}

\newcommand{\haclA}{{\widehat{\cal A}}}

\newcommand{\haR}{{\widehat{R}}}
\newcommand{\haG}{{\widehat{G}}}
\newcommand{\haU}{{\widehat{U}}}

\newcommand{\limn}{\lim_{n \rightarrow \infty}}

\newcommand{\smalloneOvern}{{\textstyle\frac{1}{n}\:}}

\newcommand{\smalloneOvert}{{\textstyle\frac{1}{t}\:}}

\newcommand{\one}{\hbox{\rm\large\textbf{1}}}

\newcommand{\haV}{{\widehat{V}}}

\def\Ree{\hbox{\rm Re\,}}

\def\bfmm{{\mbox{\protect\boldmath$m$}}}

\def\bfmW{{\mbox{\protect\boldmath$W$}}}

\def\bfPhi{\mbox{\protect\boldmath$\Phi$}}

\def\bfpsi{\mbox{\boldmath$\psi$}}
\def\bfphi{\mbox{\boldmath$\phi$}}

\def\haP{{\widehat P}}

\def\hagamma{{\hat\gamma}}

\def\til={{\widetilde =}}

\def\clA{{\cal A}}
\def\clB{{\cal B}}

\def\clD{{\cal D}}

\def\clF{{\cal F}}
\def\clG{{\cal G}}
\def\clH{{\cal H}}
\def\clJ{{\cal J}}

\def\clL{{\cal L}}
\def\clM{{\cal M}}

\def\clO{{\cal O}}

\def\clQ{{\cal Q}}

\def\clS{{\cal S}}

\def\clX{{\cal X}}

\def\half{{\mathchoice{\textstyle \frac{1}{2}}%
{\frac{1}{2}}%
{\hbox{\tiny $\frac{1}{2}$}}%
{\hbox{\tiny $\frac{1}{2}$}} }}

\def\onefourth{{\mathchoice{\textstyle {1\over 4}}%
{\textstyle{1\over 4}}%
{\hbox{\tiny $1\over 4$}}%
{\hbox{\tiny $1\over 4$}} }}

\def\eqdef{\mathbin{:=}}

\def\Prob{{\sf P}}
\def\Probsub{{\sf P\! }}
\def\Expect{{\sf E}}

 \def\eq#1/{(\ref{#1})}

\def\epsy{\varepsilon}

\def\varble{\,\cdot\,}

\newtheorem{theorem}{Theorem}[section]
\newtheorem{corollary}[theorem]{Corollary}
\newtheorem{proposition}[theorem]{Proposition}
\newtheorem{lemma}[theorem]{Lemma}

\def\Lemma#1{Lemma~\ref{t:#1}}
\def\Proposition#1{Proposition~\ref{t:#1}}
\def\Theorem#1{Theorem~\ref{t:#1}}
\def\Corollary#1{Corollary~\ref{t:#1}}

\def\Section#1{Section~\ref{s:#1}}

\def\Figure#1{Figure~\ref{f:#1}}

\def\eq#1/{(\ref{e:#1})}

\newcommand{\beqn}[1]{\notes{#1}%
\begin{eqnarray} \elabel{#1}}

\newcommand{\eeqn}{\end{eqnarray} }

\newcommand{\beq}[1]{\notes{#1}%
\begin{equation}\elabel{#1}}

\newcommand{\eeq}{\end{equation}} 

\def\bdes{\begin{description}}
\def\edes{\end{description}}

\newcommand{\oo}{\overline}
 
\def\baromega{\oo\omega}

\def\bara{{\oo {a}}}

\def\barB{{\oo {B}}}

\def\barG{{\oo {G}}}

\def\barbeta{{\oo{\beta}}}

\def\proof{\paragraph{\sc Proof. }} 
\def\proofo{\paragraph{\sc Proof Outline. }} 
 
\def\notes#1{}

\def\LV{L_\infty^V}

                \def\cP{P}        %  NOTE no check on P_\alpha, V_\alpha etc
                \def\cV{V} 
                \def\cpi{\pi} 
                \def\cclA{\clA} 
                
\def\bfmw{{\mbox{\protect\boldmath$w$}}}
\def\bfgamma{{\mbox{\protect\boldmath$\gamma$}}}

\begin{document}
\bibliographystyle{plain}

\title{Spectral Theory and Limit Theorems\\
for Geometrically Ergodic Markov Processes}

\author
{
        I. Kontoyiannis\thanks{Division
                of Applied Mathematics and Department
		of Computer Science, Brown University,
                Box F, 182 George St., Providence, RI 02912, USA.
                Email: {\tt yiannis@dam.brown.edu}
                Web: {\tt www.dam.brown.edu/people/yiannis/}.
                Work supported in part
                by NSF grants \#0073378-CCR 
		and DMS-9615444.
                              }
\and
        S.P. Meyn\thanks{Department of Electrical and Computer 
                Engineering and the Coordinated Sciences Laboratory, 
                University of Illinois at Urbana-Champaign,  
                Urbana, IL 61801, U.S.A. Email: {\tt s-meyn@uiuc.edu}.
                Part of the research for this paper was done while S.M.
                was a Fulbright research scholar and visiting professor 
                at the Indian Institute of Science, and a visiting
                professor at the Technion. Work supported in part by
                NSF grants ECS 940372, ECS 9972957.
                        }
}

\maketitle

\thispagestyle{empty}
\setcounter{page}{0}

\centerline{\textbf{Short Title:} Geometrically Ergodic Markov Processes}
\bigskip

\begin{abstract}

\noindent
Consider the partial sums $\{S_t\}$ of a real-valued 
functional $F(\Phi(t))$ of a Markov chain 
$\{\Phi(t)\}$ with values in 
a general state space. Assuming only that the Markov 
chain is geometrically ergodic and that the functional 
$F$ is bounded, the following conclusions are obtained:

\begin{description}

\item \textit{Spectral theory}: 
Well-behaved solutions $\cf$ 
can be constructed for the 
``multiplicative Poisson equation'' $(e^{\alpha F}P)\cf=\la\cf$,
where $P$ is the transition kernel of the 
Markov chain, and $\alpha\in\Co$ is a constant.
The function $\cf$ is an eigenfunction, with 
corresponding eigenvalue $\la$,
for the kernel $(e^{\alpha F}P)=e^{\alpha F(x)}P(x,dy)$.

\item \textit{A ``multiplicative'' mean ergodic theorem}:
For all complex $\alpha$ in a neighborhood of 
the origin, the normalized mean of 
$\exp(\alpha S_t)$ (and not the logarithm of 
the mean) converges to $\cf$
exponentially fast, where $\cf$ 
is a solution of the multiplicative 
Poisson equation.

\item \textit{Edgeworth Expansions}: Rates are 
obtained for the convergence of the distribution 
function of the normalized partial sums $S_t$
to the standard Gaussian distribution. The 
first term in this expansion is of order 
$(1/\sqrt{t})$, and it depends on the initial 
condition of the Markov chain through the 
solution $\haF$ of the associated Poisson 
equation (and not the solution $\cf$ of
the multiplicative Poisson equation).
      
\item \textit{Large Deviations}:     
The partial sums are shown to 
satisfy a large deviations principle 
in a neighborhood of the mean.
This result, proved under geometric ergodicity 
alone, cannot in general be extended to the
whole real line.

\item \textit{Exact Large Deviations Asymptotics}:
Rates of convergence are obtained for the large deviations
estimates above.
The polynomial pre-exponent is of order 
$(1/\sqrt{t})$, and its coefficient depends 
on the initial condition of the Markov chain
through the solution $\cf$
of the multiplicative Poisson equation. 
\end{description}

\noindent
Extensions of these results to continuous-time Markov 
processes are also given.



\end{abstract}

%%%%%%%%%%%%%%%%%%%%%%%%%%%%%%%%%%%%%%%%%%%%%%%%%%%%%%%%%
% \notes{Note: notes defs here}

\newlength{\noteWidth}
\setlength{\noteWidth}{.6in}
\long\def\notes#1{\ifinner
             {\tiny #1}
             \else
              \marginpar{\parbox[t]{\noteWidth}{\raggedright\tiny #1}}
               \fi}

% \def\notes#1{\smallbreak \begin{quote}
% {\it\noindent\textbf{\Large !} #1 } \end{quote} }
% \def\notes#1{}     %%%%%%%%%% For final version
%%%%%%%%%%%%%%%%%%%%%%%%%%%%%%%%%%%%%%%%%%%%%%%%%%%%%%%%%

\newpage
\section{Introduction}
\slabel{intro}
%%%%%%%%%%%%%%%%%%%%%%%%%%% INTRODUCTION %%%%%%%%%%%%%%%%%%%%%%%%%%% 

Consider a Markov process $\bfPhi=\{\Phi(t)\,:\, t \in \TT\}$
taking values in a general state space $\state$, 
and with time being either continuous, $\TT=[0,\infty)$,
or discrete, $\TT=\{0,1,\ldots\}$.
Let $F:\state\to\Re$ be a given functional
on the state space of $\bfPhi$. 

Our interest 
lies in the long-term behavior of 
\begin{equation}
S_t
    =   \int_{[0,t)} F (\Phi(s)) \, ds 
        ,    \qquad t\in\TT,
\elabel{psums}
\end{equation}
where in discrete-time the integral
is a sum, and $S_t$ are simply the
partial sums 
\begin{equation}
S_n
    =   \sum_{i=0}^{n-1} F (\Phi(i)) ,    \qquad n\geq 1.
\elabel{discretePS}
\end{equation}

\subsection{Multiplicative Ergodic Theory}
%%%%%%%%%%%%%%%%%%%%%%% MULT ERGODIC THEORY %%%%%%%%%%%%%%%%%%%%%%%
\slabel{ergodic}

For simplicity we first discuss 
the case of a discrete-time Markov 
chain $\bfPhi$ with a countable 
state space $\state$. If $\bfPhi$ is 
positive recurrent with invariant 
probability measure $\pi$, then for 
any $F$ with finite mean 
$\pi(F)=\sum_x\pi(x)F(x)$,
\begin{eqnarray}
\elabel{meanLLN}
\smalloneOvern \!
\Expect_x [
S_n ]
&\ \to & \pi (F), 
\quad n\to\infty,
\end{eqnarray}
where $x=\Phi(0)$ is the initial condition,
$S_n$ are the partial sums defined
above, $\Probsub_x$ is the law of
$\bfPhi$ conditional on $\Phi(0)=x,$ and
$\Expect_x$ is the corresponding expectation.

% Moreover, for any fixed 
% $\theta \in \state$, the steady-state 
% mean of $F$ can be written as
% \begin{equation}
% \pi (F) = \inf \bigl\{\eta: 
% \Expect_\theta \bigl[ S_{\tau_\theta}  
% - \eta\tau_\theta  \bigr] \le 0 \bigr\}
% \,,
% \elabel{etaDef}
% \end{equation} 
% where $\tau_\theta$ is the first time
% $\bfPhi$ visits state $\theta$.

Often we can quantify the rate of 
convergence in \eq meanLLN/ by 
showing that 
the following limit exists,
\begin{equation}
\haF(x) = \limn \Expect_x \bigl[ S_n  - n \pi (F)\bigr] \, ,
\elabel{meanLLNb}
\end{equation} 
where, in fact, the function 
$\haF$ solves the {\em Poisson equation}:
\begin{equation}
P\haF = \haF - F + \pi (F).
\elabel{fish}
\end{equation}
Here $P$ denotes the transition
kernel of $\bfPhi$,
$P(x,y)\eqdef\Pr\{\Phi(1)=y\,|\,\Phi(0)=x\}$,
and $P$ acts on functions $f:\state\to\RL$
via $Pf(x)=\sum_y P(x,y)f(y)$.
Results of this kind hold for a wide
class of Markov chains on a general 
state space, as shown in 
\cite{meyn-tweedie:book} 
in discrete-time and in 
\cite{meyn-tweedie:93e,meyn-tweedie:94b} in
continuous-time.

%As we outline in the following section these
%are motivated, in part, by their applications 
%in proving precise expansions for the central limit
%theorem and for the large deviations principle
%satisfied by the partial sums $S_n$.   SAID ON NEXT PAGE

In this paper we seek multiplicative versions 
of the ergodic results in \eq meanLLN/--\eq fish/.
Let $\alpha\in\Co$, and consider the product
\[
\prod_{i=0}^{n-1} \exp (\alpha F(\Phi(i)))
=
\exp(\alpha S_n).
\]
For countable state space chains in discrete-time,
multiplicative results corresponding to the ergodic 
theorems \eq meanLLN/--\eq fish/
were established in \cite{balaji-meyn}
when $\alpha$ is a real number.
The mean ergodic theorem \eq meanLLN/
corresponds to the multiplicative limit
\begin{equation}
\smalloneOvern\! \log \Expect_x
        [ \exp ( \alpha S_n  )]
        \to \LA(\alpha),
                   \qquad n\to\infty,
\elabel{expLLN}
\end{equation}
for some analytic function 
$\LA(\alpha)\in\Re$,
and the stronger limit theorem
\eq meanLLNb/ has the multiplicative
counterpart
\begin{equation}
\cf_\alpha(x) = \limn \Expect_x
                    [\exp ( \alpha S_n  -  n\LA(\alpha) ) ],
\elabel{expLLNb}
\end{equation}
where $\cf_\alpha$
solves the
natural analog of \eq fish/,
the {\em multiplicative Poisson equation}:
\begin{equation}
P \cf_\alpha = \exp \Bigl( - \alpha F + \LA(\alpha) \Bigr)
\cf_\alpha.
\elabel{expFish}
\end{equation} 

Our first aim is to
provide natural conditions under 
which the multiplicative ergodic
results \eq expLLN/--\eq expFish/
hold. As we indicate in 
several instances, our conditions
(and the results obtained
under them) are often optimal or
near-optimal (see \Proposition{stationaryLDP},
and the examples in \Section{examples}).
Equipped with these results,
we go on to prove precise expansions
for some classical probabilistic limit 
theorems satisfied by the partial
sums $S_n$. Specifically, the 
multiplicative mean ergodic theorem 
\eq expLLNb/ leads to Edgeworth 
expansions for the central limit 
theorem and to exact large 
deviations asymptotics.

There are numerous approaches to 
multiplicative ergodic theory and 
its related spectral theory in 
the literature; a brief survey 
is given at the end of this 
introduction. The conditions
given
in this paper considerably extend 
known criteria for the existence 
of solutions to the multiplicative 
Poisson equation and for the validity 
of the multiplicative mean ergodic theorem. 

%This is achieved, in part, by placing
%within a single framework results 
%from two different areas, the theory 
%of positive operators developed in 
%\cite{nummelin:book}, and the theory 
%of positive harmonic functions 
%for diffusions \cite{pinsky:book}.  SAID BELOW

Most closely related to the approach 
taken here are the results 
of \cite{balaji-meyn}, developed for
discrete-time Markov chains $\bfPhi$ 
on a discrete state space
along the following lines.
For any real $\alpha$, define
the new kernel $\haP_\alpha$ by
\begin{equation}
\haP_\alpha (x, y) = \exp  ( \alpha F(x) ) P(x,y), \qquad
x,y \in \state,
\elabel{hatP}
\end{equation}
where $P(x,y)$ is the transition
kernel of the Markov chain $\bfPhi$.
In this notation, 
the multiplicative Poisson equation
\eq expFish/ can be rewritten as
\begin{equation}
\haP_\alpha \cf_\alpha = \lambda_\alpha \cf_\alpha,
\elabel{expFishb}
\end{equation}
with $\lambda_\alpha = \exp (\LA(\alpha))$.
That is, the solutions $\cf_\alpha$
of the  multiplicative Poisson equation
\eq expFish/ are eigenfunctions for the
new kernel $\haP_\alpha$, with associated
eigenvalues $\lambda_\alpha$.
[Throughout the paper, we try to maintain 
the convention that lower-case letters
denote quantities that are exponential
versions of the corresponding upper-case
letters; e.g., $\lambda = \exp (\LA)$].

Under a monotonicity assumption on 
$F$, it is shown in 
\cite{balaji-meyn} that well-behaved 
eigenfunctions for \eq expFishb/ 
exist for real 
$\alpha$ in a neighborhood of zero.
Based on such an eigenfunction 
$\cf_\alpha$ with corresponding 
eigenvalue $\lambda_\alpha$, the
\textit{twisted kernel} $\cP_\alpha$ is defined,
\begin{equation}
\cP_\alpha (x, y) = \lambda_\alpha^{-1} \cf_\alpha^{-1}
(x) \haP_\alpha (x, y) \cf_\alpha (y)\,,
\elabel{cP}
\end{equation}
and the convergence in
\eq expLLNb/ is deduced
from the properties of $\cP_\alpha$.
 
For bounded functionals $F$, and assuming only that
$\bfPhi$ is ``geometrically ergodic,'' 
results corresponding to 
\eq expLLN/--\eq expFish/
are obtained in \Section{spectral} 
of the present paper, 
for Markov processes 
$\bfPhi$ on a general state space, 
in continuous- or discrete-time,
and for complex $\alpha$.
For our purposes,
a Markov chain $\bfPhi$ is 
{\em geometrically ergodic} 
if it is $\psi$-irreducible, 
aperiodic, and 
a Lyapunov function $V:\state\to[1,\infty]$ 
exists such that the following
condition holds:
$$
\left. \begin{array}{ll}
   & \mbox{For a ``small'' set $C\subset\state$, and constants
                $\delta>0,\,$ $b<\infty$:}\\
   & \\
   & \hspace{0.8in} PV\leq (1-\delta)V+b\IND_C\, .
 \end{array} \right\}
\hspace{0.7in}
\mbox{\bf (V4)}$$
Precise definitions and a more general 
version of condition (V4) for Markov 
processes in discrete- or continuous-time 
are given in \Section{ergodictheorems}.

Geometric ergodicity for $\bfPhi$
is our main assumption, and it 
will remain in effect throughout 
the paper. 
\Section{related}
offers a discussion comparing 
(V4) to several of the standard 
assumptions in the relevant 
literature, and in
\Section{examples}
geometric ergodicity 
is verified for several 
classes of important examples.
Note also that what we call geometric
ergodicity here is equivalent to the
notion of geometric
ergodicity used in \cite{meyn-tweedie:book},
where it is stated slightly differently.

%%IK1:
In the following section we briefly
describe the probabilistic implications 
of the spectral theory outlined above.
Along a different direction,
in \cite{ODEmethod:02} we extend
our present results to the case 
of products of random matrices.
This extension leads to an 
interesting and non-trivial
application of the present ideas
to a stability question arising 
from systems theory.
 
\subsection{Probabilistic Limit Theorems}
%%%%%%%%%%%%%%%%%%%%%%%%%%% LIMIT THEOREMS %%%%%%%%%%%%%%%%%%%%%%%%%%% 
\slabel{intro-limthms}

The multiplicative mean ergodic theorems 
in \eq expLLN/ and \eq expLLNb/ offer
precise information about the 
asymptotic behavior, as 
$n\to\infty$, of 
$$m_n(\alpha):=\Expect_x [ \exp ( \alpha S_n  )],
\quad \alpha\in\Co\,.$$
When $\alpha=i\omega$ is imaginary,
$m_n(\alpha)$ is simply the 
characteristic function 
of the partial sums $S_n$,
and it is well-known that 
information about the convergence 
of the characteristic functions 
leads to Edgeworth expansions 
related to the central limit 
theorem \cite{fellerII:book,hall-book:82,petrov-book:95}.
Similarly, when $\alpha$ is real,
$m_n(\alpha)$ is the 
moment-generating function
of the partial sums $S_n$, 
and the precise convergence 
of the corresponding log-moment 
generating 
functions to a smooth limiting 
$\LA(\alpha)$ as in \eq expLLN/ 
leads to exact large deviations 
asymptotics; see
\cite{dembo-zeitouni:book,
chaganty-sethuraman:93}.

Suppose $\bfPhi$ is a geometrically 
ergodic Markov chain, and let $F$
be a bounded, non-lattice,
real-valued functional on the
state space of $\bfPhi$.
In \Section{edgeworth}, we 
obtain an Edgeworth expansion
for the distribution function 
$G_n(y)$ of the normalized partial 
sums $[S_n-n\pi(F)]/\sigma\sqrt{n}$,
$$G_n(y)=\Probsub_x\left\{
        \frac{S_n-n\pi(F)}{\sigma\sqrt{n}}\leq y
        \right\},
        \quad
        y\in\RL,$$
where $\sigma^2$ is the asymptotic
variance of $S_n/\sqrt{n}$.
In \Theorem{EdgeworthNL}
we show that,
for all $x\in\state$,
\ben
G_n(y)
={\cal G}(y)
+
\frac{\gamma(y)}{\sigma\sqrt{n}}
\left[
        \frac{\rho_3}{6\sigma^2}(1-y^2)
        \,-\,\haF(x)
\right]
+o(n^{-1/2}),
\quad
n\to\infty,
\een
uniformly in $y\in\RL,$
where $\gamma(y)$ denotes the standard Normal density,
${\cal G}(y)$ is the corresponding distribution
function, $\haF$ is the solution to the Poisson
equation \eq fish/ given in \eq meanLLNb/,
and $\rho_3$ is a constant
related to the third moment of $S_n/\sqrt{n}$.

A similar expansion is obtained in the case
of lattice functionals $F$. These results
generalize the Edgeworth expansions in 
\cite{nagaev:61,jensen:91,mccormick-data:93},
where they are derived under much more restrictive
assumptions. In particular, in all these papers 
the conditions given are stronger than
Doeblin recurrence, which is significantly
stronger than the form of geometric ergodicity
assumed in this paper -- see the discussions
in \Section{related} and \Section{examples}.

In \Section{ldp} we discuss moderate and large
deviations for the partial sums $S_n$. Under
geometric ergodicity, the multiplicative mean
ergodic theorem \eq expLLNb/ 
implies that a moderate deviations principle
(MDP) holds for the partial sums $S_n$. 
Note that geometric ergodicity is essentially 
equivalent to the weakest conditions known to 
suffice for the MDP 
\cite{aco97a,acoche98a} 
(although weaker assumptions can be used 
to obtain the MDP lower bound).

By standard large deviations techniques
\cite{dembo-zeitouni:book}, the convergence 
of the log-moment generating functions in 
\eq expLLN/ to a smooth limiting 
$\LA(\alpha)$ can be used to prove
large deviations estimates for the
partial sums $S_n$:
Suppose $\bfPhi$ is a Doeblin
chain, and let $F$
be a bounded, real-valued functional 
on the state space of $\bfPhi$.
In \Proposition{stationaryLDP} 
we show that under the stationary
distribution $\pi$ of $\bfPhi$,
the partial sums 
$S_n$ satisfy a large deviations
principle (LDP) in a neighborhood
of the mean $\pi(F)$, i.e., for
any $c>\pi(F)$ close enough
to the mean $\pi(F)$, 
\be
\smalloneOvern\!
\log \Probsub_\pi\{S_n\geq nc\}
       \;\to\;
        -\LA^*(c),
        \quad
        n\to\infty,
\label{eq:introLDP}
\ee
where $\LA^*(c)$ is the Fenchel-Legendre 
transform of $\LA(\cdot)$. (A
corresponding result holds for the
lower tail.) 

Note that this result {\em cannot}
in general be extended to a full LDP 
on the whole real line. For example,
Bryc and Dembo \cite{bryc-dembo:96} 
have shown that the full LDP may even 
fail for the partial sums of a Doeblin
chain with a countable state space.

Further, the more precise convergence result 
\eq expLLNb/ leads to exact large deviations 
expansions analogous to those obtained by 
Bahadur and Rao \cite{bahadur-rao:60} for
independent random variables:
For geometrically ergodic chains and
non-lattice functionals $F$, 
in \Theorem{Bahadur-RaoNL}
we obtain the following:
For any $c>\pi(F)$ close enough
to the mean $\pi(F)$, and all $x\in\state$,
\be
\Probsub_x\{S_n\geq nc\}
       \;\sim\; 
        \frac{\cf_a(x)}{a\sqrt{2\pi n\sigma_a^2}}
        e^{-n\LA^*(c)},
        \quad
        n\to\infty,
\label{eq:introBR}
\ee
where $a\in\RL$ is chosen such that $\LA'(a)=c$,
$\cf_a(x)$ is the solution to the multiplicative
Poisson equation \eq expFishb/,
$\LA^*(\cdot)$ is as in (\ref{eq:introLDP}),
and $\sigma_a^2=\LA''(a)$.
A corresponding expansion is given for 
lattice functionals. 

These results generalize those obtained by 
Miller \cite{miller:61} for finite-state chains, 
and those in \cite{jensen:91}, proved
under conditions stronger than Doeblin 
recurrence (in \cite{jensen:91}
a version of the domination assumption 
in \eq sRecurrence/ below is assumed, 
together with additional 
regularity conditions).

The problem of obtaining {\em exact} large
deviations asymptotics 
(such as in (\ref{eq:introBR}) above)
has been considered
by \cite{iscneynum85a,ney-nummelin:87a},
using a ``pinned'' multiplicative mean
ergodic theorem for a $\psi$-irreducible
and aperiodic Markov chain.
It is shown that
for a ``small'' set $C\subset\state$,
\begin{equation}
\limn \smalloneOvern\! \log \Expect_x  [\exp(\alpha S_n )
\ind (\Phi(n) \in C) ]
     = \LA(\alpha),
\elabel{pinnedMMET}
\end{equation}
and from this, under additional conditions
(assuming a variant of the
``uniform domination'' condition
\eq sRecurrence/ discussed
in the following section),
large deviations expansions are 
proved along the same lines as 
indicated above. The difference 
here is that, because of the additional 
constraint imposed by the small set $C$ 
in \eq pinnedMMET/, the resulting
expansions are not for the
probabilities $\Probsub_x\{S_n\geq nc\}$
as in (\ref{eq:introBR}),
but for the ``pinned''
probabilities
$\Probsub_x\{S_n\geq nc\;\;\mbox{and}\;\;
\Phi (n) \in C \}$.

Finally note that in much of the relevant
literature authors often consider a Markov
additive process model instead of simply
the partial sums of a given Markov processes.
For simplicity (and without loss of generality),
we restrict our attention to the asymptotic
behavior of the partial sums themselves.

\subsection{Related Approaches}
%%%%%%%%%%%%%%%%%%%%%%%%%%% Related Work %%%%%%%%%%%%%%%%%%%%%%%%%%% 
\slabel{related}

In this paper we attempt to place within 
a single framework results from two 
previously disparate research areas: 
The theory of positive operators as 
developed in \cite{nummelin:book,niemi-nummelin:86}, where 
$r_\alpha=(\lambda_\alpha)^{-1}$ is
the {\em convergence parameter} for the 
semigroup generated by $\haP_\alpha$,
and from the theory of positive harmonic
functions for diffusions where 
$\LA(\alpha)=\log(\lambda_\alpha)$
is known as the 
{\em generalized principal eigenvalue} 
\cite{pinsky:book}.
The reason that the constant $\lambda_\alpha$
is given two different names is that, so far,
the discrete-time theory of $\psi$-irreducible
Markov chains and the related continuous-time
theory of positive harmonic functions have been
developed independently. Looked at
together, many of the results of the latter
continuous-time theory can be replicated,
improved, or generalized by lifting results
from the discrete-time setting.

These and some other
relevant approaches in the existing 
literature are summarized below.
As this literature is very extensive,
the following discussion is not
intended to be a complete review.

\paragraph{A. $\bfpsi$-irreducible operators.}
The most general approach to understanding 
the eigenfunction equation \eq expFishb/ 
has been developed for discrete-time Markov 
chains, based on renewal theory and the 
theory of positive, $\psi$-irreducible operators; 
see Nummelin's monograph \cite{nummelin:book}.
In this framework $\LA(\alpha) = -\log(r_\alpha)$,
where $r_\alpha$ is the \textit{convergence parameter}
for the semigroup generated by the
kernel $\haP_\alpha$ defined in \eq hatP/.

Although, in general, 
useful solutions to \eq expFish/
cannot be constructed,
if $\bfPhi$ is aperiodic and 
$r_\alpha>0$, then from the 
definitions it can be shown directly 
that for any ``small'' set $C$,
\[
\limn \smalloneOvern\! \log (\haP_\alpha^n (x, C))
         = \LA(\alpha)= -\log(r_\alpha)
        \quad \mbox{a.e.\ } x\in\state\,,
\]
where $\haP_\alpha^n$ denotes the $n$-fold
composition of the kernel $\haP_\alpha$
with itself. From this, the ``pinned'' 
multiplicative mean ergodic theorem 
\eq pinnedMMET/ is easily obtained.
The drawback to this approach is 
the restriction imposed by the small set
$C$ in \eq pinnedMMET/. As we will see,
this restriction is not necessary when 
$\bfPhi$ is geometrically ergodic.
Nevertheless, in the case of 
``first-order'' large deviations
(as opposed to more precise
estimates as in (\ref{eq:introBR})),
these methods provide what appear to
be the most general large-deviations
results to date
\cite{deacosta:90, deacosta-ney:98}.

\paragraph{B. Lyapunov functions and compact sublevel sets.}
A well-behaved solution to the multiplicative 
Poisson equation \eq expFishb/ can be shown 
to exist under suitable bounds on the 
transition kernel $P$. 
For example, \eq expFish/ will admit a 
\textit{bounded} solution $\cf_\alpha$ 
under the ``uniform domination''
assumption of \cite[Sec.~6]{stroock:book}:
For some $\epsy>0$ and all measurable
$A\subset\state$:
\begin{equation}
 P (x, A)  \ge \epsy P (y, A ), \qquad x,y\in\state.
                             \elabel{sRecurrence}
\end{equation}
Condition \eq sRecurrence/,
as well as its variants in 
\cite{iscneynum85a,ellis:88,jensen:91,dembo-zeitouni:book},
are significantly stronger than geometric ergodicity,
and are rarely satisfied for non-compact state 
spaces.  In particular, they imply that the 
process is Doeblin recurrent, 
a property that is equivalent to 
geometric ergodicity with a 
{\em bounded} Lyapunov function $V$;
see \cite[Chapter~16]{meyn-tweedie:book}. 

Similar conditions are used in Donsker
and Varadhan's classic papers; see 
\cite{varadhan:book} for a general 
exposition. 
% Assumption~(3) 
% in \cite[p.~34]{varadhan:book} requires 
% that $\LA(1)$ is finite for an 
% unbounded $F\colon\state\to[0,\infty)$ 
% that has compact sublevel sets
% [recall the definition of $\LA(\cdot)$ in
% \eq expLLN/],
% and assumptions~(I) and~(II)
% in \cite[p.~34]{varadhan:book} imply that 
% a bound similar to \eq sRecurrence/ 
% holds for some $\epsy(K)$ and all $x,y\in K$, 
% whenever $K$ is compact. 
Variations on their assumptions 
are used throughout the large 
deviations literature (including
the recent work by Wu -- see
\cite{wu:00} and the references
therein), and they all imply the 
validity of a condition stronger 
than geometric ergodicity, 
the {\em multiplicative regularity}
condition (mV3), stated and discussed
in \Section{ergodictheorems}.
In particular, Varadhan in \cite{var85a} 
assumes directly that (mV3) holds.

\paragraph{C. Spectral gap.}
In all of the aforementioned works, only the case
where $\alpha\in\Re$ is considered. 
Specifically,
the positivity
of the semigroup generated by the kernel
$\haP_\alpha$ in \eq hatP/ is exploited in
constructing solutions $(\lambda_\alpha,\cf_\alpha)$ 
to the eigenvalue problem \eq expFishb/.
Nagaev in \cite{nag57a} treats the special
case of ergodic Markov chains
that converge to the stationary
distribution at a {\em uniform} geometric rate,
\be
|P^t(x,A)-\pi(A)|\leq B_0e^{-b_0t},
\qquad \mbox{for all}\;x=\Phi(0),\; 
\mbox{all measurable}\;A\subset\state. 
\label{eq:nagaev}
\ee
This condition
is equivalent to Doeblin recurrence.
A version of the multiplicative mean
ergodic theorem is proved,
under (\ref{eq:nagaev}),
for purely imaginary $\alpha=i\omega$ 
in a neighborhood of zero.
The gist of this approach is to
formulate the problem in a 
vector-space setting similar to that 
considered here. Noting that the transition
semigroup $\{P^n\}$ of the Markov chain
converges in operator norm to the invariant 
probability measure (as $n\to\infty$,
where $P^n$ is viewed as a linear operator 
from $L_\infty\to L_\infty$),
the continuity of the norm is exploited 
to obtain convergence of the semigroup 
$\{\haP_\alpha^n\}$.

Operator-theoretic approaches have been 
extensively used in the classical theory 
of Markov chains, and the assumption of 
{\em uniform} geometric ergodicity 
(\ref{eq:nagaev}) is
traditionally used to ensure a spectral 
gap, and hence convergence, as in 
\cite{nag57a}. Generalizations have 
typically involved an alternative 
vector-space setting, such as an $L_p$ 
space for $p<\infty$; see 
\cite{wei84a,jen87a,bolthausen-et-al:95}
and also \cite{fleming-sheu:97,fleming:97}.
%%IK2:
In particular, under the assumption of
{\em hypercontractivity}, Deuschel
and Stoock \cite{deuschel-stroock:book}
derive large deviations properties for
Markov chains. Note that, as hypercontractivity 
implies $L_2$-ergodicity at an exponential rate, 
it also implies that (V4) holds 
\cite{meyn-tweedie:book}.

In a different vain, in
\cite{meyn-tweedie:book,meyn-tweedie:94b} 
the \textit{weighted}-$L_\infty$
space is considered,
\[
   L_\infty^V \eqdef 
        \{ g\colon\state\to\Co \,:\,\sup_x [|g(x)|/V(x)]<\infty\},
\]
with $V\colon\state\to [1,\infty)$ 
being the Lyapunov function 
in condition (V4). The
convergence of the semigroup 
$\{P^n\}$ in the induced 
operator norm on this space 
is  equivalent to 
geometric ergodicity 
\cite{meyn-tweedie:book,meyn-tweedie:94b}, and 
based on this equivalence
we show in this paper that (V4) 
leads to multiplicative mean ergodic 
theorems of the type 
\eq expLLN/--\eq expFish/ 
for complex $\alpha$,
and also to criteria for the 
existence of solutions to \eq expFishb/
under conditions far weaker
than those used in, for example, 
\cite{pinsky:book,varadhan:book}.
We also substantially strengthen the conclusions of both
\cite{pinsky:book}
and \cite{ney-nummelin:87a,ney-nummelin:87b}
since we can apply the $V$-uniform ergodic theorem
of \cite{meyn-tweedie:book} to obtain uniform 
geometric convergence in \eq expLLNb/.

In earlier work related to the ergodic theory of
Markov processes (as opposed to the multiplicative
ergodicity and large deviations issues considered
here), Kartashov considered weighted norms in 
\cite{kar85a,kar85b}, and 
a version of the $V$-uniform ergodic theorem 
for countable state space chains
first appeared in \cite{horspi92a}.

\paragraph{D. Nonlinear semigroups.}
For a continuous-time Markov process $\bfPhi$
(typically a diffusion), Fleming \cite{fle78a} 
and Feng \cite{fen99a} consider a nonlinear 
operator $\clH$ defined as a modification
of the generator $\clA$ of the process $\bfPhi$:
\ben
\clH(G) \eqdef  \log((g^{-1}) \clA g)\, ,\quad \hbox{where $g=e^G$.}
\een

For any function $F\in L_\infty$, the 
multiplicative Poisson equation is given 
in continuous time as $\clA \cf = \exp(-F+\LA) \cf$, 
where $\LA = \LA(1)$
[recall the definition of $\LA(\cdot)$ in
\eq expLLN/].
If $g=\cf$ is a solution for a given $F$, 
then
\[
\clH(G) =  \log[(g^{-1})e^{\LA - F}g] = -F+\LA.
\]
Define the functional $\clG$ on $L_\infty$
as $\clG(F)=\log(c\cf)$,
where $\cf$ solves the multiplicative Poisson equation and
$c=\pi(\cf)^{-1}$ is a normalizing constant. 
The operator $\clG$ is an inverse of $-\clH$ 
in the sense that $\clH\circ\clG= -I$ on some 
appropriately defined domain.

Under (V4), the results of the present paper imply that
$\clG$ is a {\em bounded} nonlinear operator, whose domain 
contains an open ball in $L_\infty$ centered at the origin. 
In particular, our results provide methods for verifying 
the structural assumptions of \cite{fen99a,feng-kurtz:preprint}.
%%IK3:
A thorough investigation of this nonlinear structural theory
and its intimate relationship to large deviations properties
is carried out in the subsequent work 
\cite{kontoyiannis-meyn:II}, for Markov processes 
satisfying the stronger assumption 
of multiplicative regularity.

\newpage

\paragraph{Organization.}
The rest of the paper is organized as follows.
In \Section{background} we collect the basic 
notation and definitions that will remain in 
effect throughout the paper. We present 
background results from the ergodic theory
of Markov chains and processes, and 
briefly discuss several different conditions
for ergodicity and the relationships between
them. 

In \Section{gpe} we collect some results 
about the convergence parameter of a 
positive semigroup. \Section{spectral} 
develops the spectral theory and 
multiplicative ergodic theory
along the lines discussed above. 
Analogs of \eq expLLN/--\eq expFish/ are
proved for geometrically ergodic Markov
processes.

Sections~5 and~6 contain
the probabilistic results 
outlined in \Section{intro-limthms}.
Finally in \Section{examples} 
we give numerous examples of 
Markov chains and processes 
satisfying the assumption of 
geometric ergodicity.

\section{Ergodicity}
%%%%%%%%%%%%%%%%%%%%%%%%%%% BACKGROUND %%%%%%%%%%%%%%%%%%%%%%%%%%% 
\slabel{background}

In this and the following section we review some
necessary background results from certain parts 
of the ergodic theory of Markov chains and processes 
\cite{meyn-tweedie:book}, and some results regarding 
the convergence parameter of a positive semigroup 
as defined in \cite{nummelin:book}. All of this 
concerns a $\psi$-irreducible and aperiodic 
chain or process  
$\bfPhi$ on a general state space 
$\state$
(see below for precise definitions).
We assume that $\state$ is equipped
with a sigma-field 
$\clB$, and that $\clB$
is countably generated.
The distribution of $\bfPhi$
is described 
by a transition semigroup
$\{P^t\,:\, t\in \TT\}$, where $\TT$ is 
taken to be either
the nonnegative integers $\nat$ (in discrete-time)
or the nonnegative reals $\Re_+$ (in continuous-time),
and where for each $t$, $P^t$ is the 
transition kernel
$$P^t(x,A)\eqdef\Pr\{\Phi(t)\in A\,|\,\Phi(0)=x\},
\quad x\in\state,\,A\in\clB.$$
Recall that $P^t$ acts on functions
$f:\state\to\RL$ and signed measures $\nu$ 
on $\clB$, via
\be
P^tf(\cdot)=\int_{\state}P^t(\cdot,dy)f(y)
\;\;\;\mbox{and}\;\;\;
\nu P^t(\cdot)=\int_{\state}\nu(dx)P^t(x,\cdot),
\label{eq:act}
\ee
respectively.

\subsection{$\psi$-Irreducibility}
%%%%%%%%%%%%%%%%%%%%%%%%%%% IRREDUCIBILITY %%%%%%%%%%%%%%%%%%%%%%%%%%%
\slabel{irr}

For any $\theta >0$, we defined the {\em resolvent kernel}
$R_\theta$ by,
\begin{equation}
R_\theta \eqdef
\left\{
\begin{array}{lr} 
\displaystyle
\sum_0^\infty (1-e^{-\theta})  e^{-\theta n} P^n 
          &\hbox{discrete-time}
\\
\\ 
\displaystyle
\int_{[0,\infty)} \theta e^{-\theta t} P^t
\,dt  
            &\hbox{continuous-time,}
\end{array}
\right.
\elabel{resolve}
\end{equation} 
and we write $R$ for $R_1$.

If for some $\sigma$-finite measure $\psi$ on 
$\clB$, some $\theta >0$, and all functions 
$s: \state \to [0,\infty)$ with
$\psi(s) =\int s(x)\, \psi(dx) > 0$,
we have
\[
R_\theta(x, s) \eqdef \int_\state R_\theta(x,dy) s(y)  > 0, \quad x\in \state,
\]
then the semigroup $\{P^t : t\in\TT\}$
is called 
\textit{$\psi$-irreducible}, and $\psi$ 
is called an {\em irreducibility measure.}
If the transition semigroup $\{P^t\}$ associated
with the Markov process $\bfPhi$ is 
$\psi$-irreducible,
then we say that 
$\bfPhi$ is $\psi$-irreducible.
The set of functions $s\colon\state\to\Re_+$ with
$\psi(s) =\int s(x)\, \psi(dx) > 0$
is denoted by $\clB^+$, and all
such $s$ are called $\psi$-positive.

Throughout the paper, we
will assume that $\bfPhi$ is $\psi$-irreducible.
Moreover, without loss 
of generality we assume 
that $\psi$ is 
\textit{maximal} in the sense that 
any other irreducibility measure $\psi'$ 
is absolutely continuous with respect
to $\psi$ \cite{meyn-tweedie:book}. 
We will also assume that 
the semigroup $\{P^t : t\in\TT\}$
is \textit{aperiodic}, that is,
for any $s\in\clB^+$ and any initial condition $x$,
\[
P^t(x,s) > 0\qquad \hbox{for all $t$ sufficiently large.}
\]
If the semigroup associated
with the Markov process $\bfPhi$ is
aperiodic, then we say that $\bfPhi$ is
aperiodic.

A measurable subset $C$ 
of $\state$ is called 
\textit{full} if $\psi(C^c) =0$, 
and it is called \textit{absorbing} if 
$R_\theta (x, C^c)=0$
for $ x\in C$ (for some $\theta$). 
We recall that,
for a $\psi$-irreducible
$\bfPhi$, a non-empty absorbing set
is always full 
\cite[Proposition~4.2.3]{meyn-tweedie:book}.

A function $s\in\clB^+$ 
and a measure $\nu$ on $\clB$ are 
called \textit{small} if, for some 
$\theta>0$, 
\begin{equation}
R_\theta (x, A) \ge s(x) \nu (A),\qquad 
x\in \state, A \in \clB\,.
\elabel{small}
\end{equation} 
In \cite[Proposition~5.5.5]{meyn-tweedie:book} 
it is shown that for a $\psi$-irreducible $\bfPhi$,
one can always find a $\theta$ and a pair 
$(s,\nu)$ satisfying the bound \eq small/, 
such that $s(x) > 0$ for all $x$, and with 
$\nu$ equivalent 
to the maximal irreducibility measure 
$\psi$ (in the sense that they are mutually
absolutely continuous). A similar 
construction works in continuous-time
as well.

If a small function $s$ is of the form
$s= \epsy \ind_C$ for some $\epsy>0
$ and $C\in \clB$, then the set $C$ 
is called {\em small}. 
We denote by $\clB_p^+$ the set of all small 
functions $s\in\clB^+$, and we denote by
$\clM_p^+$ the set of all (positive) small 
measures $\nu$ which satisfy \eq small/ for 
some $s\in \clB_p^+$. Both $\clM_p^+$ and 
$\clB_p^+$ are positive cones, 
and they are closed under addition. 
 
\subsection{Ergodicity Conditions}
%%%%%%%%%%%%%%%%%%%%%%%%%%% ERGODICITY CONDITIONS %%%%%%%%%%%%%%%%%%%%%%%%%%%
\slabel{ergodictheorems}
 
Let $V:\, \state \to [0, \infty]$ be an extended-real 
valued function, with $V(x_0) < \infty$ for at least 
one $x_0\in\state.$ Let $S_V$ denote the (nonempty) set:
\begin{equation}
\elabel{SV}
S_V = \{x\,:\, V(x) < \infty\}.
\end{equation} 
In most of the results below
our assumptions will guarantee
that $S_V$ is absorbing,
hence full, so that 
$V(x) < \infty$ a.e.\ $[\psi]$.

Let $L^V_\infty$ denote the vector space 
of measurable functions $h:\,\state\to\Co$ 
satisfying
\[
\|h\|_{V} \eqdef \sup_{x\in\state} \frac{|h(x)|}{V (x)} < \infty.
\]
Similarly, $L^f_\infty$ will denote the corresponding
space for an arbitrary nonnegative (measurable)
function $f$ on $\state$. We define the 
$V$-norm $\lll \haP \lll_{V}$
of an arbitrary kernel $\haP=\haP(x,dy)$
by
\begin{equation}
\lll \haP \lll_{V}
\eqdef  \sup \frac{\|\haP h\|_{V}}{\|h\|_{V}},
\elabel{Vnorm}
\end{equation}
where the supremum is over all $h\in L^V_\infty$
with $\|h\|_{V}\neq 0.$

In what follows,
it will be convenient 
to describe some important 
properties of $\bfPhi$ 
in terms of its {\em generator} $\clA$
rather than in terms of its transition
semigroup $\{P^t\}$. For a function $g:\state\to\Co$,
we write $\clA g=h$
if for each initial condition 
$\Phi(0)=x\in\state$ the process 
$\{m(t) \,:\, t \in \TT\}$ defined by
\begin{equation}
m(t) \eqdef \int_{[0,t)} h(\Phi(s))\,  ds   
                              -  g(\Phi(t)),\qquad t\in\TT,
\elabel{martGen}
\end{equation}
is a local martingale with respect to the 
natural filtration $\{\clF_t
=\sigma(\Phi(s),\, 0\le s\le t)\,:\,t\in\TT\}$.
In discrete-time the generator is simply $\clA=P-I$.

Next we introduce two different regularity
conditions on $\bfPhi$, taken from \cite{meyn-tweedie:book}.
As we will see, the first one guarantees 
the validity of ergodic results as in
equations \eq meanLLN/--\eq fish/,
whereas the second one will be used
to prove their multiplicative 
counterparts \eq expLLN/--\eq expFish/;
see \Section{spectral}.

Throughout this paper we assume that 
the function $V$ is finite for at least
one $x\in\state$. 

$$
\parbox[b]{3.5in}{For a function $f:\state\to[1,\infty)$, 
        a probability measure $\nu$ on $\clB$, 
        a constant $b<\infty$, 
        a function $s:\state\to(0,1],$
        and a $V:\state\to(0,\infty]$:}
\qquad 
\begin{array}{rcl}
\clA V &\le & -f + b s
\\
R&\ge & s\otimes\nu.
\end{array}
\eqno{\hbox{\bf (V3)}}
$$

$$
\parbox[b]{3.5in}{For a probability measure $\nu$ on $\clB$, 
        some constants $b<\infty$ and $\delta>0$,
        a function $s:\state\to(0,1],$
        and a $V:\state\to[1,\infty]$:}
\qquad 
\begin{array}{rcl}
\clA V &\le & -\delta V + b s
\\
R&\ge & s\otimes\nu.
\end{array}
\eqno{\hbox{\bf (V4)}}
$$ 

Note that condition (V4) is stronger than (V3):
When (V4) holds, (V3) also
holds with  $f=V$, $V'=V/\delta$,
and $b'=b/\delta$.
The assumption that (V4) holds for a Markov chain
$\bfPhi$ is the main condition required for most
of our results, and it will remain in effect 
essentially for the rest of the paper. To 
formalize this assumption we introduce the
following definition:

\begin{quote}
{\bf Geometric Ergodicity.} A Markov process
$\bfPhi$ is called {\em geometrically
ergodic} (with Lyapunov function $V$),
if it is $\psi$-irreducible,
aperiodic, and it satisfies condition (V4)
(with this $V$).
\end{quote}

In \Section{examples}, numerous examples 
are given for which the validity of (V4) is
explicitly verified; see also 
\cite[Chapter~16]{meyn-tweedie:book}.
%%IK4: 
For comparison, we also introduce
the following related condition
for continuous time Markov processes,
which we think of as the natural 
multiplicative analog of (V3):
$$
\parbox[b]{3.5in}{For a function $f:\state\to[1,\infty)$, 
        a probability measure $\nu$ on $\clB$, 
        constants $\delta>0$ and $b<\infty$, 
        a function $s:\state\to (0,1],$
        and a $V:\state\to[1,\infty]$:}
\qquad 
\begin{array}{rcll}
\log\Big(e^{-V}\clA e^V\Big) &\le & -\delta f + b s \quad  \\
R&\ge & s\otimes\nu.
\end{array}
\eqno{\hbox{\bf (mV3)}}
$$ 

As discussed in the introduction, 
condition (mV3) is very closely
related to the conditions in the
well-known Donsker-Varadhan
large deviations results.
In particular, under the conditions 
of \cite{varadhan:book}, especially 
Assumption~(3) in \cite[p.~34]{varadhan:book},
% that requires $\LA:=\lim_n \smalloneOvern \!\log \Expect_x
%         [ \exp \{\int_{[0,t)}F(\Phi(t))dt\}]<\infty$
% for an unbounded $F\colon\state\to[0,\infty)$ 
% with compact sublevel sets, 
it follows from 
\Theorem{evector} below that (mV3) is 
satisfied.
%  with $f=F$, $s=\LA \ind_S$, and $V=\cf$, 
% where $S=\{ x \,:\, F(x)\le \LA\}$,
% and $\cf$ is the eigenfunction
% constructed in \Theorem{evector}. 
Moreover, in the general case
where the state space $\state$
is not compact,
Varadhan's conditions
imply that (mV3) holds
with an unbounded $f$ with compact 
sublevel sets, an assumption already
stronger than (V4) as the following
Proposition shows. A detailed study
of Markov processes satisfying (mV3)
is given in \cite{kontoyiannis-meyn:II},
where the analog of (mV3) for discrete-time
Markov processes is also given.
In the context of diffusions, more 
specific results can be found in
\cite{huimeysch01a}. \Proposition{v5}
is proved in \cite{kontoyiannis-meyn:II}.

%%% ---------------------------------
%%% OLD DISCUSSION
%%%
%%% The function 
%%% $\cf$ is bounded below uniformly 
%%% in $x$ since the Markov process $\bfPhi$ 
%%% is assumed to be recurrent (the argument 
%%% used in Theorem~3.1 of \cite{balaji-meyn} 
%%% may be applied). Hence, by scaling $\cf$, 
%%% we may assume that the required bound 
%%% $V\ge 1$ holds.  Conversely, if (mV3) holds,
%%% it follows that $\Lambda\le b<\infty$.
%%% 
%%% ---------------------------------

\begin{proposition}
\tlabel{v5}
Suppose $\bfPhi$ is $\psi$-irreducible and aperiodic.
If (mV3) holds, then so does (V4).
\end{proposition}

\newpage

\subsection{Ergodic Theorems}
%%%%%%%%%%%%%%%%%%%%%%%%%%% ERGODIC THEOREMS %%%%%%%%%%%%%%%%%%%%%%%%%%%
\slabel{ergodictheorems2}

Under either (V3) or (V4),
there exists 
a unique invariant probability measure
$\pi$ on $\clB$ (see below). Given such a $\pi$,
we define $\Pi$ as the kernel
\[
\Pi=\one\otimes \pi,
\]
so that $\Pi(x,A)=\pi(A)$, $x\in\state$,
$A\in\clB$.  If 
$\pi(V)\eqdef\int_\state\pi(dx)V(x)<\infty$, 
then
$\Pi$ acts on $L^V_\infty$
as a bounded linear operator.

A \textit{fundamental kernel} is a 
linear operator 
$Z\colon L_\infty^f \to L_\infty^V $
(for some measurable
functions $f\ge 1$, $V\ge 1$),
satisfying
\begin{equation}
   \clA Z = - (I-\Pi)\, .
\elabel{Ainverse}
\end{equation}
That is, for any $F\in L_\infty^f$, 
the function $\haF=Z F\in L_\infty^V$ 
solves the {\em Poisson
equation},
\begin{equation}
    \clA \haF =  -F + \pi(F),
\elabel{Poisson}
\end{equation}
where $\pi(F)=\int_\state\pi(dx)F(x)$.  Equivalently, the stochastic
process 
\begin{equation}
m(t) = \haF(\Phi(t)) - \haF(\Phi(0)) + \int_{[0,t)} (F(\Phi(r)) - \pi(F))\, dr,
\qquad t\ge 0,
\elabel{MartFish}
\end{equation}
is a martingale with respect to $\{\clF_t\}$.

The following two theorems give equivalent 
conditions for $\bfPhi$ to be ergodic
or geometrically ergodic, respectively.
\Corollary{ethms} states that a fundamental 
kernel exists, and the ergodic 
results \eq meanLLN/--\eq fish/
given in the introduction 
indeed hold as soon as $\bfPhi$
satisfies (V3).

For any $C\in\clB$, let $\tau_C$ denote the hitting time
\[
\tau_C \eqdef \inf \{ t \ge 1\,:\, \Phi (t) \in C\}.
\]
\begin{theorem}
\tlabel{Fergodic}
{\em (Ergodicity)}
Suppose that $\bfPhi$ is $\psi$-irreducible and aperiodic.  
For any function $f:\state\to[1,\infty)$
the following are equivalent:
\begin{description}
\item[(i)] 
The process $\bfPhi$ is positive recurrent with
invariant probability measure $\pi$, and
$\pi(f)<\infty.$
\item[(ii)]
There exists a small set $C$ such that
\[
\sup_{x\in C} \Expect_x 
\Bigl[  \int_{[0,\tau_C )} f(\Phi(t)) dt \Bigr] < \infty\,.
\]
\item[(iii)]
Condition (V3) holds with the same $f$.
\end{description}
If any of these conditions holds, then
the set $S_V$ defined in \eq SV/ is absorbing 
and full, and
\begin{equation}
\sup_{g:|g|\leq f} \,
        | P^t(x,g) - \pi(g) | \to 0,
               \qquad t\to\infty,
                  \quad x\in S_V.
\elabel{f-norm}
\end{equation}
Moreover, for any small measure $\nu$
there exists a fundamental kernel 
$Z$ which 
is a bounded linear operator,
\[
Z\colon L_\infty^f \to 
\{ h\in L_\infty^V \,:\, \nu(h) = 0\}.
\]
If $Z'$ is any other such fundamental kernel,
then $\| ZF - Z'F \|_V=0,\,$ $F\in L_\infty^f$.
\end{theorem}

\proof
The discrete-time version of (i)-(iii) is a consequence of  
the $f$-Norm Ergodic Theorem of \cite{meyn-tweedie:book}, 
and the continuous-time version follows from 
\cite[Theorem 5.3]{meytwe93b}. The construction of the 
fundamental kernel and the uniform 
bound is given in \cite[Theorem 2.3]{glymey96}.
\qed

\medskip

Recall the definition of $S_t$ in \eq psums/.
\begin{corollary}
\tlabel{ethms}
{\em (Ergodic Theorems)}
Let $\bfPhi$ be a $\psi$-irreducible, 
aperiodic Markov process that satisfies 
(V3). 
If $F\in L_\infty^V$ 
with $f$ as in (V3), then for $x\in S_V$:
\begin{description}
\item[(a)] $\;\Expect_x[\smalloneOvert S_t]\to \pi(F)\;$
        as $t\to\infty.$

\item[(b)] There exists $\haF\in L_\infty^V$ 
        with $\pi(\haF)=0$, so that $\haF$
        solves the Poisson equation
\[
\clA \haF =  -F + \pi(F).
\]

\item[(c)] If, in addition, $\pi(V)<\infty$,
        then $\haF$ satisfies
\be
\haF(x)\;=\;\lim_{t\to\infty}\Expect_x[S_t-t\pi(F)]\, .
\label{eq:Fhat}
\ee
\end{description}
\end{corollary}

\proof 
The convergence in norm \eq f-norm/ implies the convergence,
\[
\Expect_x[F(\Phi(t))] \to \pi(F),\qquad t\to\infty,\qquad x\in S_V,
\]
which gives (a).  For (b) we can take $\haF=Z F$, where $Z$ is 
given in \Theorem{Fergodic}.

When $\pi(V)<\infty$ it follows that, for some $b_1<\infty$,
\[
\int_{\TT}\, \Bigl| \Expect_x[F(\Phi(t))] - \pi(F) \Bigr| \, dt
\le 
                b_1 \|F\|_f V(x) \, ,\qquad x\in S_V \, .
\]
This is given as
Theorem~14.0.1 of \cite{meyn-tweedie:book} 
in discrete-time.  The continuous-time case follows on considering
the skeleton chain $\Phi(\delta k)$, $k=1,2,3,\dots$, as discussed
on p.~247 of \cite{meyn-tweedie:93e}.
This implies that one version of the fundamental kernel 
may be expressed as
\[
Z(x,A) = \int_{\TT}  (P^t(x,A) - \pi(A))  
        \, dt\, ,\qquad x\in S_V ,\ A\in\clB,
\]
and $Z$
is a bounded linear operator
from $L^f_\infty$ to $ L_\infty^V$.
This gives (c).
\qed

\medskip

The solution $\haF$ of the Poisson equation 
given in \Corollary{ethms}~(b)
arises in almost every limit theorem considered 
below. In particular, it can be used to define 
the {\em asymptotic variance}
$\sigma^2$ in the central 
limit theorem; see Theorem~17.4.5 
of \cite{meyn-tweedie:book}.
In the discrete-time case, $\sigma^2<\infty$
as soon as $\pi(\haF^2)<\infty$,
and equation~(17.44) of \cite{meyn-tweedie:book} 
gives the representation,
$$
\sigma^2 = \lim_{n\to\infty}\frac{1}{n}\VAR_x\{S_n\}
        = \Expect_{\pi} \Bigl[ \haF (\Phi (n))^2 - 
                (P \haF (\Phi (n)))^2 \Bigr],\qquad x\in S_V.
$$

Next we obtain a characterization
of the case when $\sigma^2 = 0$.
In discrete time, a similar result
is derived in \cite{bezhaeva-o:96}
using different methods.

\begin{proposition}
\tlabel{clt}
{\em (Variance Characterization)}
Suppose that $\bfPhi$ satisfies (V3) with $\pi(V^2)<\infty$.
Then,  for any $F\in L_\infty^f$, the 
{\em asymptotic variance}
\begin{equation}
\sigma^2\eqdef \lim_{t\to\infty}\frac{1}{t}\VAR_x\{S_t\}
\elabel{var1}
\end{equation}
exists for any initial condition $\Phi(0)=x\in S_V$.
Writing  $\barF = F - \pi (F)$, $\sigma^2$ satisfies,
for all $t>0$,
\begin{equation}
\sigma^2 = \frac{1}{t}
\Expect_\pi\left[\left(\haF(\Phi(t)) - \haF(\Phi(0))
        + \int_{[0,t)} \barF(\Phi(s)) ds\right)^2
        \right]<\infty\, .
\elabel{var2}
\end{equation}
Moreover:
\begin{description}
\item[(i)] 
If $\sigma^2 = 0$ then there exists $G\in L^V_\infty$,
satisfying, 
\begin{equation}
\int_{[0,t)} \barF(\Phi(s)) ds  =  G(\Phi(t)) - G(\Phi(0)) \, ,
\qquad a.s.\;[\pi].
\elabel{fishEq4}
\end{equation}
When time is discrete, this can also 
be expressed 
as $P (x, S_x) = 1$,
$x \in S_V$, where $S_x =\{ y\in\state \,:\, G (y) = G (x) - \barF (x) \}.$

\item[(ii)]
Conversely, if \eq fishEq4/ holds for some  $G\in L^V_\infty$,
then $\sigma^2 = 0$.
\end{description}
\end{proposition}

\proof
Result (i) is an immediate consequence of \eq var2/, 
which follows from the martingale characterization 
of $\haF$ (see \eq MartFish/).
Result (ii) is immediate from \eq var1/.  
\qed

\begin{theorem}
\tlabel{Vuniform}
{\em (Geometric Ergodicity)}
Suppose that $\bfPhi$ is $\psi$-irreducible and aperiodic.  
The following
are equivalent:
\begin{description}
\item[(i)]
There
exists a probability measure $\pi$ and a
$V\colon\state\to [1,\infty]$, such that 
$P^t$ converges to $\pi$
in the $V$-norm, 
$$
P^t\to \one\otimes \pi,\quad t\to\infty.
$$
\item[(ii)]
There exists a small set $C$ and $\epsy > 0$ such that
\[
\sup_{x\in C} \Expect_x \Bigl[ \exp\Bigl(\epsy \tau_C  \Bigr) \Bigr] < \infty.
\] 
\item[(iii)]
Condition (V4) holds for some $V\colon\state\to [1,\infty]$.
\end{description}
If any of these conditions holds, then
the set $S_V$ defined in \eq SV/ is absorbing and full
for any function $V$ satisfying (iii),
and there exist constants  
$b_0>0$, $B_0,B'_0<\infty$
and an invariant probability 
measure $\pi$ on $\clB$,
such that 
\ben
\lll P^t - \one\otimes \pi \lll_{V} 
&\le&  B_0e^{-b_0t}
        \, , \qquad t\in \TT\, ,\\
|\Expect_x[S_t-t\pi(F)]-\haF(x)|
&\le&  B'_0 \|F\|_V e^{-b_0t}\, ,
                \qquad 
                        F\in L_\infty^V
                \, ,
                        x\in S_V
                \, , 
                        t\in \TT\, .
\een
\end{theorem}

\proof
The discrete-time result is 
\cite[Theorem~15.0.1]{meyn-tweedie:book},
and the continuous-time version is the main result of
\cite{dowmeytwe}.
\qed

\section{The Generalized Principal Eigenvalue}
%%%%%%%%%%%%%%%%%%%%%%%%%%% G.P.E. %%%%%%%%%%%%%%%%%%%%%%%%%%%
\slabel{gpe}

As in the previous section, we assume
that $\bfPhi$ is a geometrically ergodic
Markov process. We fix throughout this 
section a bounded (measurable)
function $F:\, \state\to \Re$ and 
a real number $\alpha\in\RL$,
and we define the semigroup
$\{\haP^t_\alpha\,:\, t\in \TT\}$ by 
\begin{equation}
\haP^t_\alpha (x,A) = \Expect_x [\exp (\alpha S_t) \ind_A (\Phi (t)) ],
\qquad t\in\TT\,,
\elabel{haP}
\end{equation} 
where $S_t$ is defined as before
by \eq psums/.
In this section we 
consider the general properties of this
positive semigroup;
we therefore suppress the dependency 
on $\alpha$ and $F$ and simply
write $\{\haP^t\}$ for
$\{\haP^t_\alpha\,:\, t\in\TT\}$. 

We say that an arbitrary positive kernel $\haP$ is
{\em probabilistic} if $\haP(x,\state)=1$
for all $x\in\state$. Similarly, a semigroup 
$\{\haP^t\}$ is called {\em probabilistic}
if $\haP^t$ is probabilistic for all $t$.
Clearly, the semigroup $\{\haP^t\}=\{\haP^t_\alpha\}$
is, in general, {\em non-probabilistic}.
The definitions of $\psi$-irreducibility and
aperiodicity carry over to non-probabilistic
semigroups immediately. 
[Further extensions to kernels defined
for {\em complex} numbers $\alpha\in\Co$ 
will be treated in \Section{spectral}.]
Note that $\{\haP^t\}$ 
is irreducible and aperiodic as soon as $\{P^t\}$ is.
We also define a family of resolvent kernels 
$\haR_\theta$ for $\theta>0$ exactly as
in \eq resolve/, with $\haP^t$ in place of $P^t$.
To ensure that these are finite for all
$x\in\state$ and a suitable class of $A\in \clB$,
we usually consider $\theta$'s 
in the range $\theta>|\alpha|\|F\|_\infty.$
(And as before, we write $\haR=\haR_1$.)
With $\haR_\theta$ replacing $R_\theta$,
the definitions of small functions, measures, 
and sets carry over verbatim. 

Finally, we define the {\em generator}
$\haclA$ of the semigroup
$\{ \haP^t\,:\,t\in\TT\}$: We
write $\haclA g=h$ if
\be
\haP^t g\, (x)
= g(x) + \int_{[0,t)} \haP^s h\, (x) \, ds,\qquad t\in\TT,\ x\in\state.
\label{eq:generator}
\ee
The following resolvent equations 
will play a central role in a lot of
what follows:
\begin{equation}
\begin{array}{rclr}
\haclA \haR_\theta &=& (e^\theta - 1) (\haR_\theta -I)
&\hbox{discrete-time}
\\
\\
\haclA \haR_\theta &=&  \theta (\haR_\theta - I) 
&
\hbox{continuous-time}
\end{array}
\elabel{ResEqn}
\end{equation}
In continuous-time, 
the resolvent equation 
can be used to establish 
the following identity, whenever
the sum and integral converge absolutely,
\begin{equation}
\sum_{n=1}^\infty \haR_\theta^n z^{-n} =  \theta z^{-1}
\int_{[0,\infty)} e^{-(1-z^{-1}) \theta t} \haP^t\, dt  \, ,\qquad z\in\Co.
\elabel{gpekey}
\end{equation}
This is tremendously valuable in consolidating 
continuous- and discrete-time theory.

Since the semigroup 
$\{\haP^t\}$
is $\psi$-irreducible,
the kernel $\haR_\theta$ satisfies the 
following minorization condition: There
are $s\in \clB_p^+$ and  $\nu \in \clM_p^+$
such that
\[
\haR_\theta \ge s \otimes \nu.
\]
[Note that, since the semigroup $\{\haP^t\}$
is derived from $\{P^t\}$, the above domination
condition is satisfied with $s$ and $\nu$ that
are small with respect to $\{P^t\}$.]
Let $\{\kappa_n \,:\, n\ge 1\}$ 
denote the positive sequence
defined by:
\[
\kappa_n = \nu (\haR_\theta)^{n-1} s, \quad n \ge 1\,.
\]
This sequence 
is supermultiplicative,
\begin{eqnarray*}
\kappa_{n+m} & = &
\nu \haR_\theta^{n-1} \haR_\theta  \haR_\theta^{m-1} s \\
&\ge & \nu (\haR_\theta)^{n-1} (s\otimes \nu )(\haR_\theta)^{m-1} s \\
& = & \kappa_n \kappa_m,
\end{eqnarray*}
so there exists some 
$L(\theta) \in (-\infty,\infty]$ such that
\[
\smalloneOvern\! \log (\nu \haR_\theta^n s) 
= \smalloneOvern\! \log (\kappa_{n+1}) \to L(\theta), \quad
n \to \infty\,.
\]
The constant
\be
r_\theta \eqdef \exp (- L(\theta))
\label{eq:rtheta}
\ee
is called the \textit{convergence
parameter} for the kernel 
$\haR_\theta$ \cite{nummelin:book}.  
It satisfies:
\[
\sum_{n=0}^\infty [\haR_\theta^n s(x)] r^n 
\quad
\left\{
\begin{array}{llr}
= \infty , &\mbox{for all\ } x\in\state, 
        &  \mbox{if\ } r >r_\theta 
\\ 
< \infty , & \mbox{for a.e.\  } x\in\state\ [\psi], 
        & \mbox{if\ } r<r_\theta.
\end{array}
\right.
\]

To move from the resolvent back to the original semigroup,
we apply the resolvent equations \eq ResEqn/.
These relations establish the major part of the following theorem.

\begin{theorem} 
\tlabel{gpe}
{\em (Generalized Principal Eigenvalue)}
Suppose $\bfPhi$ is
$\psi$-irreducible and aperiodic.
Then there is a $\lambda_\circ \in (0,\infty]$  
such that, for any $s\in\clB_p^+$:
\begin{description}
        \item[(i)] 
\begin{equation} 
\int_\TT \lambda^{-t} \haP^t s(x) \, dt
\quad
 \left \{ 
\begin{array}{llr}
=\infty &   \hbox{for all $x\in\state $,}        
    & \lambda < \lambda_\circ  
\\  \\
< \infty & \hbox{for a.e.\ $x\in\state$\ $[\psi]$,}
    & \lambda > \lambda_\circ .
\end{array}
\right.
\elabel{gpeInt} 
\end{equation} 
        \item[(ii)] \quad
$\displaystyle
\frac{1}{t} \log (\haP^t s(x)) \to \Lambda_\circ 
 \eqdef \log (\lambda_\circ ) 
 \qquad \mbox{a.e.}\   x\in \state \ [\psi],\quad t\to\infty\, .
$
\end{description}
\end{theorem}

% save this comment for future reference:   
% If we define consistently $R_z=z[Iz-A]^{-1}$ then
% one must replace $e^{\theta}$ by $\theta +1$ 
% in the formula for $\lambda_\circ$.

\proof 
Result~(i) follows from Theorem~3.2 of \cite{nummelin:book} 
for discrete-time chains, and from \eq gpekey/ 
for the continuous-time case where we may translate to the discrete-time
case using the resolvent $\haR_\theta$.  Then with
$\lambda_\theta = r_\theta^{-1}$,
\begin{equation}
\begin{array}{rcll}
 \lambda_\circ &=& e^\theta + (1-e^\theta) \lambda_\theta
                                        \quad & \hbox{discrete-time}
\\[.2cm]
 \lambda_\circ &=& \exp(\theta(1-\lambda_\theta^{-1}) )
                                 \quad & \hbox{continuous-time.}
\end{array}
\elabel{evaluesRes}
\end{equation}
The second part follows from an argument similar 
to that used in the proof of Lemma~3.2 of 
\cite{balaji-meyn}.  
\qed 
    
\medskip    

We call the constant $\lambda_\circ$  the 
\textit{generalized principal eigenvalue} (g.p.e.)
of the semigroup $\{\haP^t\,:\, t \in\TT\}$.
This generalizes the corresponding
definition of \cite{pinsky:book},
and, as we will see in 
\Theorem{evector} below, $\lambda_\circ$
does indeed play the role of an eigenvalue.
The interpretation of (ii) 
is the ``pinned''
multiplicative mean ergodic theorem
\eq pinnedMMET/ discussed in the 
introduction,
\begin{equation}
\smalloneOvert\! \log 
\Expect_x [\exp (\alpha S_t ) \ind_C (\Phi (t) ) ] \to \LA(\alpha) 
\,,\qquad t \to \infty\,,
\elabel{miniMMET}
\end{equation} 
for a.e.\ $x\in \state$ \ $[\psi]$.
This follows from taking $s=\epsy\ind_C$
in (ii) with  $C$ small, and  $\epsy>0$.

\Theorem{gpe} leaves open what happens in \eq gpeInt/ when
$\lambda=\lambda_\circ $.  The semigroup $\{\haP^t\}$ is called:
\begin{description}
        \item[(i)] 
\textit{$\lambda_\circ $-transient} if
\[
\int_\TT \lambda_\circ ^{-t} \nu \haP^t s dt< \infty
\]
        \item[(ii)] 
\textit{$\lambda_\circ $-recurrent} if
\[
\int_\TT \lambda_\circ ^{-t} \nu \haP^t s dt= \infty 
\] 
        \item[(iii)] 
\textit{$\lambda_\circ $-geometrically recurrent} if the function
\[ 
(\lambda_\circ -z)\int_\TT  (\nu \haP^t s ) z^{-t}dt
\]
is analytic in a neighborhood of $z=\lambda_\circ$.
\end{description}
In  (i)--(iii),   $(s,\nu)$ is any pair with 
$s\in \clB_p^+$ and $\nu\in \clM_p^+$.  
The particular small function or small 
measure chosen is not important \cite{nummelin:book}.

% Note to self: Twice I have gone to Nummelin's book to eliminate the
% proof of the next Lemma part (i) - I know that this is the construction 
% of the 'minimal harmonic function'.  But I couldn't find a result 
% which handled the non-recurrent case.

The construction of $h$ in part (ii) of the following Lemma is an 
extension of the \textit{minimal harmonic function} in
\cite[Proposition~3.13]{nummelin:book}, where here we allow the semigroup 
$\{P^t\}$ to possibly be transient.

\begin{lemma}
\tlabel{small}
\begin{description}
\item[(i)]
Suppose that 
$\{\haP^t \,:\, t\in\nat\}$ 
has g.p.e.\ $\lambda_\circ<\infty$, and suppose that the following   
minorization condition holds for some 
$s\in \clB_p^+$ and  $\nu\in \clM_p^+$:
\[
 \haP \ge s\otimes\nu.
\]
Then,
\[
\sum_{t=0}^\infty \lambda_\circ ^{-t-1} \nu (\haP-s\otimes\nu)^t s \le 1,
\]
with equality if and only if the semigroup is $\lambda_\circ$-recurrent.

\item[(ii)]
If in (i) we take 
$\{ \haP^t \,:\, t\in\nat\}$ to
be the probabilistic semigroup
$\{ P^t \,:\, t\in\nat\}$, then,
\[
\begin{array}{rcl}
h(x)\eqdef \sum_0^\infty (P - s \otimes \nu)^n s\, (x)
&\le &  1;
\\[.2cm]
\sum_0^\infty (P - s \otimes \nu)^n P s\, (x)
& = &  -s(x) + (1+\nu(s)) h(x) \leq 2,
 \quad 
\mbox{ for all  $x\in\state$.}
\end{array}
\]
In this case, if  $\{P^t \,:\, t\in\nat\}$ is 
$1$-recurrent, $h(x)=1$ 
for a.e.\ $x\in\state$ $[\psi]$.
\end{description}
\end{lemma}

\proof 
Part (i) is Proposition 5.2 of \cite{nummelin:book}. 
The essence of this result is the inversion formula,   
\begin{equation}
\Bigl[Iz - \haP \Bigr]^{-1}  =  
\Bigl[ Iz - (\haP  - s \otimes \nu) \Bigr]^{-1}
\Bigl(I + \frac{1}{1-\kappa} s \otimes \nu \Bigr )
\elabel{inversion}
\end{equation}
where
\[
\kappa =    \nu [Iz - (\haP  - s \otimes \nu )]^{-1} s \, .
 \]
 From \eq inversion/ it may be seen that, for $z >0$,
\begin{equation}
 \hbox{$ \kappa = 1$ 
 \textit{if and only if}
 $ \nu [Iz -  \haP ]^{-1} s =\infty$.}
\elabel{inversionB}
\end{equation}

The proof of (ii) is by induction. Define, for $n\ge 0$,  
\[
h_n  =  \sum_{t=0}^n (P  - s \otimes \nu)^t s.
\] 
For $n = 0$, $h_0=s$ and  $s\leq 1$  by assumption.
If true for $n$, then 
\begin{eqnarray*}
h_{n+1}(x) & = & (P  - s \otimes \nu) h_n\, (x) + s(x) \\ 
        & \leq & (P  - s \otimes\nu)  \one \, (x) + s(x) \\  
        & = & [P  (x,\state) - s (x) \nu (\state)]  + s(x)\\  
        & = & 1,
\end{eqnarray*}
where in the last equation we have 
again used the fact that $s\le 1$. 
It follows that $h(x)=\lim h_n(x) \le 1$
for all $x$.

To see the second bound, write 
\[
( P  - s \otimes \nu)^n P s = ( P  - s \otimes \nu)^{n+1} s 
+ ( P  - s \otimes \nu)^n [s\otimes\nu] s.
\]
Summing over $n$ gives the desired result.
\qed

\medskip
 
The reason we call the constant $\lambda_\circ$ 
a generalized eigenvalue is clarified by 
the next theorem, where it shown that, if the 
semigroup $\{\haP^t\,:\, t\in\TT\}$ is 
$\lambda_\circ$-recurrent,
then there is a function 
$\cf :\ \state \to[0,\infty)$
so that $(\cf,\la_\circ)$ solve
the eigenvalue problem,
\begin{equation}
\haP \cf  = \lambda_\circ  \cf \, .
\elabel{evector}
\end{equation} 
Equation \eq evector/ is an instance
of the {\em multiplicative Poisson equation}.
Conditions for the existence of a solution 
to \eq evector/ based upon \Lemma{small}~(i) are 
well-known in the discrete-time case.  
A candidate solution is given by
\begin{equation}
\cf  =   \sum_{n=0}^\infty r_\theta^n \haR_\theta 
        (\haR_\theta- s \otimes \nu)^n s,
\elabel{cfdef}
\end{equation} 
where $(\theta,s,\nu)$ satisfy  $\theta>\lambda_\circ$,  
$s\in \clB_p^+$,   $\nu \in \clM_p^+$,
$\haR_\theta \ge s \otimes \nu$,
and $r_\theta$ is the convergence 
parameter defined in (\ref{eq:rtheta}).

\begin{theorem} 
\tlabel{evector}
{\em (Existence of an Eigenfunction $\cf$)}
Suppose that $\bfPhi$ is $\psi$-irreducible, 
and that the g.p.e.\ $\lambda_\circ$ of the
positive semigroup
$\{\haP^t\}$ is finite. Then
the function $\cf$ given in \eq cfdef/ is finite a.e.\ $[\psi]$, and
\begin{description}
\item[(i)]
If $\{\haP^t\,:\, t\in\TT\}$ is $\lambda_\circ$-recurrent 
then $\cf$ solves the multiplicative Poisson equation:
\begin{equation}
\haP^t \cf  = \lambda_\circ^t  \cf \, ,
\qquad t\in\TT
\, .
\end{equation}

\item[(ii)]
If $\{\haP^t\,:\, t\in\TT\}$ is $\lambda_\circ$-transient 
then for any small function $s\in\clB_p^+$, 
there exists $\delta>0$ such that
\begin{equation}
\haP \cf  =\lambda_\circ  \cf  -  \delta s . 
\elabel{fIdentity}
\end{equation} 
Hence, in the $\lambda_\circ $-transient  
case there is a solution $\oo f$ 
to the pointwise inequality
\begin{equation}
\haP \oo f \le \lambda_\circ  \oo f\,
\elabel{esubvector} 
\end{equation} 
with $\oo f$ finite a.e. $[\psi]$, and
where the inequality is strict whenever 
$\oo f(x) < \infty$.

\item[(iii)]
The solution \eq cfdef/  is minimal and essentially unique:
If  $\oo{f}\colon\state\to(0,\infty)$
is any solution to the inequality
\eq esubvector/, then there exists $c\in\Re_+$ such 
that $ \oo{f} (x) \ge c\cf (x)$
for all $x$.

If $\{\haP^t\,:\, t\in\TT\}$ is $\lambda_\circ$-recurrent,
then we have $\oo f = c \cf $ a.e.\ $[\psi]$, 
and  $ \oo{f} (x) \ge c\cf (x)$ for all $x$.
\end{description}
\end{theorem}
 
\proof
These results are all based on 
Theorem~5.1 of \cite{nummelin:book} 
in the discrete-time case.
 
If $\theta > 0$ is taken large enough,
then the resolvent $\haR_\theta$ satisfies
$\nu \haR_\theta s<\infty$ for any small $s,\nu$ satisfying 
the domination condition $\haR_\theta \ge s\otimes\nu$.  We then set
\[
\cf _\theta  = \sum_{n=0}^\infty r_\theta^n (\haR_\theta - s \otimes
\nu)^n s\,,
\]
where $r_\theta$ is the convergence parameter for $\haR_\theta$.  
We have
\[
r_\theta (\haR_\theta - s \otimes \nu) \cf _\theta = \cf _\theta - s,
\]
and hence
\[
\cf \eqdef \haR_\theta \cf_\theta 
  = r_\theta^{-1} \cf_\theta -  \delta_\theta s,
\]
where $\delta_\theta =  r_\theta^{-1}- \nu (\cf_\theta) \ge 0$.
This constant is strictly positive if and only if 
the semigroup
$\{\haR^n_\theta\}$
is $r_\theta^{-1}$-transient (see \Lemma{small}~(i)).

Results (i)--(iii) then follow from the resolvent equation
in discrete or continuous-time.
\qed

\section{Spectral Gap and Multiplicative Mean Ergodic Theorems}
%%%%%%%%%%%%%%%%%%%%%%%%%% Spectral Theory %%%%%%%%%%%%%%%%%%%%%%%%%%% 
\slabel{spectral}  

The following assumptions will be held
throughout the remainder of this paper:
\begin{equation}
\parbox{.63\hsize}{\raggedright
\begin{description}
\item[(i)]
The Markov process $\bfPhi$ is geometrically 
ergodic with a Lyapunov function $V:\state \to [1, \infty)$,
such that $\pi(V^2)<\infty$.
\item[(ii)]
The (measurable) function 
$F\colon \state \to [-1,1]$
has zero mean $\pi(F)=0$,
and non-trivial asymptotic
variance
$\sigma^2:=\lim_t\VAR_x\{S_t/\sqrt{t}\}>0$.
\end{description}}
\elabel{assumption}
\end{equation}
Note that the additional assumption 
$\pi(V^2)<\infty$ can be made without
any loss of generality:
(V4) implies (V3) with $f=V$ as
discussed above, which implies 
that $\pi(f)<\infty$
\cite[Theorem~14.0.1]{meyn-tweedie:book}.
Moreover, Lemma~15.2.9 
of \cite{meyn-tweedie:book}
says that (V4) also holds
with respect to $\sqrt{V}$
(and some, possibly different,
small function $s$), so we can
always take $V$ in (V4) such that 
$\pi(V^2)<\infty$.

Until \Section{SpectralCts} 
we specialize to the discrete-time 
case for the sake of clarity.

\medskip

With $V$ as in \eq assumption/,
the \textit{spectrum} $\clS(\haP)\subset \Co$ 
of a bounded linear operator $\haP\colon L_\infty^V\to L_\infty^V$ is 
defined to be the set of nonzero $\lambda\in\Co$ 
for which the inverse $(I \lambda - \haP)^{-1}$ 
does not exist as a bounded linear operator 
on $L_\infty^V$.

Recall that an arbitrary kernel $\haP$
acts on functions (on the right) and
on signed measures (on the left)
as in (\ref{eq:act}). With that in
mind, we think of a kernel $\haP$
as an operator acting on a appropriate
function space. 
The kernel $\haP$ is a bounded linear 
operator on $L_\infty^V$ provided its $V$-norm 
$\lll\haP\lll_V$ is finite,
since this is precisely the induced operator norm.
For an arbitrary linear operator
$\haP\colon L_\infty^V\to L_\infty^V$ we 
continue to define the norm 
$\lll\haP\lll_V$ as in \eq Vnorm/.
Also we recall that $\haP$ acts on
a suitable space of measures (on the
left) as
$$\nu\haP(A)\eqdef\nu(\haP\ind_A),\qquad
A\in\clB.$$

For $\alpha\in\Co$ the kernel
$\haP_\alpha$ defined in \eq haP/
yields an operator
$\haP_\alpha\colon L_\infty^V\to L_\infty^V$ 
acting via
\be
\haP_\alpha g\, (x) = \exp(\alpha F(x)) Pg\, (x)
,
\qquad          x\in\state,\ g \in L_\infty^V\, .
\elabel{haP-discrete}
\ee
Its spectrum is denoted $\clS_\alpha = \clS(\haP_\alpha)$. 
The $n$-fold composition of the kernel
$\haP_\alpha$ with itself acts on $L_\infty^V$ as
\[
\haP^n_\alpha g\, (x) = \Expect_x[\exp(\alpha S_n) g(\Phi(n))]
,
\qquad          x\in\state,\ g \in L_\infty^V,\ n\ge 1,
\]
where $\{S_n \,:\, n\ge 1\}$ denote the partial sums \eq discretePS/.
Letting $\clM_1^V$ denote the 
space of signed and possibly complex-valued
measures $\mu$ satisfying 
$|\mu|(V)<\infty$, we obtain analogously,
\[
\mu \haP^n_\alpha \, (A)
= \int\Expect_x[\exp(\alpha S_n) \ind(\Phi(n)\in A)]\, \mu(dx)
,               \qquad          A\in\clB,\ \mu \in \clM_1^V,\ n\ge 1\, .
\]

In this section we identify a region
$\Omega\subset\Co$ such that, 
for geometrically Markov chains,
eigenfunctions $\cf_\alpha\in L_\infty^V$
and (positive) eigenmeasures 
$\cmu_\alpha\in\clM_1^V$ exist
for $\haP_\alpha$, corresponding to 
a given eigenvalue $\lambda_\alpha\in\clS_\alpha$
and $\alpha\in\Omega$.
Suppose that such $\cf_\alpha,\cmu_\alpha$ are found,
and assume that they are normalized so that
\begin{equation}
\cmu_\alpha(\cf_\alpha) = \cmu_\alpha(\state)=1\, .
\elabel{cmu-norm}
\end{equation}
We then let  
$\haQ_\alpha\colon L_\infty^V \to  L_\infty^V$ 
denote the operator 
$\haQ_\alpha = \cf_\alpha\otimes\cmu_\alpha$, 
\[
\haQ_\alpha g\, (x)  = \cmu_\alpha (g) \cf_\alpha(x) ,  
                     \qquad  g \in L_\infty^V,\ x\in\state\, .
\] 
Note that $\haQ_\alpha$ is a {\em projection} operator,
that is, $\haQ_\alpha^2 = \haQ_\alpha$.

%%% YIANNIS-- Notes to myself:
%%% Note that for a $\psi$-irreducible and aperiodic
%%% Markov process $\bfPhi$ for the definition of 
%%% a lattice function we can equivalently assume
%%% that (\ref{eq:lattice}) holds for {\em all}
%%% $x\in\state$, since for any full 
%%% $C\subset\state$ there is an absorbing 
%%% and full $C'\subset C$.

The main results of this section are 
summarized in the following two theorems.
In particular, the multiplicative mean
ergodic theorem given in (\ref{eq:mainMMET})
will play a central role in the proofs
of all the subsequent probabilistic limit theorems.

\begin{theorem}
\tlabel{mainSpectral}
{\em (Multiplicative Mean Ergodic Theorem)}
Assume that the Markov chain $\bfPhi$ 
and the functional $F$ satisfy \eq assumption/.
With $\delta$ and $b$ as in (V4), define:
\begin{equation}
\bara\eqdef  \Bigl( \frac{e-1}{2b - \delta}\Bigr) \delta  >0.
\elabel{bara}
\end{equation}
Then there exists $\baromega>0$ 
such that,
for any $\alpha$ in the compact set
$$
\Omega=\{\alpha=a+i\omega\in\Co\,:\,
|a|\leq\bara,\;\mbox{and}\;|\omega|\leq\baromega\},
$$
there is an eigenvalue $\lambda_\alpha \in \clS_\alpha$
which is maximal and isolated, i.e.,
\[ 
|\lambda_\alpha|=
\max \{|\lambda|:\lambda\in\clS_\alpha\}
\quad
\hbox{and}
\quad
\clS_\alpha  
\cap \bigl
\{ z : |z| \ge |\lambda_\alpha |-\delta_0 \bigr\}
=\{\lambda_\alpha\}
\]
for some $\delta_0>0$.

Moreover, for any such $\alpha$,  
there exist $\cf_\alpha \in L_\infty^V$
and $\cmu_\alpha \in \clM_1^V$, satisfying 
\eq cmu-norm/, and:
\begin{description}
\item[(i)]  
The functions $\cf_\alpha$ solve the multiplicative
Poisson equation
$$\haP_\alpha \cf_\alpha = \lambda_\alpha \cf_\alpha,$$
and the $\cmu_\alpha$ 
are eigenmeasures for 
the kernels $\haP_\alpha$:
$$\cmu_\alpha \haP_\alpha  =  \lambda_\alpha \cmu  _\alpha.$$
        \item[(ii)] 
There exist constants
$b_0>0$, $ B_0 < \infty$, 
such that for all $\alpha\in\Omega$,
$x\in\state$,  $n\ge 1$,
\be
\Bigl| \Expect_x \bigl[ \exp (\alpha S_n - n \LA(\alpha))\bigr]
                - \cf_\alpha(x)\Bigr|
\leq 
                        B_0 |\alpha| V(x) e^{-b_0 n} \,,
                \label{eq:mainMMET}
\ee
where 
$\LA(\alpha)\eqdef\log(\lambda_\alpha)$
is analytic on $\Omega$, and $S_n$ are the partial 
sums defined in \eq discretePS/. More generally,
for any $g\in L_\infty^V$,  
$$
\Bigl| \Expect_x \bigl[ \exp (\alpha S_n 
                - n \LA(\alpha)) g (\Phi(n))\bigl] 
                - \haQ_\alpha g \, (x)\Bigr| 
\leq 
                B_0 \| g \|_{V} V(x) e^{-b_0 n} .$$
\end{description}
\end{theorem}

\proof
The existence of an isolated, maximal eigenvalue
$\lambda_\alpha$ is given in \Proposition{basic1}. 
It is nonzero for $\alpha=a\in[-\bara,\bara]$ by
\Proposition{finiteSR}, and since it is analytic
in $\alpha$ (by \Proposition{basic1}), we can 
pick $\baromega>0$ small enough such that 
$\lambda_\alpha$ is nonzero on $\Omega$.

The existence of an eigenfunction and eigenmeasure
as in (i) follows from   
\Proposition{basic1} combined with \Proposition{V-geo}.  
To see that $\bara>0$ note that, under (V4),
\be
\pi(V)=\pi(PV)\le (1-\delta)\pi(V) + b\pi(s).
\label{eq:bbound}
\ee
Hence, $b\ge \delta\pi(V)/\pi(s)\ge \delta$.

To prove the limit theorems in (ii), 
consider the linear operator
\[
U_{(z,\alpha)} =[ Iz - (\lambda_\alpha^{-1} 
        \haP_\alpha -\haQ_\alpha )]^{-1}.
\]
From  \Proposition{V-geo} we can find
$\epsy_0>0$ such that $U_{(z,\alpha)}$ is an analytic function of 
two variables $(z,\alpha)=(z,a+i\omega)$  
on the domain 
\[
\clD=
\bigl\{ |z| > 1-\epsy_0  ,\  |a| < \bara+\epsy_0,\ |\omega| 
                           < \baromega+ \epsy_0\bigr\}.
\]
We may also assume that $\epsy_0>0$ is suitably small 
so that, for some $\barb<\infty$, we have
$\lll U_{(z,\alpha)}\lll_{V}\le \barb$  
for all $(z,\alpha)\in\clD$.

Set $b_0= -\log(1-\epsy_0) > 0$.  The following bound then holds 
for all $(z,\alpha)\in\clD$, $g\in L_\infty^V$, $x\in\state$, and $n\ge 1$,
by representing  $U_{(z,\alpha)}$ as a power-series, and using
the fact that $\haQ_\alpha$ is a projection operator:
\[
\begin{array}{rcl}
\barb V(x)
&\ge&
\Bigl| \int_0^{2\pi}
e^{i(n+1)\phi} U_{( \exp(-b_0+i n\phi), \alpha)} g\,(x)
        \, d\phi \Bigr| 
\\[.3cm]
&=&
e^{(n+1)b_0}  \bigl |(\lambda_\alpha^{-1} \haP_\alpha  
-\haQ_\alpha )^n g\, (x) \bigr| 
\\[.2cm]
&=&
e^{(n+1)b_0}  \bigl |\lambda_\alpha^{-n} 
\haP_\alpha^n g\, (x)  -\haQ_\alpha g\, (x) \bigr|.
\end{array}
\]
This gives the second bound in (ii).  
The first one follows from the second since, 
when $\alpha=0$ and $g=\one$,
\[
\bigl |\lambda_\alpha^{-n} \haP_\alpha^n g\, (x)  -\haQ_\alpha g\, (x) \bigr| 
= 0, 
\]
for all $n\ge 1,\ x\in\state$.
\qed

\medskip

Next we give a weaker multiplicative mean
ergodic theorem for all $\alpha=a+i\omega$ in
a neighborhood of the $i\omega$-axis.
A function $F:\state\to\Re$ is called 
{\em lattice} if there are
$h>0$ and $0\le d < h$, such that
\be
\frac{F(x)-d}{h}\qquad\hbox{is an integer,}
\qquad x\in\state\, .
\elabel{lattice}
\ee
The minimal $h$ for which
\eq lattice/ holds is called the {\em span} of $F$.
If the function $F$ can be written as a sum,
\[
F=F_0 + F_\ell,
\]
where $F_\ell$ is lattice with span $h$
and $F_0$ has zero
asymptotic variance (recall \eq var1/),
then $F$ is called {\em almost-lattice} (and $h$
is its span).
Otherwise, $F$ is called {\em strongly non-lattice.}

Although these definitions are somewhat different from 
the ones commonly used when studying the partial sums
of independent random variables, in the Markov case
they lead to the natural analog of the classical
lattice/non-lattice dichotomy. This dichotomy,
which is close in spirit to the discussion in 
\cite{shurenkov:84}, is stated in \Theorem{lattice}.

\medskip

\begin{theorem}
{\em (Bounds Around the $i\omega$-Axis)}
\tlabel{mainSpectral2}
Assume that the Markov chain
$\bfPhi$ and the functional $F$
satisfy \eq assumption/.
\begin{itemize}
\item[(NL)]
Suppose that $F$ is strongly non-lattice.  
%%% YIANNIS: note to myself: keep this statement here:
%%% $\Ree (\LA(i\omega))<0$ for all nonzero $\omega\in\RL.$ 
%%% in case i DO need it later in the proofs
For any $0<\omega_0 <\omega_1< \infty$, 
there exist $b_0>0$, $B_0 < \infty$ 
(possibly different than in \Theorem{mainSpectral}), 
such that
\be
\Bigl | \Expect_x [ \exp (\alpha S_n - n\LA(a))] \Bigr|
\leq 
B_0 V(x) e^{-b_0 n}\, ,  \qquad  x \in \state,\; n \geq 1,
\label{eq:ImMMET}
\ee
for all $\alpha=a+i\omega$ with $|a| \le \bara$
and $\omega_0\leq |\omega|\leq\omega_1$.

\item[(L)] 
Suppose that $F$ is almost-lattice with span $h>0$.
For any $\epsilon>0$, there exist  
$b_0>0$, $B_0 < \infty$ 
(possibly different than above and
in \Theorem{mainSpectral}), 
such that (\ref{eq:ImMMET}) holds
for all $\alpha=a+i\omega$ with $|a| \le \bara$
and $\epsilon\leq|\omega|\leq 2\pi/h - \epsilon.$
\end{itemize}
\end{theorem}

\proof
By \Theorem{lattice} we have the bound 
$\haxi_\alpha < \haxi_a=\lambda_a$
for the range of $\alpha\in\Co$ considered in the theorem. 
This implies that there is an $\epsy_1>0$, $b_1<\infty$ such that 
\[
\lll [Iz - e^{-\lambda_a} \haP_\alpha ]^{-1}\lll_V < b_1
\]
for all $|z|\ge 1-\epsy_1$, and all $\alpha$ in this range.  
An argument similar to the proof of \Theorem{mainSpectral}~(ii) 
then gives the desired bounds.
\qed

\subsection{Spectral Radius and Spectral Gap}
%%%%%%%%%%%%%%%%%%%%%%%%%% Spectral Gap %%%%%%%%%%%%%%%%%%%%%%%%%%% 

Recalling our standing assumption \eq assumption/, 
we fix the Lyapunov function
$V:\state\to[1,\infty)$ throughout this section.

For complex $\alpha$ we wish to construct  $\LA(\alpha)\in\Co$ 
satisfying the multiplicative mean ergodic limit,
\[
\LA(\alpha)  =  
\lim_{n \to \infty} \smalloneOvern\! \log 
\Expect_x [ \exp (\alpha S_n ) ], \qquad x\in\state.
\]
This requires a generalization of the notion of the 
g.p.e.\ of \Section{gpe}. 
The  
previous definition is meaningless when 
$\alpha\not\in\Re$, since the definition of a 
small set depends on the linear ordering of $\RL$.

\paragraph{Spectral radius.}
For a bounded linear operator
$\haP\colon L_\infty^V\to L_\infty^V$
we define the \textit{spectral radius} 
of $\haP$ by
\begin{equation}
\haxi \eqdef \lim_{n\to\infty}
 (\lll \haP^n \lll_{V} )^{1/n}
=
\exp\Bigl( \lim_{n\to\infty}
\smalloneOvern\! \log\lll \haP^n \lll_{V} \Bigr).
\elabel{SpectralRadius}
\end{equation} 
Note that in the above definition 
$\haP$ is not assumed to be a positive
operator, and it is possibly 
complex-valued.
Since $\lll \cdot\lll_V$ is an operator norm,
the sequence $\{ \log(\lll \haP^n \lll_V) \,:\, n\geq 1\}$ 
is subadditive \cite{rienag55}.
Therefore $\haxi$ always 
exists, although it may be infinite.

We let $\haxi_\alpha$ denote
the spectral radius of the operator $\haP_\alpha$ 
defined in \eq haP-discrete/.  
When $\alpha=a$ is real, from the definitions we have 
that $\haxi_\alpha\ge \lambda_\alpha$
where $\lambda_\alpha$ 
is the g.p.e.\ of the positive
kernel $\haP_\alpha$.
One of the main goals of this section is to
show that the spectral radius 
$\haxi_\alpha$ coincides with $\lambda_\alpha$ 
for real $\alpha$ in a neighborhood of $\alpha =0$.
We first establish upper and lower bounds:

\begin{proposition}
\tlabel{finiteSR}
Under \eq assumption/, the spectral 
radius $\haxi_\alpha$ of $\haP_\alpha$
is finite and
\[
\haxi_\alpha\le (b +1)\exp(|a| )\, ,
\]
for all $\alpha=a+i\omega\in\Co$. 
Moreover, for $\alpha=a\in\RL$,
$$\haxi_a\geq e^{-a}>0.$$ 
\end{proposition}

\proof
The function $s$ in (V4) is necessarily bounded by one.  Consequently, 
under (V4) we have for any $g\in L_\infty^V$,  $\alpha = a+ i\omega\in \Co$,
\[
|\haP_\alpha g\, (x)|  \le \exp(|a|) \|g\|_V PV 
\le  \exp(|a|) \|g\|_V (1 + b) V \, .
\]
This implies that $\lll \haP_\alpha \lll_V \le e^{|a|}(1+b)$.  
The operator norm $\lll\varble\lll_V$ is submultiplicative,
\[
\lll \haP_\alpha^n \lll_V \le
\lll \haP_\alpha \lll_V^n \le e^{|a|n}(1+b)^n, \qquad n\ge 1,
\]
giving the upper bound.

When $\alpha=a$ is real,
for any $g\in L_\infty^V$, $g\ge 0$, we have,
\[
\haP_a^n g\, (x) \ge e^{-an} P^n g(x).
\]
It follows immediately that
$\haxi_a \ge \haxi_0 e^{-a} = e^{-a}.$
\qed

\paragraph{Spectral gap and $V$-uniform operators.}

Recall the following 
classical result from \cite[p.~421]{rienag55}:

\begin{theorem}
\tlabel{decomposition}
{\em (Decomposition Theorem)}
Let $\haP:\LV\to\LV$ be a bounded
linear operator, and suppose that $z_0\in \clS(\haP)$ is 
isolated, i.e., for some $\epsy_0>0$,
\[
\clS(\haP)\cap D = \{z_0\}
\qquad\hbox{\it where} \qquad D=\{ z\in\Co : |z-z_0|\le \epsy_0\}\, .
\]
Then, the following bounded operator 
on $L_\infty^V$ is well-defined,
\[
\haQ = \frac{1}{2\pi i} \int_{\partial D} [Iz-\haP]^{-1} \, dz
\,,
\]
and moreover:
\begin{description}
\item[(i)]
$\haQ\colon L_\infty^V\to  L_\infty^V$ is a projection operator,  
that is, $\haQ^2 =  \haQ$;
\item[(ii)]
$\haP\haQ = \haQ\haP =z_0\haQ$;

\item[(iii)]
 $\clS(\haQ)=\{1\}$, and
 $\clS(\haP-z_0\haQ)\cap D =\emptyset$.
\end{description}
\end{theorem}

We say that $z_0\in \clS(\haP)$ is a 
\textit{pole of finite multiplicity} if
$z_0$ is an isolated point in $\clS(\haP)$
and the associated projection operator $\haQ$ 
can be expressed as a finite linear combination 
of some $\{s_i\}\subset L_\infty^V$, 
$\{\nu_i\}\subset \clM_1^V$:
\be
\haQ = \sum_{i,j=0}^{n-1} m_{i,j} [ s_i\otimes\nu_j].
\label{eq:poles}
\ee
In particular, we call $z_0$ 
a {\em pole of multiplicity one}, 
if (\ref{eq:poles}) holds for $n=1$,
and also there exists $\epsy_0>0$ such that
\[
\clS(\haP- z_0 \haQ ) \subset
\{ z : |z|\le \haxi - \epsy_0 \} \,,
\]
where $\haxi$ is the spectral radius of $\haP$.

We say that $\haP$ admits a \textit{spectral gap} 
if there exists $\epsy_0>0$ 
such that $\clS(\haP)\cap \{ z : |z|\ge \haxi - \epsy_0 \}$ is finite, 
and contains only poles of finite multiplicity.  

Further, we say that $\haP$ is 
{\em $V$-uniform},
if it admits a spectral gap
and also there exists a unique 
pole $\lambda_\circ\in \clS(\haP)$
of multiplicity one, 
satisfying $|\lambda_\circ| = \haxi$.  
In that case, $\lambda_\circ$ is called 
{\em the generalized principal eigenvalue}
(g.p.e.),
generalizing the previous definition.
In particular, if $\haP_\alpha$ is
$V$-uniform for some $\alpha\in\Co$,
then we write $\lambda_\alpha$ for
its associated g.p.e.

Much of the  
development of this section, is based on 
properties of rank-one operators 
of the form $\haM=s_0\otimes\nu_0$  
for some $s_0\in L_\infty^V$, 
$\nu_0\in \clM_1^V$.  
The associated {\em potential operator} 
is defined as
\begin{equation}
\haU_z 
\eqdef \Bigl[Iz - (\haP - \haM )\Bigr]^{-1} ,
                                        \qquad z\in \Co\,,
\elabel{U}
\end{equation}
whenever the inverse exists.
The potential operator
is used to construct eigenfunctions and 
eigenmeasures for a $V$-uniform operator:

\begin{proposition} 
\tlabel{eig}  
Suppose that $\haP$ is $V$-uniform with g.p.e.\ $\lambda_\circ$, 
and that the associated $s_0,\nu_0$ 
in (\ref{eq:poles}) are chosen so that the 
potential operator $\haU_z$ in \eq U/ 
is bounded for $z$ in a neighborhood of $|z|\ge|\lambda_\circ|$.  
Then, setting $\cf = \haU_{\lambda_\circ} s_0$ 
and $\cmu =  \nu_0 \haU_{\lambda_\circ}$, 
we have $\cf\in L_\infty^V$,  $\cmu\in \clM_1^V$,
\[
\haP \cf = \lambda_\circ \cf\,  ,
 \quad
\mbox{\it and}
\quad
\cmu \haP  = \lambda_\circ \cmu \,  .
\]
\end{proposition}

\proof 
From $V$-uniformity we know that there exists $\epsy_0>0$ such that the
inverse    $\bigl[Iz -\haP\bigr]^{-1}$  exists and is bounded as a linear
operator on $\LV$, for all $|z|\ge \haxi- \epsy_0$, $z\neq \lambda_\circ$.
Moreover, for such $z$ we may apply the inversion formula \eq inversion/
to obtain the identity, 
\begin{equation}
\bigl[Iz -\haP\bigr]^{-1}
=
\haU_z +   \frac{ (\haU_z s_0) \otimes (\nu_0 \haU_z) }
                 {1-\nu_0 \haU_z s_0} .
\elabel{Uinv}
\end{equation}
Since $\lambda_\circ\in \clS(\haP)$, and
$\lambda_\circ\not\in\clS(\haP-s_0\otimes\nu_0)$, it 
follows from this
equation that $\nu_0 \haU_{\lambda_\circ} s_0 =1$.

Applying $[Iz-\haP]$ to \eq Uinv/ on the left, and $s_0$ on the right
then gives,
\[
s_0 = [Iz - \haP] \haU_z s_0  +  [Iz - \haP] \frac{ (\haU_z s_0)   (\nu_0
\haU_z s_0) }
                 {1-\nu_0 \haU_z s_0} .
\]
Multiplying both sides by $(1-\nu_0 \haU_z s_0) $, and then setting
$z=\lambda_\circ$ gives $0=[I\lambda_{\circ} - \haP]
\haU_{\lambda_\circ} s_0$, which shows
that $\cf$ is an eigenfunction.

The proof that $\cmu$ is an eigenmeasure is completely
analogous, and follows by applying $[Iz-\haP]$ to \eq Uinv/ 
on the right and $\nu_0$ on the left.
\qed

\medskip

The following proposition 
provides useful 
characterizations of $V$-uniformity.

\begin{proposition} 
\tlabel{equivalents}
The following are equivalent for an operator 
$\haP$ with finite spectral
radius $\haxi$.
\begin{description}
\item[(i)]
$\haP$ is $V$-uniform.
\item[(ii)]
There exists $\lambda\in\Co$ satisfying $|\lambda|=\haxi$, and 
$\epsy_0>0$, $s_0\in L_\infty^V$, $\nu_0\in \clM_1^V$, satisfying
\[
\begin{array}{rcl}
\sup\Bigl\{ 
	|z-\lambda| \lll Iz - \haP\lll_V : |z|\ge \haxi-\epsy_0 
\Bigr\}& <& \infty,
\\[.25cm]
\mbox{and}\qquad
\sup\Bigl\{  \lll \haU_z \lll_V  : |z|\ge \haxi-\epsy_0 
\Bigr\} &<& \infty,
\end{array}
\]
where $\haU_z$ is the potential operator defined in \eq U/,
with $\haM=s_0\otimes\nu_0$.

\item[(iii)]
There exists $\lambda\in\Co$ satisfying $|\lambda|=\haxi$ and
$\cf\in L_\infty^V$, $\cmu\in \clM_1^V$, such that
\[
\lambda^{-n} \haP^n \to \cf\otimes\cmu,\qquad n\to\infty,
\] 
where the convergence is in the $V$-norm.
\end{description}
\end{proposition}

\proof
If (i) holds then the matrix-inversion formula \eq inversion/ gives
\[
[Iz - \haP + \lambda_\circ \cf\otimes\cmu]^{-1}
=
[Iz-\haP]^{-1} - \frac{\lambda_\circ}{z(z-\lambda_\circ)} \cf\otimes\cmu.
\]
The left hand side is bounded for $|z|\ge \haxi - \epsy_0$ under (i).  
Hence we may set $\cf=s_0$ and $\cmu=\nu_0$ to obtain (ii).

The implication (ii) $\Rightarrow$ (i) also  
follows from the matrix inversion formula \eq inversion/ since \eq Uinv/
then holds for all $|z|\ge \haxi- \epsy_0$, $z\neq \lambda$.
This implies the limit
\begin{equation}
\haQ 
\eqdef  
\lim_{z\to\lambda} (z-\lambda)\bigl[Iz -\haP\bigr]^{-1}
=
\frac{ (\haU_\la s_0) \otimes (\nu_0 \haU_{\lambda}) }
           {      \nu_0 \haU_{\lambda}^2 s_0  } \,,
\elabel{QhatFormula}
\end{equation}
and (i) holds with this $\haQ$,  and $\lambda_\circ=\lambda$.

The equivalence of (i) and (iii) follows 
exactly as in \Theorem{mainSpectral}~(ii).
\qed

\medskip

For a probabilistic kernel $\haP$, 
the following proposition says that
$V$-uniformity implies that 
the chain with transition kernel $\haP$
is geometrically ergodic. 
The converse is also true; 
see \Proposition{Ubdd}.

\begin{corollary}
\tlabel{Vuniformity}
If $\haP$ is a $V$-uniform, probabilistic kernel, 
then
the Markov chain with transition kernel $\haP$
is geometrically ergodic.
\end{corollary}

\proof
Since $\haP$ is probabilistic,
applying the limit result
of \Proposition{equivalents}~(iii)
to the constant function $\bf 1$,
implies that $\la=1$ and that $\cf$
is constant. By rescaling we can
take $\cf={\bf 1}$ and $\cmu$ to
be a probability measure. From 
\Theorem{Vuniform} it the follows
that $\haP$ is geometrically
ergodic.
\qed

\medskip

\Proposition{eig} applied to the family of kernels 
$\{\haP_\alpha\}$ gives the following:

\begin{proposition}
\tlabel{V-geo}
Suppose that $\haP_{\alpha_0}$ is $V$-uniform for 
a given $\alpha_0\in\Co$. Then there exists $\epsy_0>0$ 
such that $\haP_\alpha$ is $V$-uniform (with associated
g.p.e.\ $\la_\alpha$) for all $\alpha\in\Co$,
$|\alpha-\alpha_0| < \epsy_0$. Moreover, for each
such $\al$ there exist $\cf_\alpha  \in L_\infty^V$ and
$\cmu_\alpha\in \clM_1^V $ such that:
\begin{description}
\item[(i)]
$\cf_\alpha$ solves the multiplicative 
Poisson equation,
$\,\haP_\alpha \cf_\alpha =\lambda_\alpha \cf_\alpha$.

\item[(ii)]
$\cmu_\alpha$ is an eigenmeasure for $\haP_\alpha$,
$\cmu_\alpha \haP_\alpha =\lambda_\alpha \cmu_\alpha$.
 
\item[(iii)]
The g.p.e.\ $\lambda_\alpha$ 
is an analytic function of 
$\alpha$,
and so is $\cf_\alpha(x)$ for any fixed
$x\in\state$.
\end{description}
\end{proposition}
 
\proof
The existence of eigenvectors in (i) and (ii)
is immediate from \Proposition{eig} when $\alpha = \alpha_0$. 
Define  $\haU_{z}=\haU_{z,\alpha}$ by 
\eq U/ with $\haP= \haP_\alpha$, and 
$\haM=s_0\otimes \nu_0 $:
\be
\haU_{z,\alpha}
 \eqdef \Bigl[Iz - (\haP_{\alpha} - s_0\otimes \nu_0 )\Bigr]^{-1}\, .
\label{e:Uzalpha}
\ee
From $V$-uniformity we know that $\haM$ can be chosen so that 
$\haU_{z,\alpha_0}$ is a bounded
linear operator for $z$ in a neighborhood of $\lambda_{\alpha_0}$.
Since $\haP_\alpha$  is continuous in $V$-norm, it then follows that 
$\haU_{z,\alpha}$ is a bounded linear operator for $(z,\alpha)$ in a 
neighborhood $O$ of $(\lambda_{\alpha_0},\alpha_0)$.  
This combined with \Proposition{eig} proves (i) and (ii).

Write  $\clJ(z,\alpha) = \nu_0(\haU_{z,\alpha} s_0)$, $z\in \Co, \alpha\in O$, 
so that
\[
\begin{array}{rcll}
\clJ(\lambda_\alpha,\alpha) &=&1,\qquad &
\\
\frac{\partial}{\partial z}\clJ(z,\alpha)\Big|_{z=\lambda_\alpha} 
                  &=&   \nu_0(\haU_{z,\alpha}^2) s_0
                   = \cmu_\alpha(\cf_\alpha)
                            \neq 0
                           ,\qquad &\alpha\in O\, ,
\end{array}
\] 
where $\cf_\alpha,\cmu_\alpha$ are the eigenfunction and eigenmeasure 
given in \Proposition{eig}.
We conclude that $\lambda_\alpha$ is an analytic function by 
the implicit function theorem. 

The proof that $\cf_\alpha(x)$ is analytic 
in $\alpha$ for $x\in S_V=\state$ 
follows from the expansion
\begin{equation}
\cf_\alpha = \haU_{\lambda_\alpha,\alpha}s_0 
= 
\sum_{n=0}^\infty \lambda_\alpha^{-n-1} (\haP_{\alpha} - s_0\otimes \nu_0 )^n s_0 .
\elabel{cf-unnorm}
\end{equation}
This expression for $\cf_\alpha$ converges uniformly 
for $\alpha\in O$, and for each $n$ the finite 
sum is analytic, which completes the proof of (iii).
\qed

\medskip

The eigenfunction \eq cf-unnorm/ will not in general satisfy the 
required normalization \eq cmu-norm/.  The following eigenfunction and 
eigenmeasure do satisfy this condition, and are the 
\textit{unique} such solutions,
\begin{equation}
\cmu_\alpha 
= 
\frac{ \nu_0\haU_{\lambda_\alpha,\alpha}}{\nu_0 
\haU_{\lambda_\alpha,\alpha}\One} 
\,
\in \clM_1^V
\qquad\qquad
\cf_\alpha =
 \frac{\haU_{\lambda_\alpha,\alpha}s_0}{\cmu_\alpha 
\haU_{\lambda_\alpha,\alpha}s_0} 
\,
 \in L_\infty^V\,.
\elabel{cf-norm}
\end{equation}
Given such $\cf_a$ and $\lambda_a$ for some real $a$,
we define the {\em twisted kernel} $\cP_a$ by
\[
\cP_a (x, dy) = \lambda_a^{-1} \cf_a^{-1}
(x) \haP_a (x, dy) \cf_a (y)\,  ,
\]
(cf.\ \eq cP/ in the introduction), and we 
let $\cV_a= V/\cf_a$. (As we will see below,
$\cf_a$ is bounded away from zero for real
$a$ in the range of interest.)
The following proposition describes
the relationship between the transition 
kernels, the eigenfunctions, and the eigenmeasures
$\{\cP_a, \cf_a,\cmu_a :a\in\Re\}$.

\begin{proposition}
\tlabel{V-geo2}
Suppose that $\haP_{a_0}$ is $V$-uniform for
a given real $a_0$. Then there is an open 
set $O\subset\RL$ containing $a_0$,
such that, for all $a\in O$,
with $\cf_a,\cmu_a$ given in \eq cf-norm/
and with $\cP_a$ equal to the associated twisted kernel,
we have:
\begin{description}
\item[(i)]
The operator
$\cP_a$ is $\cV_a$-uniform.

\item[(ii)]
$\displaystyle \frac{d}{da} \LA(a) 
= \frac{d}{da} \log(\lambda_a)
=\cpi_a(F)$, where $\pi_a$ is the 
invariant probability
measure for $\cP_a$.

\item[(iii)]
 $\displaystyle 
\haF_a\eqdef\frac{d}{da} \log(\cf_a)$ 
is a solution to the Poisson equation,
\begin{equation}
\cP_a \haF_a = \haF_a - F + \pi_a(F)\, .
\elabel{PoissonCheck}
\end{equation}
For $a=0$, this is the unique solution satisfying
$\pi(\haF)=0$.

\item[(iv)]
 $\displaystyle \frac{d}{da} \cf_a \in L_\infty^V,\qquad
 \frac{d}{da} \cmu_a \in \clM_1^V$.

\item[(v)]
 $\displaystyle 
\haF_a \in L_\infty^{1+\log(V_a)}$.

\end{description}
\end{proposition}

\proof
The existence of $O$ follows from \Proposition{V-geo},
and from its proof we know that
$\lll \haU_{\lambda_a,a}\lll_V<\infty$ when $a\in O$, 
where $\haU_{z,a}$ is given 
in \eq Uzalpha/.

The linear operators $\haP_a$ and $\cP_a$ 
are related by the scaling $\lambda_a$ 
and a similarity transformation,
\[
\cP_a  = \lambda_a^{-1} (I_{\cf_a})^{-1} 
\haP_a I_{\cf_a} \, ,
\]
where $I_g$, for an arbitrary function $g$,
denotes the kernel $I_g(x,\cdot)\eqdef g(x)\delta_x(\cdot)$.
Hence $\haP_a$ is $V$-uniform if and only if $\cP_a$ is 
$(I_{\cf_a}^{-1} V)$-uniform.  Result (i)  immediately follows.

Consider the unnormalized eigenfunction given in \eq cf-unnorm/.
Differentiating the expression $
\cf_a = \haU_{\lambda_a,a}s_0 $
and applying the quotient rule gives,
\begin{equation}
\begin{array}{rcl}
\cf_a' &=& \frac{d}{da}   
[ I \lambda_a -  \haP_a + s_0\otimes\nu_0 ]^{-1} s_0
\\[.2cm]
&=&
- \haU_{\lambda_a,a} [ I \lambda_a' - I_F \haP_a ] \haU_{\lambda_a,a} s_0
\\[.2cm]
&=&
-  \haU_{\lambda_a,a}  [ ( \lambda_a' 
              - \lambda_a F )\cf_a ]
=
\lambda_a \haU_{\lambda_a,a} I_{\cf_a } [ F-  \LA'(a) ]\, .
\end{array}
\elabel{fishBdd1}
\end{equation}
The right hand side of \eq fishBdd1/ 
lies in $ L_\infty^V$ since $F\in L_\infty$, $\cf_a\in L_\infty^V$, 
and  $\haU_{\lambda_a,a}\colon L_\infty^V\to L_\infty^V$ is a bounded linear operator.  
This proves the first bound in (iv) since the two versions of $\cf_a$ 
are related by a smooth normalization.  
The proof that $ \frac{d}{da} \cmu_a \in \clM_1^V$ is identical.

Differentiating both sides of the eigenfunction equation gives
\[
F \lambda_a \cf_a + \haP_a \cf'_a 
   =  \lambda_a' \cf_a +\lambda_a  \cf_a'  .
\]
Dividing this identity by $\lambda_a\cf_a $ 
shows that \eq PoissonCheck/ does indeed hold.  To conclude that 
$\cpi_a(F)=\LA'(a)$ we will show that 
$\cpi_a(|\haF_a|)<\infty$.  The invariant probability measure
$\cpi_a$ may be expressed as
\[
\cpi_a = k_a \cmu_a I_{\cf_a},
\]
where $k_a$ is a normalizing constant.
Hence,
\[
\cpi_a(\haF_a)
=
\cpi_a\Bigl(\Bigl|\frac{\cf_a'}{\cf_a}\Bigr|\Bigr) 
= 
k_a \cmu_a(|\cf_a'|)<\infty\, .
\]
Finiteness follows from (iv) and
the fact that the eigenmeasure $\cmu_a$ lies in $\clM_1^V$.  
This proves (ii) and the identity in (iii).

To complete the proof of (iii) we 
must show that $\pi(\haF_0)=0$.  This follows
from the normalization \eq cmu-norm/ (assumed to hold for all $a$) which 
implies the limits,
\[
\cf_a \to \One,\quad \cmu_a \to \pi,\qquad a\to 0.
\]

To prove (v) we obtain an alternative expression for $\haF_a$.  We again
consider the unnormalized eigenfunction \eq cf-unnorm/.
Observe that a fundamental kernel is derived from $\haU_{\lambda_a,a}$ 
through a scaling and a similarity transformation,
\[
Z_{a} =
\lambda_a I_{\cf_a}^{-1} \haU_{\lambda_a,a} I_{\cf_a}
        =
      [ I -  \cP_a + s_a \otimes\nu_a ]^{-1},
\]
with $s_a = \lambda_a^{-1} \cf_a^{-1} s_0$, and
$\nu_a = \nu_0 I_{ \cf_a }$.
We have $\cP_a Z_a G = Z_a G - G$ whenever $\pi_a(G)=0$.

Using \eq fishBdd1/ then gives,
\[
\haF_a=
\frac{\cf_a'}{\cf_a}
= Z_a (F - \LA'(a)) = Z_a (F - \cpi_a(F)) .
\]
It again follows that $\haF_a$ solves the Poisson equation: 
It is the unique solution in $\LV$ with $\nu_a(\haF_{a})=0$.

The desired bound on $\haF_\alpha$ is obtained as follows.  Using Jensen's 
inequality we know that $\haV_\alpha=\log(V/\cf_\alpha)=\log(\cV_\alpha)$ 
solves a  version of (V3),
\[
\cP_\alpha \haV_\alpha\le \haV_\alpha -\epsy + bs,
\]
where $\epsy>0$ and $b$ is a finite constant.  Using the bound 
$\cf_\alpha\in L_\infty^V$ it follows directly that 
the function $\haV_\alpha$ is uniformly bounded below.
The bound on $\haF_\alpha$ then follows from 
\cite[Theorem 2.3]{glymey96}.
\qed

\medskip

\begin{proposition}
\tlabel{Ubdd}
Suppose that \eq assumption/ holds. 
Take $\haP=R$, and define the potential
operator $\haU_z$ as in \eq U/ 
with $\haM=s\otimes\nu$. Then 
$B_1 \eqdef \lll \haU_1\lll_V \le 2b\delta^{-1}$, and 
\[ 
\lll \haU_z \lll_V \le B_1(1 - |z-1|B_1)^{-1} ,\qquad |z-1|\le B_1^{-1} ,
\qquad z\in\Co.
\]
Hence both $R$ and $P$ are $V$-uniform.
\end{proposition}

\proof
Under (V4) we have, by \eq resolve/,
\[
(e-1)(R-I)V = R(P-I) V \le -\delta RV + bRs .
\]
Rearranging terms then gives
\[
(R-I)V  \le -\Bigl(\frac{\delta}{e-1+\delta}\Bigr) V 
                       + b \Bigl(\frac{1}{e-1+\delta}\Bigr)Rs,
\]
which we write as
\[
(R-s\otimes\nu) V \le V - \delta_1 V 
                 -\nu(V)s  + b_1 Rs\,,
\]
where $\delta_1 =  \delta(e-1+\delta)^{-1}$ and $b_1 = b (e-1+\delta)^{-1}$.

Iterating gives, for all $n\ge 1$,
\begin{eqnarray*}
(R - s\otimes \nu)^nV &\le& V - \delta_1 
        \sum_{i=0}^{n-1} (R - s\otimes \nu)^i V
\\
&&\quad + \sum_{i=0}^{n-1} 
        (R - s\otimes \nu)^i (b_1 Rs - \nu(V) s).
\end{eqnarray*} 
Letting $n\to\infty$, and applying \Lemma{small}~(ii) yields
\[
\delta_1 \haU_1V 
\le  V -\nu(V) + 2b_1
\le 2b_1 V,
\]
or $\lll \haU_1\lll_V\le 2b_1/\delta_1= 2b/\delta$.

To obtain a bound for $z\sim 1$ write 
\[
\haU_z=\Bigl[I(z-1) + [I-(R-s\otimes\nu)] \Bigr]^{-1}
=
\Bigl[\haU_1(z-1) + I \Bigr]^{-1}\haU_1\, .
\]
Provided $\lll \haU_1 \lll_V |z-1| <1$, we can write
$\haU_z = \sum_n (1-z)^n \haU_1^{n+1},$
and
$\lll \haU_z\lll_V \le B_1/[1-|z-1| B_1].$
\qed

\begin{proposition}
\tlabel{Hcriterion}
Suppose that \eq assumption/ holds, let 
$a\in\Re$ satisfy 
$|a|\le|\log(1-\delta)|$,  
and suppose that there exists $
g\colon\state\to (0,\infty)$, satisfying
$g\in L_\infty^V$ and $\haP_a g\le \lambda_a g$.    
Then $\haP_a$ is $V$-uniform.
\end{proposition}

\proof
The conditions of the proposition imply that there exists $b_1<\infty$
such that
\[
\haP_a V \le e^{|a|} (1-\delta) V +   e^{|a|} b s
                     \le   V + b_1 s.
\]
From the resolvent equation \eq resolve/ we then have, for some $b_2<\infty$,
\[
\haR_\theta V \le V +b_2 \haR_\theta s,
\]
where $\haR_\theta$ is the resolvent kernel defined through $\haP_a$.

We also have $\LA(a) > 0$ for all $a\neq 0$ under \eq assumption/,
and hence the g.p.e.\  $\gamma_\theta $ 
for $\haR_\theta$ is also strictly greater than one 
when $\theta > |a| >0 $ (see \eq evaluesRes/).  
Choosing $s_0 \in \clB_p^+$ and  $\nu_0\in \clM_p^+$ so that
$R_\theta \ge s_0\otimes\nu_0$, we find that 
\[
\gamma_\theta^{-1}(\haR_\theta -s_0\otimes\nu_0) V \le V 
        - \epsy V +b_2 \haR_\theta s_0,
\]
where $\epsy = 1-\gamma_\theta^{-1}>0$.
Exactly as in the proof of \Proposition{Ubdd} we conclude that
\[
\epsy \haU_{\gamma_\theta} V \le b_2 \haU_{\gamma_\theta} \haR_\theta s_0 
\le 2b_2 \haU_{\gamma_\theta} s_0 ,
\]
where $\haU_z = \sum z^{-n-1}(\haR_\theta -s_0\otimes\nu_0)^n$.

From \Theorem{evector}~(iii) and the conditions of the
proposition we know that 
$\cf = \haU_{\gamma_\theta} s_0 $ satisfies
$\cf \le c g$ for some constant $c$, and hence $\cf\in L_\infty^V$.  
It follows
that $ \lll\haU_{\gamma_\theta}\lll_V<\infty$,
from which $V$-uniformity of 
$\haR_\theta$, and hence of $\haP_a$, immediately
follow.
\qed

\begin{proposition} 
\tlabel{basic1}
Suppose that \eq assumption/ holds.
Then there exists $\epsy_0>0$, $\barb<\infty$ such that:
\begin{description}
\item[(i)]
$S_0 \cap \{ z\in\Co : |z| \ge 1- \epsy_0 \} = \{1\}$.

\item[(ii)]
$\lll [Iz - (P- \one\otimes\pi)]^{-1}\lll_{V} 
\le \barb$ when $|z| \geq 1-\epsy_0$.

\item[(iii)]
$\haP_\alpha$ is $V$-uniform for all $\alpha=a+i\omega\in\Co$ 
satisfying
\[
|\omega| \le \epsy_0
\;\;\;\mbox{and}\;\;\;
|a|\le \bara \eqdef \Bigl( \frac{e-1}{2b-\delta}\Bigr)\delta\,.
\]
Moreover, the associated g.p.e.\ $\lambda_\alpha$ 
is an analytic function of $\alpha$ in this range, and
so is the corresponding eigenfunction $\cf_\alpha(x)$
(for each fixed $x\in\state$).

\item[(iv)]
The eigenfunctions $\cf_a$
are (uniformly) bounded from below when $a$ is real: 
\[
\inf_{-\bara\le a\le \bara} \cf_a(x) > 0,\qquad x\in\state\, .
\]
\end{description}
\end{proposition}

\proof
Results (i) and (ii) follow immediately from \Proposition{Ubdd}.
To prove (iii) we must establish an appropriate 
range of real $a$ for which $\haP_a$ is $V$-uniform.
From \Proposition{Ubdd} we know that
$P=\haP_0$ is $V$-uniform.

For any function $G_0\in L_\infty$, set $g_0=\exp(G_0)$, and
consider the  kernel $I_{g_0} R$, where,
as before, $I_{g_0}$ denotes the kernel
$I_{g_0}(x,\cdot)\eqdef g_0(x)\delta_x(\cdot)$.
We assume that the convergence
parameter for this kernel is equal to one.  It then follows from 
\Proposition{Ubdd} that the function below lies in
$L_\infty^V$ provided $\|g_0\|_\infty^{-1}> 1-B_1^{-1}$,
\[
\cg_r(x) = 
  \sum_{k=0}^\infty  [I_{g_0} (R-s\otimes \nu)]^k I_{g_0}s\, ,
\]
and it is clear that $I_{g_0} R$ is in fact $V$-uniform in this case.
Applying \Lemma{small}~(i) we know that $\nu(\cg_r) = 1$.

The function $\cg_r$ solves the eigenfunction equation,
\[
\begin{array}{rcl}
 [I_{g_0} (R-s\otimes \nu)] \cg_r
&=& \cg_r - I_{g_0} s
\\[.3cm]
\Longrightarrow
\quad
R \cg_r
&=& g_0^{-1} \cg_r .
\end{array}
\]
Setting $\cg=R \cg_r= \cg_r g_0^{-1}$ and applying the resolvent equation
\eq resolve/ then gives
\[
\begin{array}{rcl}
(P-I)\cg = (P-I)R \cg_r = (e-1) (R-I)\cg_r
&=&
(e-1)\cg - (e-1) g_0 \cg.
\end{array}
\]
Hence $\cg$ is the solution to the multiplicative Poisson equation for
the function $G =\log(g)\eqdef - \log ( e - (e-1)g_0)$.
The map $g_0\mapsto g$ is one to one.

We have already remarked that
$\cg_r\in L_\infty^V$ provided $\|g_0\|_\infty^{-1} > 1-B_1^{-1}$, 
and hence and 
$\cg\in L_\infty^V$ whenever $g_0$ satisfies this bound.
If $B_1\leq e$, then this constraint is trivially
satisfied. For $B_1>e$, equivalently
the function $g$ must satisfy,
\begin{equation}
\|g\|_\infty 
<
 \frac{1}{e- (e-1)(1-B_1^{-1})^{-1} }  =  \frac{B_1-1}{B_1-e}.
\elabel{gBdd}
\end{equation}
From the inequality $\log(1+x) < x$, $x\neq 0$, we obtain
\[
\log\Bigl( \frac{B_1-1}{B_1-e}\Bigr)
=
-
\log\Bigl( 1- \frac{e-1}{B_1-1}\Bigr)
> 
\Bigl( \frac{e-1}{B_1-1}\Bigr)
\ge
\Bigl( \frac{e-1}{2b-\delta}\Bigr)\delta,
\]
where the last inequality uses  the bound $B_1\le 2b /\delta $.

This gives the sufficient condition, $\| G\|_\infty \le \bara$.
\Proposition{Hcriterion} implies that $I_gP$ is $V$-uniform, and
$\cg\in L_\infty^V$ when this uniform bound holds.

The function $G$ falls outside of the class of functions 
$F$ satisfying \eq assumption/, since $\LA(a) > 0$ 
for all $a\neq 0$ when $\pi(F)=0$, and we have already 
noted that the spectral radius $\xi(g)$ of $I_g P$ is equal to $1$.
However, given any $a$, the function $G = aF -\LA(a)$ satisfies 
$\xi(g)=1$ and $G(x) \le |a|-\LA(a)< |a|$, $x\in\state$, so that the
normalized function satisfies \eq gBdd/ when $|a|\le \bara$.
This transformation immediately gives the desired conclusion in (iii).

To see (iv),  take any $\barlambda$ satisfying 
$\barlambda \ge \max(\lambda_{\bara},\lambda_{-\bara})$, and set
\[
\haG_a = \sum_0^\infty \barlambda^{-n-1} \haP^n_a\, .
\]
By irreducibility we can find $s_0\colon\state \to (0,1)$ and a probability 
distribution $\nu_0$ on $\clB$ satisfying the uniform bound,
\[
\haG_a(x,A) \ge R_\theta(x,A)\ge s_0(x)\nu_0 (A),\qquad x\in\state,\ A\in\clB,\ a\in [-\bara,\bara]\, ,
\]
where $\theta = \bara + \log(\barlambda)$.
We may assume that $\nu_0$ is equivalent to the irreducibility measure $\psi$.

It follows that for all $a\in [-\bara,\bara]$ and all $x$,
\[
(\barlambda-\lambda_a)^{-1} \cf_a(x)
=
\haG_a \cf_a\, (x) \ge s_0(x)\nu_0(\cf_a)>0\, .
\]
By continuity of $\cf_a$ we obtain the desired uniform bound.
\qed

\medskip

We now develop the consequences of the 
lattice condition. Our main conclusion 
is contained in \Theorem{lattice}:  
The function $F$ is almost-lattice 
{\em if and only if} the
spectral radius $\haxi_{i\omega}$ attains 
its upper bound (i.e. $\haxi_{i\omega} =1$)
for some $\omega>0$.

Some of the spectral theory for complex 
$\alpha$ is most easily developed
in a Hilbert space setting. Define 
$L_2 \eqdef \{ f\colon\state\to\Co
\;\mbox{such that}\;
\|f\|_2^2=\pi(|f|^2)<\infty\}$, 
with the natural associated inner product, 
$\langle h,g\rangle = \pi(h^*g)$, $h,g\in L_2$.  
We note that $V\in L_2$ under our standing 
assumption \eq assumption/.
For any $n$, the induced operator norm of 
$\haP_\alpha^n\colon L_2\to L_2$
may be expressed,
\[
\lll \haP_\alpha^n \lll_2 
=
\sup \frac{ \| \haP_\alpha^n g \|_2 }{\|g\|_2}
=  \sup\Bigl\{
| \Expect_{\pi}[ h^*(\Phi(0)) \exp(\alpha S_n) g(\Phi(n))] |
		: \|h\|_2\le 1, \|g\|_2 \le 1\Bigr\} \, .
\]
We let $\hagamma_\alpha$ denote 
the $L_2$-spectral radius, 
\[
\hagamma_\alpha\eqdef \lim_{n\to\infty}\lll \haP_\alpha^n \lll_2^{1/n}\,.
\]
When $\alpha = i\omega$,
the linear operators $\{\haP^n_{i\omega} \}$ 
are contractions on $L_2$, so that $\hagamma_{i\omega}\le 1$.

\Theorem{lattice} provides several 
characterizations of the almost-lattice condition. 
It is analogous to the variance 
characterization given in \Proposition{clt}.

\begin{theorem} 
\tlabel{lattice}
{\em (Characterization of Lattice Condition)}
The following are equivalent under \eq assumption/, for any given
$\omega>0$, 
$-\bara\le a\le \bara$:
\begin{description}
\item[(i)] 
$\haxi_{i\omega} =1$;

\item[(ii)] 
$\hagamma_{i\omega} =1$;

\item[(iii)] 
$\haxi_{a + i\omega} =\haxi_a$;

\item[(iv)]
There exists a bounded function $\Theta\colon\state\to [0,2\pi)$ and   
$d_0>0$ such that for a.e.\ $x\in\state$ $[\psi]$,
\begin{equation}
\exp\Bigl(i\omega \int_{[0,t)} (F(\Phi(s)) -  d_0)\, ds\Bigr)  
=  \exp\Bigl(i\Theta(\Phi(t)) - i \Theta(\Phi(0)) \Bigr)
 \, ,
\qquad a.s.\ [\Prob_x].
\elabel{multFishEq4}
\end{equation}
\item[(v)] 
$F$ is an almost-lattice function whose span is an integer
multiple of $2\pi/\omega$. 
\end{description}
\end{theorem}

\proof  
We first note that by \Proposition{clt} the existence of
$\Theta,\omega,d_0$ satisfying (iv) is equivalent to the 
almost-lattice condition (v).  To prove the proposition 
it remains to show that (i)--(iv) are equivalent.

The implications (iv) $\Rightarrow$ (i), (ii), (iii) 
are obvious since, under (iv), we have for all $n\ge 1$,
\[
\haP^n_{a+i\omega}(x,\varble) 
= e^{in d_0} I_{i\Theta} \haP_a^n(x,\varble)I_{i\Theta}^{-1}, \qquad
\hbox{for  a.e.\ $x\in\state\;[\psi]$.}
\]

We now establish implication (i) $\Rightarrow$ (ii).
We first note that if $\haxi_{i\omega}=1$ then, from
the fact that the $V$-norm is submultiplicative 
(as it is an operator norm), we must have 
$\lll \haP^n_{i\omega}\lll_V\ge 1$ 
for all $n$.  Note also that for any $g\in \LV$,
\[
\begin{array}{rcl}
|\haP_{i\omega}^{n+m} g\, (x) | &=&
|\Expect_x[\exp(\alpha S_n) \Expect_{\Phi(n)}[ \exp(\alpha S_m) g(\Phi(m))] ] |
\\[.25cm] 
&\le&
\Expect_x[| \Expect_{\Phi(n)}[ \exp(\alpha S_m) g(\Phi(m))]|] 
\\[.25cm]
&\le&
\int| \Expect_{y}[ \exp(\alpha S_m) g(\Phi(m))]| \, \pi(dy) 
                         + O(V(x) e^{-b_0 n}),\qquad n,m \ge 1,\ x\in\state,
\end{array}
\]
where $b_0>0$ exists by $V$-uniformity of $P$.
This implies the bound,
\begin{equation}
1=\haxi \le \liminf_{m\to\infty} \Bigl( \sup\Bigl\{
| \Expect_{\pi}[ h^*(\Phi(0)) \exp(\alpha S_m) g(\Phi(m))] |
		: \|h\|_\infty \le 1, \|g\|_V \le 1\Bigr\}\Bigr).
\elabel{omegaBdd}
\end{equation}

We have already remarked that $\haP_\alpha$ is a contraction on $L_2$.
It follows that either $\lll\haP^n_\alpha\lll_2 \to 0$ geometrically fast, 
or   $\lll\haP^n_\alpha\lll_2=1$ for all $n$. 
We may conclude the latter using \eq omegaBdd/, and this establishes 
the implication (i) $\Rightarrow$ (ii).

We now show that (ii) implies (iv).
The supremum in the definition of $\lll\haP^n_{i\omega}\lll_2$ is attained 
since $\hagamma_{i\omega}=1$.
To see this, construct for any $N\ge 1$ functions $h^N,g^N$ with $L_2$-norm 
equal to one, with 
\[
1
\ge
\| h^N\|_2 \| \haP^n_{i\omega} g^N\|_2
\ge
\langle h^N , \haP^n_{i\omega} g^N\rangle 
     \ge \lll\haP^n_\alpha\lll_2 -1/N
     = 1-1/N.
\]
Part of the construction ensures that the inner product above is real-valued.
These bounds imply that $\|h^N - \haP^n_{i\omega} g^N\|_2\to 0$, $N\to\infty$,
which is equivalently expressed,
\[
 \Expect_{\pi}[ | (h^N(\Phi(0)) )^* \exp(i\omega S_n) g^N(\Phi(n))]  -1 |]
              \to 0,\qquad N\to\infty.
\]
It then follows that $\exp(i\omega S_n)\in \sigma(\Phi(0),\Phi(n))$, and that
there exist $h_n,g_n\in L_2$ such that 
\begin{equation}
 h_n^*(\Phi(0)) \exp(i\omega S_n) g_n(\Phi(n)) = 1 \quad a.s.\ [\Prob_\pi]\, .
\elabel{xi1-id}
\end{equation}
We may assume without loss of generality that $|g_n(x)|=|h_n(x)|=1$ for all 
$x$ since $| \exp(i\omega S_n)|=1$.

Note that \eq xi1-id/ is almost the desired conclusion (iv).  
In particular, on 
dividing the expressions for $n$ and $(n+1)$ we obtain 
the suggestive identity,
\begin{equation}
\exp(i\omega F(\Phi(n))) =\left( \frac{h_{n+1}(\Phi(0))}{ h_n(\Phi(0)) }\right)
                       \left( \frac{g_n(\Phi(n))}{ g_{n+1}(\Phi(n+1)) }\right)
\qquad n\ge 0.
\elabel{suggF}
\end{equation}
To establish (iv) we show that $\{g_n, h_n\}$ may be \textit{chosen}
as follows: $\{g_n\}$ is  independent of $n$, with common value 
$g\in L_\infty$, and we may construct
$\theta_0\in\Re$ such that 
$h_n = e^{i\theta_0 n}g$ for all $n$.  
The required function $\Theta$ in (iii) can then be taken  as a version of 
$-\log(g)$.

Applying \eq xi1-id/ and appealing to stationarity, we conclude that
for any $n,m\ge 1$,
\[
\begin{array}{rrcl}
& h_{n+m}^*(\Phi(0)) \exp(i\omega S_{n+m}) g_{n+m}(\Phi(n+m)) &=&1  
\qquad a.s.\ [\Prob_\pi],
\\[.25cm]
\hbox{and} &&&
\\[.25cm]
& h_m^*(\Phi(n)) \exp\left\{(i\omega \sum_{k=n}^{n+m-1} F(\Phi(k))) 
	g_m(\Phi(n+m)) \right\}
&=&
\vartheta^n[ h_m^*(\Phi(0)\exp(i\omega S_m) g_m(\Phi(m)) ] 
\\[.2cm]
&&=&1
\qquad a.s.\ [\Prob_\pi],
\end{array}
\]
where $\vartheta^n$ denotes the $n$-fold shift operator on 
the sample space.

Combining \eq xi1-id/ with these two identities then gives,
\[
 h_{n+m}^*(\Phi(0))  h_n(\Phi(0)) 
g_n^*(\Phi(n)) h_m(\Phi(n))   g_m^*(\Phi(n+m))  g_{n+m}(\Phi(n+m)) 
=1.
\]
On taking conditional expectations with respect to $\Phi(0)=x$ we see that
for a.e.\ $x\in\state$ $[\psi]$,
\[
\begin{array}{rcl}
 h_{n+m}^*(x)  h_n(x) 
&=&
\Expect_\pi\Bigl[ 
g_n(\Phi(n)) h_m^*(\Phi(n))   g_m(\Phi(n+m))  g_{n+m}^*(\Phi(n+m)) \Bigr]
+ O(V(x) e^{-b_0 n} )
\\[.25cm]
&=&
\pi(g_nh_m^*)  \pi( g_m g_{n+m}^*) + O(V(x)e^{-b_0 n}+e^{-b_0 m} )
,\qquad n,m\ge 1.
\end{array}
\]
Since $|g_n(x)|=|h_n(x)|=1$ for all $x$ we  
conclude from Jensen's inequality that  for all $n,k\ge 1$,
\[
\begin{array}{rcl}
g_n^*(x)h_{n+k}(x)&=&\pi(g_n^*h_{n+k})  + \epsilon_1(x)
\\[.2cm]
 h_{n+k}^*(x)  h_n(x)  &=&
 \pi( h_{n+k}^*h_n) + \epsilon_2(x)
\end{array}
\]
where $|\epsilon_1(x)|+|\epsilon_2(x)| = O(V(x) e^{-b_0 n})$. 
This, combined with \eq suggF/,
shows that the desired expression can be obtained as an approximation:
For any $\epsilon >0$ we can find a function $\Theta$ 
(of the form $-\log(g_n)$ for large $n$) and $\theta_0\in\Re$ 
such that for a.e.\ $\Phi(0)=x\in\state$,
\[
\Bigl|
\exp(i\omega F(\Phi(0))  
-  \exp\Bigl(i(\theta_0+\Theta(\Phi(1)) - \Theta(\Phi(0)) )\Bigr) 
\Bigr|
\le \epsilon V(x)
 \quad a.s.\ [\Prob_x].
\]
This easily gives (iv).

Finally we show that (iii) implies (iv).
Observe first that we have already
established the equivalence of (i) and (iv).  
Moreover, (iii) is equivalent to the statement (i) 
for the transition kernel $\cP_a$, from which we 
deduce the implication (iii) $\Rightarrow$ property (iv) 
for the Markov chain with transition law $\cP_a$.  
This is equivalent to (iv) for
the original Markov chain.
\qed

\subsection{Continuous Time}
%%%%%%%%%%%%%%%%%%%%%%%%%% CONT'S TIME %%%%%%%%%%%%%%%%%%%%             
\slabel{SpectralCts}

We now translate the definitions and results of the previous 
section to the continuous-time case.
Suppose that $\{ \haP^t \,:\,t\in\Re_+\}$ is a semigroup 
of operators on $L_\infty^V$, with generator $\haclA$, 
and with finite spectral radius given by
\[
\haxi := \lim_{t\to\infty} 
\lll \haP^t \lll_{V}^{1/t}\, .
\]
[Note that the definition of the
generator of a positive semigroup
given in (\ref{eq:generator}), 
immediately generalizes to general
(not necessarily positive) 
semigroups.]

Consider the eigenvector equation $\haclA h = \LA h$.
The functions $h$ we consider will always be of the form 
$h=\haR_\theta h_0$,
usually with $h_0\ge 0$, where $\haR_\theta$ is defined 
as in
\eq  resolve/.  When all the integrals are well-defined 
we have the resolvent
equation \eq ResEqn/, so that 
\[
\haclA h = \theta  (\haR_\theta - I ) h_0.
\]
This identity allows us to lift all of the previous results to the
continuous-time setting.  In particular, under (V4), Fatou's lemma implies
that the resolvent $R=R_1$ satisfies,
\[
R V \le (1 -\delta_1) V + b_1 R s,
\]
with $\delta_1= \delta(1+\delta)^{-1}$, and  $b_1= b(1+\delta)^{-1}$.
It then follows as in the discrete time case that
$$\lll [I-R + s\otimes\nu]^{-1} \lll_V\le 2b_1/\delta_1.$$

%%IK5: 
Recall the definition of the semigroup $\{ \haP^t_\alpha \,:\,  t\in\RL_+\}$ 
from \eq haP/, where we now allow $\alpha$ to be possibly
complex. The next lemma offers an expression for the
generator of this semigroup, analogous to the classical
Feynman-Kac formula for diffusions. The result is easy to
check via the martingale representation \eq martGen/.

\begin{lemma}
{\em (Feynman-Kac Formula)}
The generator  $\haclA_\alpha$ of the semigroup  
$\{ \haP^t_\alpha \,:\,  t\in\RL_+\}$ satisfies,
\begin{equation}
\haclA_\alpha = \clA + \alpha F,
\elabel{FK}
\end{equation}
where $\clA$ is the generator of $\{P^t\}$.
\end{lemma}

Although none of the generators we consider 
are linear operators on $L_\infty^V$,
we may still define the \textit{spectrum} 
of $\haclA$, $\clS(\haclA)\subset \Co$,
as  the set of $z\in\Co$ such that 
the inverse $[Iz-\haclA]^{-1}$ does not exist as a bounded linear operator.
We have the generalized resolvent equation,
\[
 z [Iz-\haclA]^{-1}  = \haR_z =
\int_{[0,\infty)} z e^{-z t} \haP^t \,dt  \, , \qquad z\in \Co\, ,
\]
where the integral converges in norm for $z\not\in \clS(\haclA)$ 
such that $|e^z| \ge\haxi$.
The generator $\haclA$  is called \textit{$V$-uniform} 
if it admits a spectral gap and there is a unique 
pole $\LA_\circ\in\clS(\haclA)$ 
of multiplicity one, satisfying 
$|\lambda_\circ| =\exp(\LA_\circ) = \haxi$.

\begin{proposition} 
\tlabel{AgpeR}
If $\{ \haP^t \,:\,t\in\Re_+\}$ has finite spectral radius $\haxi$, then:
\begin{description}
\item[(i)]
The following statements are equivalent:
\begin{description}
\item[(a)]
The generator $\haclA$ has eigenvalue $\LA_\circ\in \Co$ 
and associated eigenfunction $\cf\in L_\infty^V$.

\item[(b)]
For $\theta > \haxi$, the resolvent $\haR_\theta$ has eigenvalue
$\lambda_\theta = \theta^{-1} \LA_\circ- 1$,
and eigenfunction $\cf$.
\end{description}
\item[(ii)]
The following statements are equivalent:
\begin{description}
\item[(a)]
$\haclA$ is $V$-uniform.
\item[(b)]$\haR_\theta$ is $V$-uniform.
\end{description}
\end{description} 
\end{proposition} 

The proof of \Proposition{AgpeR} is obvious from \eq gpekey/.
 
Using these identities, the following results may 
be proven as in \Theorem{mainSpectral} 
and \Theorem{mainSpectral2}.  The definition of 
$\bara$ is given in \eq bara/.

\begin{theorem}
\tlabel{mainSpectral-cts}
{\em (Multiplicative Mean Ergodic Theorem)}
Suppose that the Markov process 
$\bfPhi=\{\Phi(t) \,:\, t\in\Re_+\}$ 
and the functional $F$ satisfy \eq assumption/,
and write $\bara =\Bigl( \frac{e-1}{2b - \delta}\Bigr)$
as before.
Then there exists $\baromega>0,\delta_0>0$, 
such that for any $\alpha=a+i\omega\in\Co$ 
with
$| a | \leq \bara $,  $| \omega | \leq \baromega$,
there exists $\LA(\alpha) \in \clS_\alpha$ 
which is maximal and isolated:
\[ 
\Ree(\LA(\alpha))=
\max \{\Ree(\LA) :\LA\in\clS_\alpha\}
\quad
\hbox{and}
\quad
\clS_\alpha  
\cap \bigl\{ z : \Ree(z) \le \Ree(\LA(\alpha)) - \delta_0 \bigr\}  
=\LA(\alpha) .
\]

Moreover, for any such $\alpha$,  
there exist   $\cf_\alpha \in
L_\infty^V$ and $\cmu_\alpha \in \clM_1^V$, satisfying 
\eq cmu-norm/, and 
\begin{description}
\item[(i)]  
For all $x\in\state$,  $A\in\clB$, $t\in\Re_+$,
\begin{eqnarray*}
\haP_\alpha^t \cf_\alpha\, (x)   &=& \lambda_\alpha^t \cf_\alpha (x);
 \\
\cmu_\alpha \haP_\alpha^t \, (A)
 &=&  \lambda_\alpha^t \cmu  _\alpha (A).
 \end{eqnarray*}
The function $\cf_\alpha$ is also an eigenfunction for $\haclA_\alpha$:
\[
\haclA_\alpha\cf_\alpha = \LA(\alpha) \cf_\alpha.
\]

\item[(ii)] 
There exist $b_0>0$, $B_0 < \infty$, such that for all 
$x\in\state$,  $t>0$,
\begin{eqnarray*}
\Bigl| \Expect_x \bigl[ \exp (\alpha S_t
                      - t \LA(\alpha)) g (\Phi(t))\bigl] 
                                              - \haQ_\alpha g \, (x)\Bigr| 
& \leq & B_0 \| g \|_{V} e^{-b_0 t} V(x)
 \\
\Bigl| \Expect_x \bigl[ \exp (\alpha S_t - t \LA(\alpha))\bigr]
                                           - \cf_\alpha(x)\Bigr|
                     & \leq &
        B_0 | \alpha | e^{-b_0 t}  V(x).
\end{eqnarray*}
\end{description}
\end{theorem}

\begin{theorem}
\tlabel{mainSpectral2-cts}
{\em (Bounds Around the $i\omega$-Axis)}
Assume that the Markov process 
$\bfPhi=\{\Phi(t) \,:\, t\in\Re_+\}$ 
and the functional $F$
satisfy \eq assumption/.
\begin{itemize}
\item[(NL)]
Suppose that $F$ is strongly non-lattice.
For any $0<\omega_0 <\omega_1< \infty$,
there exist $b_0>0$, $B_0 < \infty$
(possibly different than above), such that
\be
\Bigl | \Expect_x [ \exp (\alpha S_t - t\LA(a))] \Bigr|
\leq
B_0 V(x) e^{-b_0 t}\, ,  \qquad  x \in \state,\; t>0,
\label{eq:ImMMET2}
\ee
for all $\alpha=a+i\omega$ with $|a| \le \bara$
and $\omega_0\leq |\omega|\leq\omega_1$.
 
\item[(L)]
Suppose that $F$ is almost-lattice with span $h>0$.
For any $\epsilon>0$, there exist
$b_0>0$, $B_0 < \infty$
(possibly different than above),
such that (\ref{eq:ImMMET2}) holds
for all $\alpha=a+i\omega$ with $|a| \le \bara$
and $\epsilon\leq|\omega|\leq 2\pi/h - \epsilon.$
\end{itemize}
\end{theorem}

\section{Edgeworth Expansions for the CLT}
\slabel{edgeworth}
%%%%%%%%%%%%%%%%%%%%%%%%%%% Edgeworth Expansions %%%%%%%%%%%%%%%%%%%%%%%%%%%

Here we show how the multiplicative mean
ergodic theorems of the previous section 
can be used to obtain Edgeworth expansions
for the central limit theorem (CLT) 
satisfied by the partial sums of
a geometrically ergodic Markov chain;
see, e.g., 
\cite[Ch.~17]{meyn-tweedie:book} 
for the standard CLT.

Throughout this section we 
consider a discrete-time Markov
chain $\bfPhi$ and a bounded functional 
$F:\state\to\RL$. Recall our standing 
assumptions \eq assumption/ about
$\bfPhi$ and $F$. To avoid repetitions
later on, we collect below a number of
properties that will be used repeatedly 
in the proofs of the results in this and
the following section. They are proved
in the Appendix.

\paragraph{Properties.} Assume that the
discrete-time Markov chain $\bfPhi$ and 
the function $F:\state\to\RL$ satisfy 
\eq assumption/, and let
$S_n$ denote the partial sums as before.
Choose and fix an arbitrary $x\in\state$, and let       
\be
m_n(\alpha):=\Expect_x[\exp(\alpha S_n)]\,,\qquad n\ge 1,\
        \alpha\in\Co.
\label{eq:mgf}
\ee
\begin{itemize}
\item[P1.] 
        There is a sequence $\{\epsilon_n\}$ 
        such that
        $$m_n(\alpha)=\exp(n\LA(\alpha))
        [\cf_\alpha (x)+ |\alpha| \epsilon_n]\, ,
        \quad n\geq 1\, ,$$
        and $|\epsilon_n|\to 0$ 
        exponentially fast
        as $n\to\infty$,
        uniformly over all $\alpha\in\Omega$
        (with $\Omega$ as in \Theorem{mainSpectral}).
\item[P2.]
        If $F$ is strongly non-lattice, then for any 
        $0<\omega_0 <\omega_1< \infty$
        there is a sequence $\{\epsilon_n'\}$ 
        such that
        $$m_n(\alpha)=\exp(n\LA(a))\epsilon'_n\, ,
        \quad n\geq 1\,,$$
        and $|\epsilon'_n|\to 0$ 
        exponentially fast
        as $n\to\infty$,
        uniformly over all 
        $\alpha=a+i\omega$ with $|a| \le \bara$
        and $\omega_0\leq|\omega|\leq\omega_1$
        (with $\bara$ as in \Theorem{mainSpectral}).
\item[P3.]
        If $F$ is lattice (or almost lattice) with span $h>0$,
        then for any $\epsilon>0$, as $n\to\infty$,
        $$\sup_{
        \epsilon\leq|\omega|\leq 2\pi/h - \epsilon
        }|m_n(i\omega)| \to 0
        \qquad\mbox{exponentially fast.}$$
\item[P4.] 
        The function
        $\LA(\cdot)$ is analytic in $\Omega$ with $\LA(0)=\LA'(0)=0$,
                and $\LA''(0)=\sigma^2>0.$
                Moreover, $\sigma^2_a:=\LA''(a)$
                is strictly positive for all 
                real $a\in[-\bara,\bara]$.
\item[P5.] 
        The third derivative $\rho_3:=\LA'''(0)$ 
        can be expressed as
        \ben
        \rho_3
        &=&
 \Expect_\pi[ F^3(\Phi(0))]
+3 \sum_{{i=-\infty \atop i\neq 0}}^\infty \Expect_\pi[F^2(\Phi(0))F(\Phi(i))]
\\
&&\qquad 
+6\sum_{i,j=1}^\infty \Expect_\pi[ F(\Phi(0)) F(\Phi(i)) F(\Phi(i+j)) ].
        \een
\item[P6.] 
        Let $\haF$ be the solution of the Poisson equation
        given by (\ref{eq:Fhat}), and write
        $$\Delta_n:=\Expect_x[S_n] - \haF(x)\,.$$
        Then $|\Delta_n|\to 0$
        exponentially fast
        as $n\to\infty$.
\item[P7.] 
                The eigenfunction $\cf_\alpha$ 
                is analytic
                in $\alpha\in\Omega$, 
                it satisfies 
                $\cf_\alpha\big|_{\alpha=0} \equiv 1$,
                and it is strictly positive
		for real $\alpha$. Moreover,
		there is some $\baromega_0\in(0,\baromega]$
                (depending on $x$),
                such that 
                $$\delta(i\omega):=|\log \cf_{i\omega}(x)
                        -i\omega\haF(x)|\leq (\mbox{Const})\omega^2,$$
                for all $|\omega|\leq\baromega_0$,
                where $\haF$ is as in P6.
\end{itemize}

\noindent
The following two results generalize
those in 
\cite{nagaev:61,jensen:91,mccormick-data:93}.

\begin{theorem}
\tlabel{EdgeworthNL}
{\em (Edgeworth Expansion for Non-Lattice Functionals)}
Suppose that $\bfPhi$ and the strongly-non-lattice 
functional $F$ satisfy assumption \eq assumption/, and
let $G_n(y)$ denote the 
distribution function of the 
normalized partial sums 
$S_n/\sigma\sqrt{n}$:
$$G_n(y)\eqdef\Probsub_x\left\{
        \frac{S_n}{\sigma\sqrt{n}}\leq y
        \right\},
        \quad
        y\in\RL.$$
Then, for all $x\in\state$,
\be
G_n(y)
={\cal G}(y) 
+ 
\frac{\gamma(y)}{\sigma\sqrt{n}} 
\left[
        \frac{\rho_3}{6\sigma^2}(1-y^2)
        \,-\,\haF(x)
\right]
+o(n^{-1/2}),
\quad
n\to\infty,
\label{eq:EdgeNL}
\ee
uniformly in $y\in\RL,$
where $\gamma(y)$ denotes the standard Normal density
and ${\cal G}(y)$ is the corresponding distribution
function.
\end{theorem}

It is perhaps worth noting 
the way in which the convergence in (\ref{eq:EdgeNL})
depends on the initial state $x$ of the 
Markov chain: This dependence is only
manifested via the solution
$\haF(x)$ to the Poisson equation.
Also observe that,
since (\ref{eq:EdgeNL}) holds 
for all $y\in\RL$, 
the restriction on $F$ being $|F|\leq 1$ 
can clearly be relaxed to $\|F\|_\infty<\infty$.

For the proof of the theorem --
given in the Appendix -- it is convenient to consider 
the zero-mean version of the normalized partial 
sums, 
$$\frac{\oo{S}_n}{\sigma\sqrt{n}}:=
\frac{S_n-\Expect_x\{S_n\}}{\sigma\sqrt{n}}.$$
Let $\barG_n(y)$ denote the corresponding
distribution function. In the proof we show instead that
\be
\barG_n(y)
={\cal G}(y)
+
        \frac{\rho_3}{6\sigma^3\sqrt{n}}
        (1-y^2)
        \gamma(y)
+o(n^{-1/2}),
\quad
n\to\infty,
\label{eq:EdgeNL2}
\ee
uniformly in $y\in\RL$.
From this it is a straightforward
calculation to deduce 
(\ref{eq:EdgeNL}) via
a Taylor series expansion
and using property P6.

Before stating our next result we 
recall the following notation.
If $G$ is the distribution 
function of a lattice random 
variable with values on the lattice
$\{d+kh, \;\;k\in\IN\}$, {\em the 
polygonal approximation $G^\#$ to $G$} 
is the piecewise-linear distribution
function $G^\#(y)$ that agrees with
$G(y)$ at the mid-points of the
lattice, $y=d+(k+1/2)h$,
$k\in\IN$, and is linearly interpolated
between these points. The function
$G^\#$ is precisely the convolution
of $G$ with the uniform distribution
on $[-h/2, h/2]$.

\begin{theorem}
\tlabel{EdgeworthL}
{\em (Edgeworth Expansion for Lattice Functionals)}
Suppose that $F$ is a lattice functional
with span $h>0$, and assume that $F$ and $\bfPhi$
satisfy assumption \eq assumption/.
With $G_n(y)$ as in 
\Theorem{EdgeworthNL}, let 
$G_n^\#(y)$ denote its polygonal 
approximation.
Then, for all $x\in\state$,
\be
G^\#_n(y)
={\cal G}(y)
+
\frac{\gamma(y)}{\sigma\sqrt{n}}
\left[
        \frac{\rho_3}{6\sigma^2}(1-y^2)
        \,-\,\haF(x)
\right]
+o(n^{-1/2}),
\quad
n\to\infty,
\label{eq:EdgeL}
\ee
uniformly in $y\in\RL.$
In particular, 
writing $h_n=h/\sigma\sqrt{n}$,
(\ref{eq:EdgeL}) holds
with $G_n(y)$ in place of $G^\#_n(y)$ at
the points $\{y=(k+1/2)h_n,\;k\in\IN\},$
and with 
$[G_n(y)+G_n(y-)]/2$
in place of $G^\#_n(y)$ at the
points $\{y=kh_n,\;k\in\IN\}$.
\end{theorem}

The proof is given in the Appendix.
As with 
\Theorem{EdgeworthNL},
it is more 
convenient to prove
a version of (\ref{eq:EdgeL}) 
in terms of $\barG^\#_n(y)$ 
rather than $G^\#_n(y)$, where
$\barG^\#_n$ is the polygonal
approximation to $\barG_n$.
In the proof we show that
\be
\barG^\#_n(y)
={\cal G}(y)
+
        \frac{\rho_3}{6\sigma^3\sqrt{n}}
        (1-y^2)
        \gamma(y)
+o(n^{-1/2}),
\quad
n\to\infty,
\label{eq:EdgeL2}
\ee
uniformly in $y\in\RL.$
Then (\ref{eq:EdgeL}) 
follows from (\ref{eq:EdgeL2})
in the same way that 
(\ref{eq:EdgeNL}) follows from 
(\ref{eq:EdgeNL2}).

Before moving on to large deviations we note
that, although we shall not pursue these
directions further in this paper, using the
multiplicative mean ergodic theorems of
\Section{spectral} it is possible to 
prove higher-order Edgeworth expansions,
as well as precise local limit theorems
for the density (or the pseudo-density,
when a density does not exist) of $S_n$.
The Edgeworth-expansion proofs follow 
the same outline as those in the case 
of independent random variables; 
cf. \cite[p.~541]{fellerII:book}.
For the local limit theorems, one can
apply directly the general results of
\cite[Sec.~2]{chaganty-sethuraman:93}.

\section{Moderate and Large Deviations}
\slabel{ldp}
%%%%%%%%%%%%%%%%%%%%%%%%%%% LDPs %%%%%%%%%%%%%%%%%%%%%%%%%%%

In this section
we use the multiplicative mean ergodic theorems of
\Theorem{mainSpectral}
and \Theorem{mainSpectral2}
to prove moderate
and large deviations
results for the
partial sums of a
Markov chain.
As in \Section{edgeworth},
we consider the partial sums
$\{S_n\}$ of a bounded
functional $F$ of the 
discrete-time,
geometrically ergodic Markov 
chain $\bfPhi$.

First we note that the multiplicative mean
ergodic theorem together with the
analyticity of $\LA(\alpha)$ in
a neighborhood of the origin
(see properties P1 and P4 in
the previous section) immediately 
imply that the partial sums $S_n$ 
satisfy a moderate deviations 
principle (MDP). We state this
MDP, without proof, in 
\Proposition{MDP}. Its
proof is based on an application 
of the G\"{a}rtner-Ellis theorem, 
exactly as in the proof of Theorem~3.7.1 
in \cite{dembo-zeitouni:book}.

\begin{proposition}
\tlabel{MDP}
{\em (Moderate Deviations)
\cite{aco97a,acoche98a} }
Suppose the Markov chain $\bfPhi$ and the functional
$F$ satisfy \eq assumption/, and
let $\{b_n\}$ be a sequence of constants such that
$$\frac{b_n}{\sqrt{n}}\to\infty\quad\mbox{and}
\quad
\frac{b_n}{n}\to 0,\quad\;\; n\to\infty.$$
Then,
for all $x\in\state$ and 
any measurable $B\subset\RL$,
\ben
-\inf_{y\in B^\circ}\left(\frac{y^2}{2\sigma^2}\right)
&\leq&
\liminf_{n\to\infty}\frac{1}{b_n^2/n}\log\Probsub_x
        \left\{\frac{S_n}{b_n}\in B\right\}\\
&\leq&
\limsup_{n\to\infty}\frac{1}{b_n^2/n}\log\Probsub_x
        \left\{\frac{S_n}{b_n}\in B\right\}
        \;\leq\;
        -\inf_{y\in \barB}\left(\frac{y^2}{2\sigma^2}\right),
\een
where $B^\circ$ denotes the interior of $B$ and
$\barB$ denotes its closure.
\end{proposition}

Note that 
the same result
holds for the centered random variables 
$[S_n-\Expect_x\{S_n\}]/b_n$
in place of $S_n/b_n$.

\subsection{Large Deviations for Doeblin Chains}

Suppose that $\bfPhi$ is a Doeblin recurrent
chain, that is, suppose that for some $m\geq 1$,
$\epsilon'>0$, and a probability measure $\nu'$,
we have that $P^m\geq \epsilon'\nu'$. Equivalently,
the Doeblin condition can be stated as
\be
R\geq \epsilon\nu\,,
\qquad\mbox{for some $\epsilon>0$ and a  probability measure $\nu$}\,,
\label{eq:doeblin}
\ee
and this, in turn can be seen to be
equivalent to geometric ergodicity with a
{\em bounded} Lyapunov function $V$ in (V4);
see
\cite[Theorem~16.0.2]{meyn-tweedie:book}.
Then the state space $\state$ is small,
and the results of \cite{ney-nummelin:87b} 
can be applied to get large deviations 
results for the partial sums $S_n$.
For example, for a Doeblin chain
with a countable state space $\state$
and with $\psi$=counting measure, 
the partial sums $S_n$ satisfy 
a large deviations principle
(LDP) under the distributions
$\Probsub_x$, for any $x\in\state$.

But the situation is more complicated
when $\bfPhi$ is stationary,
i.e., when $\Phi(0)\sim\pi.$
In the following proposition we consider
the LDP for the partial sums $S_n$ under
the stationary distribution $\Probsub_\pi$.

\begin{proposition}
\tlabel{stationaryLDP}
{\em (Large Deviations)} 
Suppose the
Doeblin chain $\bfPhi$ and the functional
$F$ satisfy \eq assumption/,
and let $\bara=(\frac{e-1}{2-\epsilon})\epsilon$,
where $\epsilon$ is as in (\ref{eq:doeblin}).
\begin{description}
\item[(i)]
The partial sums $S_n$ satisfy an LDP
in a neighborhood of the origin: For any 
$c\in (0,\LA'(\bara))$ and 
any $c'\in (\LA'(-\bara),0)$, we have
\ben
\lim_{n\to\infty} \smalloneOvern\!\log\Probsub_\pi\{S_n\geq nc\}
        &=&
        - \LA^*(c)
        \\
\lim_{n\to\infty} \smalloneOvern\!\log\Probsub_\pi\{S_n\leq nc'\}
        &=&
        - \LA^*(c')\,,
\een
where 
$$\LA^*(c)\eqdef\sup_{-\bara < a < \bara}[ac-\LA(a)].$$
\item[(ii)] Part (i) cannot in general be extended to
a full LDP on the whole real line.
\end{description}
\end{proposition}

\proof
Integrating the multiplicative mean
ergodic theorem in (\ref{eq:mainMMET})
with respect to $\pi$ and noting
that $\pi(\cf_a)\in(0,\infty)$
for all $|a|\leq\bara$, we get
that
$$\smalloneOvern\! \log \Expect_\pi[\exp(a S_n)]\to\LA(a),
\qquad n\to\infty,$$
for all real $a\in[-\bara,\bara].$ Since
$\LA(a)$ is analytic, 
(i)
follows from the G\"{a}rtner-Ellis theorem
\cite[Theorem~2.3.6]{dembo-zeitouni:book}.
To see that in the Doeblin case 
$\bara=(\frac{e-1}{2-\epsilon})\epsilon$, note that 
in (V4) we can set $V\equiv 1$, 
$s\equiv \epsilon$, take $0<\delta<1$ 
be arbitrary, and define $b=\delta/\epsilon$.
We then have a version of (V4),
\[
PV = V = (1-\delta) V + b s.
\]
Using the definition of $\bara$ given in 
\Theorem{mainSpectral} then gives,
\[
\bara\eqdef  
\Bigl( \frac{e-1}{2b - \delta}\Bigr) \delta   
=
\Bigl( \frac{e-1}{2(\delta/\epsilon) - \delta}\Bigr) \delta   
=
\Bigl( \frac{e-1}{2 - \epsilon}\Bigr) \epsilon.
\]

Part~(ii) follows from the counter-example
in Proposition~5 of \cite{bryc-dembo:96}.
\qed

\subsection{Exact Large Deviations for Geometrically Ergodic Chains}

Next we consider the more general case of geometrically
ergodic Markov chains, satisfying our standing
assumptions \eq assumption/.
With $\bara$ as in \Theorem{mainSpectral},
let $(A',A)$ denote the interval
$$(A',A):=\{\LA'(a)\,:\,-\bara<a<\bara\},$$
and note that $0=\pi(F)=\LA'(0)\in(A',A)$.
Recall the definition of $\LA^*(c)$ in 
\Proposition{stationaryLDP}.

\begin{theorem}
\tlabel{Bahadur-RaoNL}
{\em (Exact Large Deviations for Non-Lattice Functionals)}
Suppose that $\bfPhi$ and the strongly-non-lattice 
functional $F$ satisfy \eq assumption/, and
let $c\in(0,A).$ 
Then, for all $x\in\state$,
\ben
\Probsub_x\{S_n\geq nc\}
       \;\sim\; 
        \frac{\cf_a(x)}{a\sqrt{2\pi n\sigma_a^2}}
        e^{-n\LA^*(c)},
        \quad
        n\to\infty,
\een
where $a$ is chosen so that
$\LA'(a)=c$, and 
$\sigma_a\eqdef \LA''(a)$.
A corresponding 
result holds for the lower tail.
\end{theorem}

It is perhaps worth pointing out that 
the way in which the 
large deviations probabilities 
$\Probsub_x\{S_n\geq nc\}$ depend on the initial 
state $x$ of the Markov chain is via
the solution $\cf_a(x)$ to the multiplicative
Poisson equation.

Although the proof (given next) relies on an 
application of a general result from 
\cite{chaganty-sethuraman:93}, the main 
idea is similar to the proof of the
corresponding result for independent 
random variables \cite{bahadur-rao:60}:
First, as in the case of finite state
space \cite{miller:61}, we perform 
a change of measure that maps the 
transition kernel $P$ to the twisted 
kernel $\cP_a$.
Since $\bfPhi$ is geometrically ergodic,
by \Proposition{basic1} $\haP_a$ is $V$-uniform.
Therefore $\cP_a$ is $V_a$-uniform by \Proposition{V-geo2},
and hence it is geometrically ergodic 
by \Corollary{Vuniformity}. Therefore we 
can apply the Edgeworth expansions of 
\Section{edgeworth}, and complete the 
proof along the lines of the corresponding
argument in the case of independent 
random variables; see, e.g., 
\cite[Theorem~3.7.4]{dembo-zeitouni:book}.

\proof
Choose and fix an arbitrary $x\in\state$.
The result of the theorem will follow by an
application of 
\cite[Theorem~3.3]{chaganty-sethuraman:93}.
We consider the moment generating functions
$m_n(\alpha)$ of $S_n$, defined in (\ref{eq:mgf})
for $\alpha$ in the interior
of the compact set $\Omega$ in \Theorem{mainSpectral}.
[Note that, although our $\Omega$ is different from the 
open disc used in \cite{chaganty-sethuraman:93},
a close examination of the proof of 
\cite[Theorem~3.3]{chaganty-sethuraman:93}
shows that the result continues to hold
when the open disc of radius $\bara$ is replaced
with the interior $\{\alpha=a+i\omega:|a|<\bara,\,
|\omega|<\baromega\}$ of the strip $\Omega$,
as long as $\baromega>0$.]

We will make repeated use of the properties
P1 -- P7 stated in \Section{edgeworth}.
From the definition of $m_n(\alpha)$ it is
easily seen that it is an analytic function
of $\alpha$, and 
from P1 and P4 it follows that $m_n(\alpha)$
is nonzero on $\Omega$, for all $n$ large
enough (uniformly in $\alpha$).

Let $\Lambda_n(\alpha)$ be the 
normalized log-moment
generating function
$$\LA_n(\alpha)\eqdef\frac{1}{n}\log m_n(\alpha)\,,
\qquad \alpha\in\Omega\,,$$
and 
$$\LA^*_n(c)\eqdef\sup_{-\bara<a<\bara}
[ac-\LA_n(a)]\,,
\qquad c\in\RL.$$
The main step in the proof is the verification
of the assumptions of 
\cite[Theorem~3.3]{chaganty-sethuraman:93}.
Most of them, plus some other 
technical properties, are established in
the following lemma (proved in the
Appendix).

\newpage

\begin{lemma}
\tlabel{br-lemma} 
Under the assumptions of
the theorem:
\begin{description}
\item[(i)] For $n$ large enough
there is a unique $a_n\in(0,\bara)$
such that
$\LA_n'(a_n) = c$
and 
$\LA^*_n(c)=a_nc-\LA_n(a_n).$
\item[(ii)]
Similarly, there is a unique
$a\in(0,\bara)$ such that 
$\LA'(a) = c$
and $\LA^*(c)=ac-\LA(a).$
\item[(iii)]
$a_n\to a$ as $n\to\infty$, and, in fact,
$a_n-a = O(\smalloneOvern\!).$
\item[(iv)]
$\LA_n''(a_n)\to\sigma^2_a,$ as
$n\to\infty$.
\item[(v)]
$\LA^*_n(c)\to \LA^*(c)$ as $n\to\infty$, and, in fact,
$$\LA^*_n(c) = \LA^*(c)
-\smalloneOvern\!\log \cf_{a}(x)+
o(\smalloneOvern\!).$$
\end{description}
\end{lemma}

The theorem follows from
\cite[Theorem~3.3]{chaganty-sethuraman:93},
upon verifying condition~$(c)$ of
\cite[p.~1685]{chaganty-sethuraman:93}.
For that, it suffices to show 
that for all $0<\omega_0<\omega_1<\infty$,
$$\sup_{\omega_0\leq|\omega|\leq \omega_1}
\left|\frac{m_n(a'+i\omega)}{m_n(a')}\right|=o(n^{-1/2})\,,$$
uniformly in $a'$ in a neighborhood of $a$.
But the above convergence actually
takes place exponentially fast,
as can be easily verified
using properties~P1 and~P2 from 
\Section{edgeworth}. 
\qed

\begin{theorem}
\tlabel{Bahadur-RaoL}
{\em (Exact Large Deviations for Lattice Functionals)}
Suppose that $\bfPhi$ and the 
lattice functional $F$
satisfy \eq assumption/, and assume
that $F$ has span $h>0$.
Let $\{c_n\}$ be a sequence of real 
numbers in $(\epsilon,A-\epsilon)$,
for some $\epsilon>0$, and assume 
(without loss of generality) that,
for each $n$, $c_n$ is in the support of $S_n$.
Then, for all $x\in\state,$
\be
\Probsub_x\{S_n\geq nc_n\}
       \;\sim\;
        \frac{h}
        {(1-e^{-ha_n})\sqrt{2\pi n\LA_n''(a_n)}}
        e^{-n\Lambda^*_n(c_n)},
        \quad
        n\to\infty,
\label{eq:brL1}
\ee
where each $a_n\in(0,\bara)$ is chosen 
so that $\LA_n'(a_n)=c_n$.
A corresponding 
result holds for the lower tail.
\end{theorem}

Note that in the lattice case we have given
a slightly more general version of the result
given in \Theorem{Bahadur-RaoNL}. If it 
turns out to be the case that the $c_n$
converge to some $c\in(\epsilon, A-\epsilon)$,
so that the corresponding $a_n$ converge
to some $a\in(0,\bara)$ at a rate $O(1/n)$, 
then applying \Lemma{br-lemma} as before, 
from (\ref{eq:brL1}) we obtain, 
\ben
\Probsub_x\{S_n\geq nc_n\}
       \;\sim\;
        \frac{h\cf_a(x)}
        {(1-e^{-ha})\sqrt{2\pi n\sigma_a^2}}
        e^{-n\Lambda^*(c)},
        \quad
        n\to\infty,
\een
where $\sigma_a^2=\LA''(a)$.

\proof
Choose and fix an arbitrary
$x\in \state$.
The proof 
parallels that of \Theorem{Bahadur-RaoNL},
relying on an application of Theorem~3.5 from
\cite{chaganty-sethuraman:93}. 
A close examination of its proof in 
\cite{chaganty-sethuraman:93} shows that,
as in the case of Theorem~3.3 above,
Theorem~3.5 remains valid if we replace 
the open disc of radius $\bara$ by
the interior of the strip $\Omega$. Proceeding
as in the proof of \Theorem{Bahadur-RaoNL},
we now need to verify condition $(c')$ on
\cite[p.~1686]{chaganty-sethuraman:93}.
For that, it suffices to show that for 
that for all $\omega_0\in(0,\pi/h)$,
$$\sup_{\omega_0<|\omega|\leq \pi/h}\left|
\frac{m_n(a'+i\omega)}{m_n(a')}\right|=o(n^{-1/2})\,,$$
uniformly in $a'\in (\epsilon,A-\epsilon)$.
Using properties~P1 and~P3 from
\Section{edgeworth}, it is easy to 
see that the above convergence actually
takes place exponentially fast,
and this completes the proof.
\qed

\section{Examples}
%%%%%%%%%%%%%%%%%%%%%%%%%%% Examples %%%%%%%%%%%%%%%%%%%%%%%%%%%
\slabel{examples}

\subsection{Countable State Space Models}
%%%%%%%%%%%%%%%%%%%%%%%%%% Countable State %%%%%%%%%%%%%%%%%%%%%%%%%%% 

Let $\bfPhi$ be a discrete-time Markov chain 
with a countable set $\state$ of states, and
let $\psi$ be counting measure. Suppose
$\bfPhi$ 
is irreducible in the usual sense that 
$R(x,y) > 0$ for all $x,y \in \state$.
Then $\ind_\theta$ is a small function
for any $\theta\in\state$, with associated small 
measure $\nu=P(\theta,\cdot)$.  Using this small function 
and measure in \Lemma{small}~(i)
leads to the following characterization of $\LA(a)$ for real $a$,
\begin{equation}
\LA(a)
=\inf\Bigl\{ \LA : 
\Expect_\theta
\Bigl[ \exp\Bigl(\sum_{k=0}^{\tau_\theta-1} 
 [a F(\Phi_k) -\LA]  \Bigr)  \Bigr]  \le  1 \Bigr\}\,;
\elabel{convPar3}
\end{equation}
see \cite{balaji-meyn}
for details.
When the infimum is attained and we may 
justify differentiation with respect to $a$, 
then 
\[
1=
\Expect_\theta
\Bigl[ \exp\Bigl(\sum_{k=0}^{\tau_\theta-1} 
 [a F(\Phi_k) -\LA(a)]  \Bigr)  \Bigr]
\quad\Longrightarrow\quad
0=
\Expect_\theta
\Bigl[ \sum_{k=0}^{\tau_\theta-1} [ F(\Phi_k) -\LA'(0)]   \Bigr]\, .
\]
This gives a more transparent proof of the identity 
$\Lambda'(0)=\pi(F)$.

\subparagraph{The simple queue.}
For our purposes, the simplest interesting example 
of a countable state space chain 
is the M/M/1 queue. This is the 
reflected random walk $\bfPhi$ on 
$\state =\{0,1,2,\dots\}$, with 
\[
P(x, x+1) = p,\  P(x,(x-1)_+) = q,\qquad x\in \state, 
\]
where $p+q=1$. We assume that $\rho = p/q<1$ 
so that the chain is positive recurrent.
As we show next:
\begin{itemize}
\item[(a)]
$\bfPhi$ is geometrically ergodic;
\item[(b)]
it is not Doeblin recurrent;
\item[(c)]
with $F=\ind_{0^c}-\pi(0^c)$,
the multiplicative mean
ergodic theorem 
(\ref{eq:mainMMET})
does not hold 
for all real $\alpha$.
\end{itemize}
It is also not hard to show that $\bfPhi$
does not satisfy (mV3) for any 
$f$ with finite sublevel sets,
so that, in view of the 
discussion in 
\Section{ergodictheorems},
the Donsker-Varadhan conditions 
do not apply.
More importantly, as Wu
recently showed, not just
the conditions, but also the
large deviations conclusions
of the Donsker-Varadhan theory
fail in this case \cite{wu:00}.
Therefore, this example does not fall under
any of the standard conditions known to 
imply large deviations results. 

Below we also show that
our central technical result, the multiplicative
mean ergodic theorem (\ref{eq:mainMMET}), cannot in
general be extended to hold on the entire real line.

First note that one can 
compute directly the expectations,
\begin{equation}
\Expect_x [r^{\tau_0}] 
=
\left\{
\begin{array}{rcl}
< & \infty, \quad & 0\le r \le \barbeta \\
= & \infty, \quad & r > \barbeta
\end{array}
\right. \, ,
\elabel{MM1r}
\end{equation}
where $\barbeta = (4q p)^{-\half}>1$.  
To construct a Lyapunov function, consider
$V(x) = r_0^x$, for $r_0 > 1$: 
\begin{equation}
PV(x) =
\left\{
\begin{array}{rcl}
(p r_0 + q r_0^{-1}) V(x) && x \geq 1 \\
(p r_0 + q) V(0) && x = 0.
\end{array}
\right.
\elabel{MM1geo}
\end{equation}
Choosing a minimal value for $(1-\delta)\eqdef (p r_0 + q r_0^{-1})$ 
gives $r_0 = \rho^{- \half}$ and a solution to (V4):
\[
PV(x)  = \sqrt{4p q} V(x) =\barbeta^{-1} V(x), \qquad x \geq 1.
\]
It easily follows that, with $r = \barbeta$, 
$\Expect_x [ r^{\tau_0}] = V(x) = \rho^{-x/2}$,  $ x \geq 1$.
This gives the finite bound in \eq MM1r/, and shows that 
$\bfPhi$
is geometrically ergodic with Lyapunov function $V$.

The easiest way to see that 
$\bfPhi$ is not Doeblin recurrent
is to notice that, in $k$ time steps,
$\bfPhi$
cannot visit
more than its $2k$ neighboring states, 
which implies that the state space is 
not small; see \cite[Theorem~16.0.2]{meyn-tweedie:book}.

Now let $F=\ind_{0^c}-\pi(0^c).$ 
Using the characterization \eq convPar3/ 
with $\theta=0$,
we find that 
$\LA(a)$ is the unique
solution to the fixed-point equation,
\[
\Expect_0[\exp\{(\pi_0a-\Lambda(a)) \tau_0\}] = e^a\,,
\]
where $\pi_0:=\pi(0).$
It follows from 
\eq MM1r/ 
that $\exp\{\pi_0a-\Lambda(a)\} \le \barbeta$
for all $a\in\Re$. Also,
from the fixed-point equation
it follows that if
$a^*:=\log\Expect_0[\barbeta^{\tau_0}]$,
then
$\LA(a^*)=\pi_0a^*-\log\barbeta$.
But since $\exp\{\pi_0a^*-\LA(a^*)\}=\barbeta$,
and $(\pi_0a-\LA(a))$ is nondecreasing in $a$,
from 
\eq MM1r/ 
we conclude that $\LA(a)=a-\log\barbeta$
for all $a\geq a^*$, and hence
$\LA''(a)=0$ for $a\geq a^*$.
[To see that $(\pi_0a-\LA(a))$ is nonincreasing,
simply recall from \Proposition{V-geo2}
that $\LA'(a)=\pi(a)$
so that $\LA'(a)\leq\sup_x F(x) = \pi_0$.]
But as we saw in property P4, 
the multiplicative mean ergodic theorem 
(\ref{eq:mainMMET})
implies that $\LA''(a)>0$
for all $a$ for which it is 
valid, therefore it cannot
be valid for real $a\geq a^*$.

A plot of $ \LA(a)$ for 
$F=\ind_{0^c}-\pi(0^c)$ 
is shown in \Figure{mm1linear}.

\begin{figure}
\begin{center} 
\Ebox{.6}{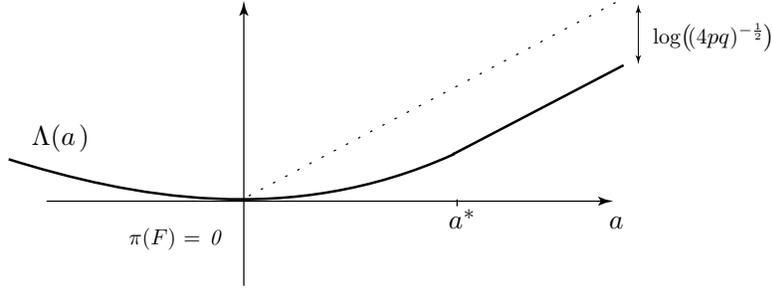}
\end{center}
\caption{The solid-curve shows the log-moment generating function 
$\LA(a),$ $a\in\RL$, for the $M/M/1$ queue with 
$F=\ind_{0^c}-\pi(0^c)$. 
It is  strictly convex for $a< a^*$,
and it is linear for $a\geq a^*$.}
\flabel{mm1linear}
\end{figure}

\subsection{Diffusions}

Consider an elliptic diffusion on a
manifold $\state$.  We assume that $\bfPhi$ is non-explosive, 
so that the sample paths are continuous on $[0,\infty)$ 
with probability one.  It is then strong Feller and 
$\psi$-irreducible, where $\psi$ is Lebesgue measure 
on $\state$ (see e.g. \cite{pinsky:book}),
and compact subsets of $\state$ are small.

Consider the special case where 
$\state=\Re^n$ and the diffusion term 
is constant,
\[
\clA  = h \cdot\nabla_x + \half \sigma^2 \Delta ,
\]
where $\Delta$ denotes the Laplacian.  
If  $\cf_a\in C^2$ solves the multiplicative
Poisson equation for some $a\in\RL$, 
we may consider the {\em twisted} process
$\bfPhi_a$, that is, the Markov process
with transition semigroup defined as before,
\[
\cP^t_a (x, dy) = \lambda_a^{-t} \cf_a^{-1}
(x) \haP^t_a (x, dy) \cf_a (y)\, ,
\qquad t>0\,.
\]
If \eq FK/ holds, then the generator $\cclA_a$ 
of $\bfPhi_a$ is given by
\begin{equation}
\begin{array}{rcl}
\cclA_a &=& ( h + \nabla_x\sigma^2\cF_\alpha)\cdot \nabla_x + \half \sigma^2 
\Delta^2,
\\[.2cm]
&=& \clA + \sigma^2\nabla_x\cF_\alpha\cdot \nabla_x \, ,
\end{array}
\elabel{twistDiff}
\end{equation}
where $\cF_\alpha = \log( \cf_a)$. 
Note that the twisted process has the same 
diffusion term as the original -- only the drift is affected 
by the twisting.  

\subparagraph{Reflected Brownian motion.}

Diffusions with reflection are currently a popular model in 
the operations-research area.  Consider for example a 
two-dimensional reflected Brownian motion
(RBM) $\bfPhi$ on $\state =\Re_+^2$,  with normal 
reflection on each boundary. We show below
that (when the drift is negative) $\bfPhi$
is geometrically ergodic. But it is not 
Doeblin recurrent, and it 
does not satisfy (mV3) for any 
$f$ with compact sublevel sets
(for the same reasons as in
the reflected random walk example above).

Within the interior of $\state$, the sample 
paths are identical to those of the affine 
stochastic differential equation (SDE) model,
\[
d\Phi_i = -\delta_i\, dt +  dW_i,\qquad i=1,2,
\]
where $\bfmW=(W_1,W_2)$ is a standard Brownian motion (BM)
on $\Re^2$, and the drift term $\delta_i$ is positive for each $i$.  
A characterization of the generator can
be found in \cite{wil85}. 

Suppose that $V\colon\Re^2\to \Re_+$ is smooth, and suppose that the following
boundary conditions are satisfied,
\begin{equation}
\frac{\partial}{\partial x_1} V \le  \frac{\partial}{\partial x_2}V 
        ,\quad x_1=0;
\qquad\qquad
\frac{\partial}{\partial x_1} V \ge  \frac{\partial}{\partial x_2}V 
        ,\quad x_2=0.
\elabel{rbmBoundary}
\end{equation}
Then, with $L =- \delta\cdot\nabla + \half \Delta$, the process below is
a supermartingale,
\[
m(t) = V(\Phi(t)) - \int_{[0,t)} L V\, (\Phi(s))\, ds.
\]

A candidate Lyapunov function is the quadratic,
$V_0(x)= \|x\|^2$, since
$LV_0 =- 2\delta\cdot x + 1$ is negative for large $x$, and  the 
boundary conditions \eq rbmBoundary/ are satisfied.
The supermartingale property implies that
\[
\Expect_x[V_0(\Phi(t))] \le - \delta_0 \int_{[0,t)}
        \Expect[\Phi_1(s)+\Phi_2(s)]\, ds + t\, ,
\]
where $\delta_0=2\min(\delta_1,\delta_2)$.  
This may be seen as a generalization of (V3), 
with $f(x)$ equal to a norm on $\Re^2$.
To obtain a version of (V4), first consider $V_1=\sqrt{V_0}$.  
We have, for some $\delta_1>0$, and some $B<\infty$,
\[
\Expect_x[V_1(\Phi(t))] \le - \delta_1 t,\qquad 0\le t\le 1,\ \|x\|\ge B.
\]
Finally, setting $V(x) = \exp(\beta V_1(x))$,
we can find $\beta>0$ 
sufficiently small
such that
\[
\Expect_x[V(\Phi(t))] \le \exp( -\beta t) V(x),\qquad 0\le t\le 1,\ \|x\|\ge B.
\]

We conclude that the RBM is geometrically ergodic, 
provided the reflection is normal and the drift is 
negative. Therefore, for any bounded functional $F$,
from Theorems~\ref{t:Bahadur-RaoNL} and~\ref{t:Bahadur-RaoL} 
we get precise large deviations bounds for the 
time-averages $\{S_t\}$, at least in some interval
around the mean $\pi(F)$ of $F$. Moreover, in view
recent results in \cite{budhiraja-dupuis:preprint}
(where a detailed study of large deviations 
properties of RBM (in one dimension) is performed),
we should not expect the limit 
theorems~\ref{t:Bahadur-RaoNL} 
and~\ref{t:Bahadur-RaoL} to hold 
on the whole real line. 

In general, geometric ergodicity depends upon the 
interaction of the drift vector and the reflection 
vectors along the boundaries.  In multidimensional 
models it is not always obvious how to choose an 
appropriate Lyapunov function, but one can devise 
numerical methods to search for a quadratic $V_0$ 
satisfying the required constraints; see 
\cite{kummey96a,schwerer97}.

\section*{Acknowledgments}
%%%%%%%%%%%%%%%%%%%%%%%%%%%%%%%%%%%%%%%%%%%%%%%%%%%%%%%%%%%%%%%%%%
Several interesting discussions with Amir Dembo 
on large deviations properties of Markov chains
are gratefully acknowledged.

\section*{Appendix}
%%%%%%%%%%%%%%%%%%%%%%%%%%%%%%%%%%%%%%%%%%%%%%%%%%%%%%%%%%%%%%%%%%

\subsection*{Proof of P1--P7}
P1, P2, P3 and the analyticity
of $\LA(\alpha)$ follow from the multiplicative
mean ergodic theorems in \Theorem{mainSpectral}
and \Theorem{mainSpectral2}.
 
To establish P4, first note that $\LA(0)=0$
follows from the uniform
convergence in \Theorem{mainSpectral}.
Similarly it follows that $\LA'(0)=\pi(F)=0$
and that $\LA''(0)=\lim_n(1/n)\VAR_x(S_n)=\sigma^2>0$
by assumption \eq assumption/ and \Proposition{clt}.

On considering the kernel $\cP_a$ for real  $a\in\Omega$,
since $\LA''(0)>0$, $\LA''(a)\geq 0$ for all such $a$, 
and $\pi_a(F)=\LA'(a)$ by \Proposition{V-geo2}, 
it follows that $\Lambda'(a)=\cpi_a(F)>0$ for all 
nonzero $a\in\Omega$. Now,
if $\LA''(a) = \sigma_a^2$
is zero, then by \eq fishEq4/ in
\Proposition{clt} it follows that $\pi(F-\cpi_a(F))=0$.
This is impossible since $\cpi_a(F)>0$.
 
The exponential convergence in P6 is given in
\Theorem{Vuniform}.
The analyticity of $\cf_\alpha$ is stated in
\Theorem{basic1}, and $\cf_\alpha|_{\alpha=0} \equiv 1$
by P1. \Proposition{V-geo2} combined with
\Proposition{basic1}
give P7.

Property P5 requires more work.
For  a neighborhood $\clO$ of zero the function $\cf_\alpha$ given below is
a constant times the normalized eigenfunction given in \eq cf-norm/:
\[
\cf_\alpha = H_\alpha^{-1} \One,
\quad
H_\alpha= I\lambda_\alpha -  \haP_\alpha + \One\otimes\pi,\qquad \alpha\in\clO.
\]
It is the unique solution in $L_\infty^V$ satisfying $\pi(\cf_\alpha)=1$.
Hence, for all $k$,
\[
\pi\left(\frac{d^k}{d\alpha^k} \cf_\alpha\right) = 0,\qquad \alpha\in\clO\,.
\]
 
We have a form of the quotient rule,
\[
\cf_\alpha' = - H_\alpha^{-1} H_\alpha' H_\alpha^{-1} \One,
\]
and after repeated differentiation we obtain
\[
\begin{array}{rcl}
\frac{d^3}{d\alpha^3} \cf_\alpha
&=&
-6 H_\alpha^{-1} H_\alpha'H_\alpha^{-1} H_\alpha'H_\alpha^{-1} H_\alpha'H_\alpha
^{-1} \One
\\[.15cm]
&&+3 H_\alpha^{-1} H_\alpha'H_\alpha^{-1} H_\alpha''H_\alpha^{-1}\One
\\[.15cm]
&&+3 H_\alpha^{-1} H_\alpha''H_\alpha^{-1} H_\alpha'H_\alpha^{-1}\One
\\[.15cm]
&& -H_\alpha^{-1} H_\alpha'''H_\alpha^{-1} \One.
\end{array}
\]
Evaluating at $\alpha=0$,
we have $H_0^{-1} = Z = [I-P+\Pi]^{-1}$ and
 $\frac{d^k}{d\alpha^k} H_\alpha\Big|_{\alpha=0}
= [I \lambda_0^{(k)} - (I_F)^k P]$,  $k\ge 1$.
Using $\Pi Z=\Pi$, and $Z\One = P\One=\One$ then gives,
\[
\begin{array}{rcl}
0=\pi\Bigl(\frac{d^3}{d\alpha^3} \cf_\alpha\Bigr)
&=&
6 \Pi I_F PZ I_F PZ F
\\[.15cm]
&&+3 \Pi I_F PZ (F^2 - \sigma^2)
\\[.15cm]
&&+3 \Pi(I_F^2-\sigma^2) P Z I_F
\\[.15cm]
&&+ \Pi(F^3-\lambda_0''').
\end{array}
\]
The proof is then complete on interpreting these formulae, 
since $\LA'''(0)=\lambda_0'''$, and
\[
PZG\, (x) = \pi(G) + \sum_{k=1}^\infty \Expect_x[G(\Phi(k)) -\pi(G)]
\,,
\]
for any function $G\in L_\infty^V$.
\qed
 
\subsection*{Proof of \Theorem{EdgeworthNL}}
We follow closely Feller's argument
in the proof of Theorem~1 in
\cite[p.~539]{fellerII:book},
leading to the statement
(\ref{eq:EdgeNL2}).
Choose and fix $x\in\state$
arbitrary. For $n\geq 1$, define
\be
M_n(\alpha):=\Expect_x[\exp(\alpha \oo{S}_n)]
        =m_n(\alpha)\exp(-\alpha\Expect_x\{S_n\})
\,,\qquad\alpha\in\Co,
\label{eq:Mn}
\ee
and the distribution functions
\be
\Psi_n(y):={\cal G}(y)-\frac{\rho_3}{6\sigma^3\sqrt{n}}(y^2-1)\gamma(y)
\,,\qquad y\in\RL,
\label{eq:psi-n}
\ee
with corresponding characteristic functions
\be
\phi_n(\omega):=
\exp(-\omega^2/2)\left(1+\frac
{\rho_3(i\omega)^3}{6\sigma^3\sqrt{n}}\right)\,,
\qquad\omega\in\RL.
\label{eq:phin}
\ee

Let $\epsilon>0$ arbitrary.
Choose $A$ large enough so that 
$A > 24(\epsilon\pi)^{-1}|\Psi'_n(y)|$ for
all $y\in\RL$, $n\geq 1$.
From Esseen's smoothing lemma given in 
\cite[p.~538]{fellerII:book}, with $T=A\sqrt{n}$
we get that, 
\be
|\barG_n(y)-\Psi_n(y)|\leq\frac{1}{\pi}\int_{-A\sqrt{n}}^{A\sqrt{n}}
\left|M_n\left(\frac{i\omega}{\sigma\sqrt{n}}\right)-
\phi_n(\omega)\right|\frac{d\omega}{|\omega|}\;
+\;\frac{\epsilon}{\sqrt{n}}\,,\qquad 
y\in\RL.
\label{eq:integral}
\ee
To prove (\ref{eq:EdgeNL2}) 
it suffices to show 
that this integral is $o(n^{-1/2})$.

We first consider the integral 
in the range $B\sqrt{n}\leq|\omega|\leq A\sqrt{n}$, 
with
$0<B<\min\{\sigma\baromega_0,A\}$ 
to be chosen later (where $\baromega_0$ is 
as in P7).  Applying the
change of variables $t=\omega/(\sigma\sqrt{n})$,
this integral is bounded above by
$$
\frac{\sigma}{B\pi}\int_{\frac{B}{\sigma}\leq|t|\leq\frac{A}{\sigma}}
        |m_n(it)|dt
+ 
\frac{\sigma}{B\pi}\int_{\frac{B}{\sigma}\leq|t|\leq\frac{A}{\sigma}}
|\phi_n(\sigma\sqrt{n}t)|dt.$$
The second integrand converges to zero
exponentially fast,
uniformly over $t$ in that range,
and the first integrand  
converges to zero
exponentially fast
by P2. Therefore, the
above expression is
certainly no larger
that $o(n^{-1/2})$.

Next we consider the integral in (\ref{eq:integral})
in the range $|\omega|\leq B\sqrt{n}$.
From the definition of $M_n$
and by properties P1 and P6, 
after the change of variables
$t=\omega/(\sigma\sqrt{n})$
this equals
$$
\frac{1}{\pi}\int_{|t|\leq \frac{B}{\sigma}}
\Big|
\exp\left(-it
\Delta_n-it\haF(x)\right)
\exp\left(n\LA(it)
\right)
[
\cf_{it} (x) + 
        it \epsilon_n]
-
\phi_n(t\sigma\sqrt{n})
\Big|\;
\frac{dt}{|t|}\, .
$$
Expanding $\LA(it)$ in a Taylor series
around zero yields
$$
\frac{1}{\pi}\int_{-\frac{B}{\sigma}}^{\frac{B}{\sigma}}
\exp(-\half n t^2\sigma^2)
\left|
\exp\left\{-it\Delta_n
+
\textstyle\frac{n}{6}(it)^3\LA'''(it)
        +\log( \cf_{it} + it \epsilon_n)
        -it\haF
\right\}
-
1-\frac{n\rho_3}{6} (it)^3
\right|\;
\frac{dt}{|t|}\,,$$
for some real $s=s(t)$ with
$|s|<B/\sigma$. Noting that,
$$
\log( \cf_{it} (x) + it \epsilon_n)
        -it\haF = \delta(it) + \log\left(1+
        \frac{it\epsilon_n}{\cf_{it}}
        \right)\,,
$$
where $\delta(\cdot)$ is as in P7,
the second exponent in the 
above integrand can be written as
$$
\textstyle\frac{n}{6}(it)^3\LA'''(it)
-it\Delta_n 
+\delta(it)
+it\epsilon_n''(it)
$$
where 
$\epsilon_n''(it):=[\log(1+
        \frac{it\epsilon_n}{\cf_{it}}
        )]/(it)$, and
\be
|\epsilon_n''(it)|
\to 0
\qquad\mbox{exponentially fast},
\qquad n\to\infty,
\label{eq:edp}
\ee
uniformly in $|t|\leq B/\sigma$ (by P1).  
Therefore, the integral we wish to bound is
\be
\frac{1}{\pi}\int_{-\frac{B}{\sigma}}^{\frac{B}{\sigma}}
\exp(-\half n t^2\sigma^2)
\left|
\exp\left\{
\textstyle\frac{n}{6}(it)^3\LA'''(is)
-it\Delta_n
+it\epsilon_n''(it)
+\delta(it)
\right\}
-
1-\frac{n\rho_3}{6} (it)^3
\right|\;
\frac{dt}{|t|}\,.
\label{eq:tobound}
\ee
To show that this is $o(n^{-1/2})$
we will apply the 
following simple inequality
from \cite[p.~534]{fellerII:book},
\be
|e^\alpha-1-\beta|\leq(|\alpha-\beta|+\half|\beta|^2)e^\gamma
\label{eq:abc}
\ee
where $\gamma\geq\max\{|\alpha|,|\beta|\}$.
First we choose $B$ small 
enough so that the following four bounds 
hold for all $|t|<B/\sigma$,
\ben
&(a)&\;\;
        |\LA'''(it)-\rho_3|<6\epsilon\;\\
&(b)&\;\;
        \frac{B}{6\sigma^3}|\LA'''(it)|\leq\frac{1}{4}\;\\
&(c)&\;\;
        |\delta(it)|\leq\frac{\epsilon}{2}\;\\
&(d)&\;\;
        \frac{B\rho_3}{6\sigma^3}\leq \frac{1}{4}\;,
\een
where $(a)$ and $(b)$ are possible by
the analyticity of $\LA(\cdot)$
and the definition of $\rho_3$ in P5, 
and $(c)$ is possible because of P7.
Then, writing
\ben
& & \alpha\,:=\, \frac{n}{6}(it)^3\LA'''(is)
-it\Delta_n
+it\epsilon_n''(it)
+\delta(it)
\\
&\mbox{and}&
\beta\,:=\,
\frac{n\rho_3}{6} (it)^3
\een
using $(b)$, $(c)$, and (\ref{eq:edp}),
we can bound
\ben
|\alpha|
&\leq& 
        nt^2\sigma^2 \frac{B}{6\sigma^3}
        \sup_{|t|<B/\sigma}|\LA'''(it)|
        +|t|[|\Delta_n| + |\epsilon_n''(it)|] +|\delta(it)| \\
&\leq& 
        \frac{1}{4}nt^2\sigma^2 + \frac{|t|}{n} + \epsilon/2\\
&\leq&
        \frac{1}{4}nt^2\sigma^2 + \epsilon\,,
\een
where the last two inequalities 
are valid after taking $n$ large enough.
Similarly, using $(a)$, $(b)$,
(\ref{eq:edp}), and P7,
\ben
|\alpha-\beta|
&\leq&
        \frac{n}{6}|t|^3|\LA'''(is)-\rho_3|+\frac{|t|}{n}+
        (\mbox{Const})t^2\\
&\leq&
        \epsilon n |t|^3 + \frac{|t|}{n}+
        (\mbox{Const})t^2\,,
\een
for $n$ large enough,
and using $(d)$,
$$|\beta|\leq \frac{1}{4}nt^2\sigma^2.$$

Applying inequality (\ref{eq:abc}) 
with $\gamma:=\frac{1}{4}nt^2\sigma^2$ and
in conjuction with the last three bounds,
the integral in (\ref{eq:tobound}) is
bounded above by
\ben
\frac{1}{\pi}\int_{-\frac{B}{\sigma}}^{\frac{B}{\sigma}}
\exp(-\onefourth n t^2\sigma^2 + \epsilon)
\left[
        \epsilon n |t|^3 
        + (\mbox{Const})t^2
        + \frac{|t|}{n}
        + \frac{1}{2}\left(\frac{n\rho_3|t|^3}{6}\right)
\right]\;
\frac{dt}{|t|}\,.
\een
and changing variables back to $\omega=t(\sigma\sqrt{n})$,
\ben
& & 
        \frac{e^\epsilon}{\pi}\int_{-B\sqrt{n}}^{B\sqrt{n}}
        e^{-\onefourth\omega^2}
        \left[
        \epsilon\left(\frac{\omega^2}{\sqrt{n}\sigma^3}\right)
        +(\mbox{Const})\frac{|\omega|}{n\sigma^2}
        +\frac{1}{n^{3/2}\sigma} 
        +\frac{\rho_3|\omega|^5}{72n\sigma^2}
        \right]
        \,d\omega\\
& & \; \leq 
        \frac{\epsilon}{\sqrt{n}}\left(\frac{e^\epsilon}{\pi}
        \int_{-\infty}^{\infty}\omega^2e^{-\onefourth\omega^2}
        \,d\omega\right)
        \;+\;O\!\left(\frac{1}{n}\right)\,.
\een
Since $\epsilon$ was arbitrary this shows that 
the integral in (\ref{eq:tobound}) is
$o(n^{-1/2})$,
and completes the proof.
\qed 

\subsection{Proof of \Theorem{EdgeworthL}}
We follow closely Feller's argument
in the proof of Theorem~2 in
\cite[p.~540]{fellerII:book}.
Choose and fix an
arbitrary $x\in \state$.
Let $\epsilon>0$ arbitrary,
and let $\Psi_n$ be the distribution
function (\ref{eq:psi-n}).
Recall that $G_n^\#=G_n*U[-h_n/2, h_n/2]$
and $\barG_n^\#=\barG_n*U[-h_n/2, h_n/2]$.
Proceeding as in \cite[p.~540]{fellerII:book}
along equations (4.9) and (4.10) 
(with $\Psi_n$ in place of ``$G$'' 
and $\Psi_n^\#:=\Psi_n * U[-h_n/2, h_n/2]$ 
in place of ``$G^\#$'')
we obtain, after taking $A>0$ large enough, 
\be
\left|\barG_n^\#(y)-\Psi_n(y)\right|
\leq\frac{1}{\pi}\int_{-A\sqrt{n}}^{A\sqrt{n}}
\left|M_n\left(\frac{i\omega}{\sigma\sqrt{n}}\right)-\phi_n(\omega)
        \right|\,\frac{|s_n(\omega)|}{|\omega|}\,d\omega
        \;+\;\frac{\epsilon}{\sqrt{n}}
        \;+\;O\!\left(\frac{1}{n}\right)
\,,\qquad
y\in\RL\,,
\label{eq:integralL}
\ee
where $M_n$ and $\phi_n$ are defined in
(\ref{eq:Mn}) and (\ref{eq:phin}),
and 
$$s_n(\omega):=\frac{\sin(\half h_n\omega)}{\half h_n\omega}\,,
        \qquad \omega\in\RL.$$
To prove (\ref{eq:EdgeL2})
it suffices to show
that the integral 
in (\ref{eq:integralL}) is $o(n^{-1/2})$.
We separately consider the integral 
over 
$|\omega|\leq B\sqrt{n}$
and over
$B\sqrt{n}\leq|\omega|\leq A\sqrt{n}$,
for some conveniently chosen 
$B<\min\{\sigma\baromega,A\}$,
where $\baromega$ is as in P7.
Noting that $|s_n(\omega)|\leq 1$ for
all $\omega$, the integral in the
former range can be shown to be
of order $o(n^{-1/2})$ 
as in the non-lattice case.

Therefore, it remains to show that 
\be
& & 
        \int_{B\sqrt{n}\leq|\omega|\leq A\sqrt{n}}
        \left|M_n\left(\frac{i\omega}{\sigma\sqrt{n}}\right)-\phi_n(\omega)
        \right|\,\frac{|s_n(\omega)|}{|\omega|}\,d\omega
        \nonumber\\
& & \leq
        \int_{B\sqrt{n}\leq|\omega|\leq A\sqrt{n}}
        \left|m_n\left(\frac{i\omega}{\sigma\sqrt{n}}\right)\right|
        \frac{|s_n(\omega)|}{|\omega|}\,d\omega
        +
        \int_{B\sqrt{n}\leq|\omega|\leq A\sqrt{n}}
        \frac{|\phi_n(\omega)|}{|\omega|}\,d\omega
        \;\;=\;\;o(n^{-1/2})\,.
        \label{eq:remainder}
\ee
The last integral above is easily seen 
to decay exponentially in $n$
(as in the proof of \Theorem{EdgeworthNL}), and hence we
concentrate on the former integral, which, 
after the change of variables
$t=\omega/(\sigma\sqrt{n})$,
becomes
$$\frac{2}{h}\int_{\frac{B}{\sigma}\leq |t|\leq\frac{A}{\sigma}}
|m_n(it)|\,|\sin(\half ht)|\,\frac{dt}{t^2}.$$

Notice that $|\sin(\half ht)|$ and $m_n(it)$
are periodic functions
of $t$ with period $2\pi/h$.
Consider, without loss of generality, the range
of $t\in[B/\sigma,A/\sigma]$ (the case of negative
$t$ is similar). Let $(k-1)$ denote the number of 
full periods of length $2\pi/h$ in that interval.
Then, since $|\sin(y)/y|\leq 1$ for all real $y$,
$$
\frac{2}{h}\int_{\frac{B}{\sigma}\leq t \leq\frac{A}{\sigma}}
|m_n(it)|\,|\sin(\half ht)|\,\frac{dt}{t^2}
\leq
        \frac{k\sigma}{B}\int_{\frac{B}{\sigma}}
        ^{\frac{2\pi}{h}-\frac{B}{\sigma}}
        |m_n(it)|dt
\;+\;
        \frac{k\sigma^2}{B^2}
        \int_{-\frac{B}{\sigma}}^{\frac{B}{\sigma}}
        |m_n(it)|\,|t|\,dt\,,
$$
where the first integral converges to zero 
exponentially fast by P3. 
Using P1 to bound $m_n(it)$ and
expanding $\LA(it)$ in a Taylor
series, 
the second integral is
$$
        C\int_{-\frac{B}{\sigma}}^{\frac{B}{\sigma}}
        |t|\,|\exp\{n\LA(it)\}|
        \,dt
\leq
        C'
        \int_{-\infty}^{\infty}
        |t|\,
        |\exp\{-\half nt^2\sigma^2\}|
        \,dt
\,=\,
        \frac{C''}{n}\,,
$$
for some constants $C$, $C'$ and $C''$.
This establishes (\ref{eq:remainder}) 
and completes the proof.
\qed

\subsection{Proof of \Lemma{br-lemma}}
Part~(ii) is immediate by the choice
of $c$ and property P4. For part~(i)
note that
by the uniform convergence of $\LA_n(a)$
to $\LA(a)$ (property P1) we also have
convergence of their derivatives, so
for $n$ large enough we can pick $a_n$
as claimed, and since $\LA''_n(a)$
eventually will be strictly positive 
for all $a\in(-\bara,\bara)$, this
$a_n$ is unique.

For part~(iii) recall that 
$\LA'_n(a_n)=c=\LA'(a)$, so the 
fact that $a_n\to a$ as $n\to\infty$
follows by
the uniform convergence of the
functions $\LA'_n$. Moreover,
expanding $\LA_n'(a)$ around $a=a_n$ 
and using P1,
\ben
0
&=& \LA'(a)-\LA'_n(a_n)\\
&=& [\LA_n'(a)-\LA'_n(a_n)] - [\LA'(a)-\LA'_n(a)]\\
&=& [(a-a_n)\LA_n''(a_n)+O(a-a_n)^2] + \frac{d}{da}\left[
        \frac{1}{n}\log(\cf_a(x)+a\epsilon_n\exp\{-n\LA(a)\})
        \right]\\
&=& [(a-a_n)\LA_n''(a_n)+O(a-a_n)^2] 
        + O\!\left(\frac{1}{n}\right)
        + O(\epsilon_n),
\een
where in the last step we used P7.
Taking $n$ large enough so that $\{\LA_n''(a_n)\}$ 
is a bounded sequence, 
bounded away from zero from below, 
this implies that $(a_n-a)=O(1/n)$.

Part~(iv) is an immediate consequence of 
P4 and of the uniform convergence in P1.
Finally for~(v) we have from~(i),~(ii), and P1,
\ben
\LA_n^*(c) 
&=& a_n c -\LA(a_n) -  \frac{1}{n}\log\left[
        \cf_a(x)+a\epsilon_n
        \right] \\
&=& \LA^*(c) + (a_n-a)c + (\LA(a)-\LA(a_n)) 
        -\frac{1}{n}\log \cf_{a_n}
        -\frac{1}{n}\log(1+a_n\epsilon_n/\cf_{a_n})\,,
\een
and, using~(iii) and P7,
\ben
\LA_n^*(c) 
&=& \LA^*(c) 
+ O\!\left(\frac{1}{n^2}\right)
- \frac{1}{n}\log \cf_a
+ O\!\left(\frac{1}{n^2}\right)
+ O(\epsilon_n)\,,
\een
as required.
\qed

\newpage

% \bibliography{../../../latex/ik}

\end{document}